\documentclass[11pt,leqno,twoside]{amsart}

  \usepackage[utf8]{inputenc} 
  \usepackage[T1]{fontenc} 
  \usepackage{lmodern}

\usepackage[scale=0.8]{geometry}

\usepackage{cite}
\usepackage{morefloats}
\usepackage{mathtools} %For the command "coloneq", used for the "definition symbol"
\usepackage{enumitem}
\usepackage[english]{babel}

\usepackage{caption}
\usepackage{subcaption}
\usepackage{svg}

\numberwithin{equation}{section}

\usepackage{tikz}
\usepackage{multirow}

\usepackage{algorithmic}
\usepackage{algorithm}

\usepackage{todonotes}

\usepackage{amsmath, amsthm}
\usepackage{amsfonts, amssymb}
\usepackage{bm}               % Lettres grecques en gras
\usepackage{hyperref}

	\newtheorem*{theorem*}{Theorem}
	\newtheorem{lemma}{Lemma}[section]
	\newtheorem{proposition}{Proposition}[section]
	\newtheorem*{proposition*}{Proposition}

\theoremstyle{remark}
\newtheorem{remark}{\textbf{Remark}}[]

\theoremstyle{definition}

\DeclareMathOperator{\sinc}{sinc}

\DeclareMathOperator{\vect}{Span}

\DeclareMathOperator{\tr}{Trace}

\DeclareMathOperator{\diverg}{div}

\newcommand{\ve}{\varepsilon}

\newcommand{\ie}{\emph{i.e. }}
\newcommand{\RR}{\mathbb{R}}

\newcommand{\ZZ}{\mathbb{Z}}

\newcommand{\SSS}{\mathbb{S}}
\newcommand{\Q}{\mathcal{Q}}
\newcommand{\G}{\mathcal{G}}
\newcommand{\Tau}{\mathcal{T}}
\newcommand{\LL}{\mathcal{L}}

\title[Hierarchical domain decomposition of the Boltzmann equation]{Hierarchical dynamic domain decomposition for the multiscale Boltzmann equation}

\author[D. Caparello]{Domenico Caparello}
\address{Domenico Caparello \\
Université Côte d’Azur, CNRS, LJAD, Parc Valrose, F-06108 Nice, France \\
Department of Mathematics and Computer Science University of Ferrara
}
\email{Domenico.CAPARELLO@univ-cotedazur.fr}
\email{domenico.caparello@unife.it}

\author[L. Pareschi]{Lorenzo Pareschi}
\address{Lorenzo Pareschi \\
Maxwell Institute for Mathematical Sciences and Department of Mathematics
Heriot-Watt University, Edinburgh \\
Department of Mathematics and Computer Science University of Ferrara \\
}
\email{L.Pareschi@hw.ac.uk}

\author[T. Rey]{Thomas Rey}
\address{Thomas Rey \\
Université Côte d’Azur, CNRS, LJAD, Parc Valrose, F-06108 Nice, France
}
\email{Thomas.REY@univ-cotedazur.fr}

\keywords{Boltzmann equation, compressible Euler, ES-BGK operator, Burnett transport coefficients, asymptotic-preserving methods, domain decomposition method, hierarchical numerical method, spectral methods}
	
\subjclass[2010]{Primary: 76P05, % Rarefied gas flows; Boltzmann equation
  82C40, % Time-dependent statistical mechanics; Kinetic theory of gases
  Secondary: 65N08, % Numerical analysis; Finite volume methods
  65N35 % Numerical analysis; Spectral, collocation and related methods 
}
\begin{document}
  
	\begin{abstract}
    In this work, we present a hierarchical domain decomposition method for the multi-scale Boltzmann equation based on moment realizability matrices, a concept introduced by Levermore, Morokoff, and Nadiga \cite{lev-mor-nad-1998}. This criterion is used to dynamically partition the two-dimensional spatial domain into three regimes: the Euler regime, an intermediate kinetic regime governed by the ES-BGK model, and the full Boltzmann regime. The key advantage of this approach lies in the use of Euler equations in regions where the flow is near hydrodynamic equilibrium, the ES-BGK model in moderately non-equilibrium regions where a fluid description is insufficient but full kinetic resolution is not yet necessary, and the full Boltzmann solver where strong non-equilibrium effects dominate, such as near shocks and boundary layers. This allows for both high accuracy and significant computational savings, as the Euler solver and the ES-BGK models are considerably cheaper than the full kinetic Boltzmann model.
    To ensure accurate and efficient coupling between regimes, we employ asymptotic-preserving (AP) numerical schemes and fast spectral solvers for evaluating the Boltzmann collision operator. Among the main novelties of this work are the use of a full 2D spatial and 3D velocity decomposition, the integration of three distinct physical regimes within a unified solver framework, and a parallelized implementation exploiting CPU multithreading. This combination enables robust and scalable simulation of multiscale kinetic flows with complex geometries.	
	\end{abstract}
  \maketitle

  \tableofcontents

  \section{Introduction}\label{sec1}
    \setcounter{equation}{0}
  Many problems in engineering and physics involve fluid flows in transitional regimes, such as atmospheric reentry of spacecraft, micro-electro-mechanical systems (MEMS), and industrial heat transfer systems. In these settings, classical macroscopic models like the Euler equations are often inadequate due to the presence of shocks, kinetic boundary layers, or other non-equilibrium effects. Accurately capturing such dynamics requires kinetic models that resolve the microscopic behavior of the gas.

Kinetic theory is a cornerstone of modern applied mathematics and physics. Its roots trace back to foundational thermodynamic works by Bernoulli, Avogadro, and Carnot in the 18th and 19th centuries. One of the earliest mathematical contributions was by James Clerk Maxwell, who, in 1859, introduced a statistical model for particle interactions in his study of Saturn's rings \cite{Maxwell:1859}. This was followed by Ludwig Boltzmann’s celebrated formulation of the Boltzmann equation \cite{Boltzmann:1872}, which remains central to the modeling of rarefied gas dynamics.

Kinetic equations describe the evolution of a particle distribution function in a high-dimensional phase space, typically six-dimensional for a three-dimensional configuration domain. These models combine linear transport with stiff, nonlinear collision operators, making them significantly more computationally expensive than their fluid-dynamical counterparts. This ``curse of dimensionality" motivates the use of macroscopic equations such as Euler or Navier–Stokes models in regions where the flow is close to equilibrium, while reserving kinetic solvers for regions where non-equilibrium effects are dominant.

This leads to the concept of hybrid multiscale models and hybrid numerical schemes, which have been widely developed in the literature. These methods blend kinetic and fluid models, often using domain decomposition and switching criteria based on physical indicators. Relevant examples include domain decomposition methods \cite{DiPa08, tiw-1998,dim-mie-ris-2014,Klar98}, asymptotic-preserving (AP) schemes \cite{fil-jin-2010,fil-jin-2011,jin-1999,jin-par-2001,DiPa13}, and Boltzmann solvers with particle/DSMC coupling \cite{tiw-kla-har-2009,bur-boy-2009,dim-par-2014,deg-dim-2012,DiPa07,DiPa08,DiPa10,CaDi08}. Our approach is inspired by these developments and by the foundational works of Levermore, Morokoff, and Nadiga \cite{lev-mor-nad-1998} and Filbet and Rey \cite{fil-rey-2014}, but differs by introducing a three-regime decomposition strategy implemented in a fully adaptive, non-intrusive fashion.

The proposed hierarchical method offers two main advantages:
\begin{itemize}
    \item \textbf{High accuracy in multiscale regimes compared to classical fluid solvers:} The method dynamically applies the appropriate physical model in each region—Euler equations in equilibrium zones, ES-BGK in transitional regimes, and the full Boltzmann equation in highly non-equilibrium areas—ensuring accurate resolution of shocks, boundary layers, and other kinetic effects where fluid models fail.
    
    \item \textbf{Reduced computational cost compared to full Boltzmann solvers:} By confining the expensive kinetic computations to a minimal subset of the domain and relying on efficient fluid solvers elsewhere, the method achieves significant speedups over fully kinetic simulations without compromising physical fidelity.
\end{itemize}

The ES-BGK model plays a crucial role in this hierarchy: it allows the method to handle moderately non-equilibrium regions where the fluid description is insufficient, yet applying the full Boltzmann model would be unnecessarily costly. This three-layer decomposition offers a computational advantage over classical two-regime methods (e.g., Euler–Boltzmann \cite{fil-rey-2014}) while preserving greater physical accuracy.

To maximize computational efficiency, we minimize the size of the full Boltzmann region during simulations. The algorithm is designed to be non-intrusive with respect to the solvers, meaning that it can work with different implementations of both kinetic and fluid models.

It is well known that solving the Boltzmann equation in stiff regimes is especially challenging. In these regimes, frequent particle collisions rapidly drive the system toward local thermodynamic equilibrium, allowing for a reduced description via macroscopic or diffusive equations \cite{dim-par-2014}. However, in transitional regimes, collisions are sufficiently strong to introduce stiffness but not strong enough to restore equilibrium, leading to severe time step restrictions for explicit schemes.

Asymptotic-preserving (AP) schemes offer a solution: these numerical methods remain stable and accurate across a range of Knudsen numbers, including the stiff limit, without the need for mesh or time step refinement \cite{cor-per-1991, jin-1995, caf-shi-rus-1995, jin-par-tos-1998, jin-1999, jin-par-2001, klar-1998, gos-tos-2002}. Our method incorporates AP strategies to ensure consistent behavior across regimes and uses fast spectral solvers for the Boltzmann operator \cite{MoPa:2006,fil-mou-par-2006} to improve performance.

In this work, we develop and implement this hierarchical  strategy in a fully two-dimensional spatial domain with three-dimensional velocity dependence. We describe in detail the moment realizability-based switching criterion, the structure of the hybrid solver, and its coupling across regimes. The method is tested on benchmark problems featuring strong shocks and boundary layers, demonstrating both improved accuracy and significant computational savings. The paper is organized as follows: in Section~\ref{secPrelim}, we introduce the kinetic and fluid models; in Section~\ref{secRegimeIndic}, we describe the hybrid decomposition indicators whereas in Section~\ref{secNumSim} we describe the details of the numerical schemes adopted; in Section~\ref{subNumSim}, we present numerical experiments; and in Section~\ref{secCon}, we conclude with a discussion of future directions.

\section{The Boltzmann equation}
    \label{secPrelim}
        Let us now present the general kinetic model we are interested in. For a given nonnegative initial condition $f_0$, we will study a particle distribution function $f^\ve = f^\ve(t,x,v)$, for $t \geq 0$, $x \in \Omega \subset \mathbb{R}^{d_x}$ and $v \in \mathbb{R}^{3}$, solution to the initial-boundary value problem
		\begin{equation} 
			\label{eqCollision}
			\left\{ \begin{aligned}
			  & \frac{\partial f^\ve}{\partial t} + v \cdot \nabla_x f^\ve \,=\, \frac{1}{\ve}\;\Q(f^\ve), 
			  \\
		  	&  \;
			  \\
			  & f^\ve(0, x, v) = f_{0}(x,v),
			\end{aligned} \right.
		\end{equation}
		where the collision operator $\Q$ is a  Boltzmann-like operator, and $\ve > 0$ is the Knudsen number (it is equal to the ratio between the mean free path of particles and the physical length scale of the system). This quantity gives an indication of the rarefaction of the gas: if $\ve \ll 1$ it means that the particle system is subject to many collisions and then the gas is \emph{dense} and in fluid regime; vice-versa, if $\ve \sim 1$ it means that the number of collisions is small and then the gas can be considered \emph{rarefied} \cite{fil-rey-2014}.\\
		The open set $\Omega$ is a bounded Lipschitz-continuous domain of $\RR^{d_x}$, which means  that the model \eqref{eqCollision} has to be supplemented with boundary conditions described later.
		
		We assume that the collision operator fulfils the three following assumptions
		\begin{enumerate}[label=\textbf{(H$\bm{_\arabic{*}}$)}, ref=\textbf{(H$\bm{_\arabic{*}}$)}]
		
      \item \label{hypConservations} 
        It preserves mass, momentum and kinetic energy. This properties can be, respectively, expressed as
			  \begin{equation*}
			    \int_{\RR^{3}} \Q(f)(v) \, dv = 0, \quad  \int_{\RR^{3}} \Q(f)(v) \, v \, dv = 0, \quad \int_{\RR^{3}} \Q(f)(v) \, |v|^2 \, dv = 0;
			  \end{equation*}
			 
			\item \label{hypEntropy}
			  It dissipates the Boltzmann entropy (H-theorem)
			  \begin{equation*}
					\int_{\RR^{3}} \Q(f)(v) \, \log(f)(v) \, dv \, \leq \, 0;
				\end{equation*}
				
	    \item \label{hypEquilib} 
	      Its equilibria are given by Maxwellian distributions
	      \begin{equation*}
	        \Q(f) \, = \, 0 \quad \Leftrightarrow \quad f = \mathcal M_{\rho, \bm u, T} := \frac{\rho}{(2 \pi T)^{3/2}} \exp \left ( - \frac{|v-\bm u|^2}{2 T} \right ),
	      \end{equation*}
	      where the \emph{density}, \emph{velocity} and	\emph{temperature} of the gas $\rho$, $\bm u$ and $T$  are computed from the distribution function $f$ as  
			\begin{equation*}
	      \rho = \int_{v\in{\RR}^{3}}f(v)\,dv, \quad \bm{u} = \frac{1}{\rho}\int_{v\in{\RR}^{3}}v f(v) \, dv, 
	      \quad T = \frac{1}{3\rho} \int_{v\in{\RR}^{3}}\vert \bm{u} - v \vert^2 f(v) \,dv.
	    \end{equation*}
	      
	  \end{enumerate}  
		According to assumptions \ref{hypEntropy}-\ref{hypEquilib}, it can be proved that, when $\ve \to 0$, the distribution $f^\ve$ converges (at least formally) to a Maxwellian distribution, whose moments are solution to the compressible Euler system
			\begin{equation}
				\label{eqHydroClosedEuler}
				\left\{ \begin{aligned}
				  & \partial_t \rho  + \diverg_x  (\rho \, \bm{u} ) \,=\, 0, 
				  \\
				  \,
				  \\
				  & \partial_t(\rho \, \bm{u} ) + \diverg_x  \left(\rho \, \bm{u} \otimes \bm{u} \,+\, \rho \, T  \,{\rm\bf I}\right) \, =\, \bm{0}_{\RR^{3}}, 
				  \\
				  \,
				  \\
				  & \partial_t E + \diverg_x \left ( \bm{u} \left ( E +\rho \, T\right )  \right ) \,=\, 0.
				\end{aligned} \right.
			\end{equation}
        This limit is typically referred to as an \emph{hydrodynamic limit}.

    In the following we are going to present two operators that satisfy these hypothesis, and that will be implemented in our hybrid domain decomposition scheme.
    
    \subsection{The Boltzmann operator}
	    \label{subBoltzOp}

        If we consider only binary elastic collisions between particles, then a dilute gas can be described by the kinetic equation with the Boltzmann operator, and in a velocity space of dimension $d$ the collision operator is

			\begin{align}
			  \label{eqQBoltz}
			  \Q (f,f)(v) & \,=\, \int_{\RR^{d}} \int_{\mathbb{S}^{d-1}}  B(|v-v_*|,\cos \theta) \, \left[ f'_* f' \,-\, f_* f \right] \, d\sigma \, dv_*,
				\end{align}
			where we used the notation $f = f(v)$, $f_* = f(v_*)$, $f ^{'} = f(v')$, $f_* ^{'} = f(v_* ^{'})$.   
                The pair $(v,v_*)$ are the pre-collisional velocities, while the quantities $(v',v'_*)$ are the post-collisional velocities. These two pairs are related by the following relations
			\begin{equation*}
			  v' \,=\, \frac{v+v_*}{2} \,+\, \frac{|v-v_*|}{2} \,\sigma, \qquad v'_* \,=\, \frac{v+v^*}{2} \,-\, \frac{|v-v_*|}{2} \,\sigma.
			\end{equation*}
			The \emph{collision kernel} $B$ is a non-negative function which by physical arguments of invariance only depends on $|v-v_*|$ and	$\cos \theta = {\widehat u} \cdot \sigma$, where ${\widehat u} =(v-v_*)/|v-v_*|$. 
									
			It can be proved that the Boltzmann collision operator has the fundamental properties \ref{hypConservations}, \ref{hypEntropy} and \ref{hypEquilib}.

            In the remainder of this paper we shall use the assumption of \emph{short range interactions}, and in particular assume that $B$ is locally integrable. This assumption is satisfied by the \emph{hard-spheres model}, which in dimension $d=3$ can be written as
            \begin{equation}\label{eqVHS}
                B(|v-v_*|, \cos{(\theta)})=|v-v_*|.
            \end{equation}
            As an important benchmark model for the numerical simulation, and for the development of fast spectral methods for the computation of the Boltzmann operator, we therefore introduce the so-called \emph{variable hard-spheres model} (VHS), which writes
            \begin{equation}
                B(|v-v_*|, \cos{(\theta)})=C_\gamma |v-v_*|^\gamma ,
            \end{equation}
            for some $\gamma \in \left[0,1\right]$ and a constant $C_\gamma >0$.
            
            It is important to remark that, for this last class of models, it is possible to split the collision operator as follows
            \begin{equation}
                \Q(f,f) = \Q^+(f,f) - L(f)f,
            \end{equation}
            where 
            \begin{equation}
                \Q^+(f,f) = \int_{R^d}\int_{S^{d-1}}B(|v-v_*|, \cos{(\theta)})f' f'_* d\sigma dv_*,
            \end{equation}
            and
            \begin{equation}
                L(f) = \int_{R^d}\int_{S^{d-1}}B(|v-v_*|, \cos{(\theta)})f_* d\sigma dv_*.
            \end{equation}        
            The computation of Boltzmann-like collision operator, performed using deterministic numerical methods, requires to work on a bounded velocity space. In most cases, the classical technique is to periodize the particle distribution function and the collision operator (which correspond to add some non physical binary collisions), and then using Fourier series to compute the truncate operator. Unfortunately this causes the loss of some local invariants, but the mass can be preserved by implementing an accurate periodization, which is the basis of spectral methods \cite{rey-2023}, that we discuss in further details in Section \ref{spectralMethods}.

		\subsection{The ES-BGK operator}
		  \label{subBGK}
            The BGK operator \cite{bha-gro-kro-1954}, and its Ellipsoidal Statistical (ES-BGK) extension \cite{and-let-per-per-2000}, are two other well known collision operators which have the properties \ref{hypConservations}-\ref{hypEntropy}-\ref{hypEquilib}. 
		  The bilinear Boltzmann collision operator $\Q$ is replaced by a nonlinear relaxation operator, which matches the same hydrodynamic limit than the Boltzmann operator.

      First of all we introduce some macroscopic quantities of the particle distribution function $f$ such as the \emph{stress tensor}
			\begin{equation*}
			  \rho_f\,\Theta_f(t,x)\, =\,  \int_{{\RR}^{3}}  (v-\bm{u}_f)\otimes (v-\bm{u}_f)\,f(t,x,v) \,dv.
			\end{equation*} 
			Therefore the \emph{translational temperature} is related to the stress tensor as $T_f = \tr(\Theta_f)/3$. We finally introduce the corrected tensor
			\begin{equation*}
			  \Tau_f(t,x) \,\,=\,\, \left[(1-\beta) \, T_f \,{\rm\bf I} \,\,+\,\,\beta \,\Theta_f\right](t,x),  
			\end{equation*} 
			which can be viewed as a linear combination of the initial stress tensor $\Theta_f$ and of the isotropic stress tensor $T_f \,{\rm \bf I}$ developed by a Maxwellian distribution. 
			The parameter $-\infty< \beta < 1$ is used to modify the value of the Prandtl number through the formula 
			\[
			0 \,\leq\, {\rm Pr}  \,=\, \frac{1}{1-\beta} \,\leq\, +\infty\quad {\rm for } \quad \beta\in (-\infty\,,\,1).
			\]
			The correct Prandtl number for a monatomic gas of hard-spheres is equal to $2/3$ \cite{str-2005}, namely obtained here for $\beta = -1/2$ whereas the classical BGK operator, obtained for $\beta = 0$, has a Prandtl number equal to $1$.
			
			To define the ES-BGK operator, we introduce a corrected Gaussian $\G[f]$ defined by
			\[
			\G[f]= \frac{\rho_f}{\sqrt{{\rm det}(2\pi\,\Tau_f)}}\,\exp\left(-\frac{(v-\bm{u}_f)\,\Tau_f^{-1}\,(v-\bm{u}_f)}{2}\right),
			\]
			and the corresponding collision operator is now
			\begin{equation} 
			  \label{eqOpES-BGK}
			  \Q_\mathcal{ES-BGK}(f) =  \nu(\rho_f, T_f) \left( \G[f] \,-\, f \right),
			\end{equation}
			where $\nu$ is the collision frequency from the Boltzmann operator. It can be shown \cite{str-2005} that it depends only on the kinetic density $\rho_f$ and temperature $ T_f $.
            %In our numerical scheme we will use the ES-BGK operator, and not the BGK operator, since the former has better properties than the latter. %(INSERIRECITAZIONE E SPIEGARE QUALI SONO LE PROPRIETA.....PRANDTL NUMBER E ASIMMETRIA)

    \section{Regime indicators}
    \label{secRegimeIndic}
    Most of the papers in literature about decomposition methods rely on the domain decomposition introduced by Boyd, Chen and Chandler in \cite{boy-che-can-1995}, in which a macroscopic criterion is used to switch between the hydrodynamic description (easy to compute numerically, but inaccurate near shocks or boundary layers) and the kinetic one (computationally expansive but accurate in most of the situations).
    
    This criterion is based on the value of the local Knudsen number: if it is below a (problem-dependent) threshold, then the kinetic description is used. In fact, it is important to remark that the Knudsen number is defined as the ratio between the mean free path of the particles and the characteristic physical length scale of the system: so when the Knudsen number is very small, then it means that the gas is rarefied and so the usage of kinetic models is necessary.
    
    Several authors have implemented methods of this kind, considering also the spatial derivatives in the numerical solution itself, among which we recall Kolobov, \emph{et al.} in \cite{kol-ars-ari-fro-zab-2007}, and Degond and Dimarco in \cite{deg-dim-2012}.
    In literature it is possible to find other works where different, but similar, indicators are used. For example Tiwari in \cite{tiw-1998} introduced a criterion based on the viscous and heat fluxes of the Navier-Stokes equation, through a Grad's 13-moments expansion.
    Similar strategies have been used, among the others, by Tiwari, Klar and Hardit in \cite{tiw-kla-har-2009}, and by Alaia and Puppo in \cite{ala-pup-2011, ala-pup-2012}.

    % To our opinion, and as it was demonstrated in \cite{fil-rey-2014}, a possible drawback of some criteria in literature is that they are based on the macroscopic description, so it is possible that the indicator itself can be wrong in such situations where the fluid is far from the thermal equilibrium. We now want to give an example of a situation of this kind, by recalling \cite{fil-rey-2014}, where it is considered the case of a gas which is distributed in the velocity space as a sum of two Gaussian with non-zero mean, and which is constant in space
    % \begin{equation}
    %   \label{eqSumMaxwell}
    %     f(x,v) = \frac{1}{2}\left ( \mathcal{M}_{1, \bm{u}_0, 1}(v) + \mathcal{M}_{1, -\bm{u}_0, 1}(v) \right ), \quad \forall x \in \TT, \, v \in \RR^3,
    % \end{equation}
    % for \eg $ \bm{u}_0 = (1,0,1)$. This distribution is very far from the \emph{thermal equilibrium} given by Maxwellian distributions according to assumption \ref{hypEquilib}. Nevertheless, both criteria from \cite{boy-che-can-1995} or \cite{tiw-1998} would detect an hydrodynamic setting, which is not the case.     

    In this work we will use two criteria: one to detect when the hydrodynamic description is not valid any more (and then it is necessary to switch to the kinetic one), and the other one to understand when to switch from the kinetic regime to the hydrodynamic one.

    The idea of using hybrid decomposition solver is not limited only to the Boltzmann equation, but other recent works exploiting the lower dimensionality of the asymptotic fluid model for different kinetic equations, like the Vlasov-BGK or the Vlasov-Poisson-BGK equations in the diffusive scaling are \cite{lai-2023, lai-rey-2023, dim-mie-ris-2014}.  This shows the growing interest and need for the development of hybrid domain decomposition schemes to speed up the calculation of the solution for kinetic equations.

        In the following, we present the criteria, and we start introducing a key mathematical tool: the Chapmann-Enskog expansion.
	
        \subsection{The Chapmann-Enskog expansion}			
        \label{subChapmannEnskog}
	 	
        The criterion used in this article was introduced by Levermore, Morokoff and Nadiga in \cite{lev-mor-nad-1998}, and it has already been used by other authors, like Tiwari in \cite{tiw-2000}, Filbet and Rey in \cite{fil-rey-2014}, Filbet and Xiong in \cite{fil-xio-2018}, Li, Song and Wang in \cite{li-son-wan-2021} and Xiong and Qiu in \cite{xio-qiu-2017}. In particular, in this last paper, Xiong and Qiu propose a hierarchical numerical scheme based on the following three levels (Euler equations, Navier-Stokes equation and kinetic equation). % In this work, we shall use it with a three layers numerical scheme (Euler equation, Boltzmann equation with ES-BGK operator and Boltzmann equation with Boltzmann operator), on a 2D physical grid and a 3D velocity grid.\\
        This approach has the main interest to depend on the closure made for obtaining the hydrodynamic model. For the hydrodynamic regions it is possible to consider either the  compressible Euler equations, or the compressible Navier-Stokes, or other models, such as Burnett or Super-Burnett \cite{str-2005}. In our code we implemented the Euler equations.

%%%%%%%%%%%%%%%%%%%%%%%%
        Assuming that the distribution $f^\ve$ is close to equilibrium thanks to the relaxation property \ref{hypEntropy}, we can do formally the \emph{Chapman-Enskog} expansion
        \begin{equation}
          \label{devChapEnsk}
			  f^\ve \, = \, \mathcal M_{\rho^\ve,\bm{u}^\ve,T^\ve}\left [1 + \ve \, g^{(1)} + \ve^2 \, g^{(2)} + \ldots \right ],
        \end{equation}
        where the fluctuations  $g^{(i)}$ for $i \geq 0$ designate a function that depends smoothly on  \[
        (\rho^\ve,\rho^\ve\bm{u}^\ve,E^\ve) \, = \, \int_{\RR^{3}} f^\ve(v) \, \left(1, v, \frac{|v|^2}{2}\right) \, dv, \quad T^\ve \, = \, \frac{1}{3 \, \rho^\ve }\int_{\RR^{3}} f^\ve(v) \, |v-\bm{u}^\ve|^2 \, dv,
        \]   and any finite number of its derivatives with respect to the $x$-variable at the same point $(t, x)$, and on the $v$-variable. According to \ref{hypConservations}, it verifies
        \[
            \int_{\RR^{3}} g^{(i)}(v) \, \left(1, v, \frac{|v|^2}{2}\right) \, dv \, = \, \bm{0}_{\RR^{5}}^\intercal.
        \]
        Without any closure, according to the conservative properties \ref{hypConservations} of the collision operator $\Q$, we have 
        \begin{equation}
            \label{eqHydroNonClosed}
            \left\{ \begin{aligned}
            & \partial_t \rho^\ve  + \diverg_x (\rho^\ve \, \bm{u}^\ve ) \,=\, 0, 
            \\
            & \,
            \\
            & \partial_t(\rho^\ve \, \bm{u}^\ve ) + \diverg_x \left (\rho^\ve \bm{u}^\ve \otimes \bm{u}^\ve + \rho^\ve T^\ve \left ( \bm{I} + \bm{\bar{A}}^\ve\right )\right ) \,=\, \bm{0}_{\RR^{3}}, 
            \\
            & \,
            \\
            & \partial_t E^\ve + \diverg_x \left (\frac 12 \rho^\ve |\bm{u}^\ve|^2 \bm{u}^\ve + \rho^\ve T^\ve \left (  \frac{5}{2}\bm{I} + \bm{\bar{A}}^\ve \right ) \bm{u}^\ve + \rho^\ve (T^\ve)^{3/2} \bm{\bar{B}}^\ve \right ) \,=\, 0,
            \end{aligned} \right.
        \end{equation}
        where the traceless matrix $\bm{\bar{A}}^\ve \in M_{3}$  and the vector $\bm{\bar{B}}^\ve\in \RR^{3}$ are given by \cite{fil-rey-2014}
        \begin{equation}
            \label{eqDefABV}
	      \left\{\begin{aligned}
            & \bm{\bar{A}}^\ve := \frac{1}{\rho^\ve} \int_{\RR^{3}} \bm{A}(\bm{V}) f^\ve(v)
            %-\mathcal{M}_{\rho^\epsilon, \bm{u}^\epsilon, T^\epsilon}
            \, dv, && \bm{A}(\bm{V}) = \bm{V} \otimes \bm{V} - \frac{|\bm{V}|^2}{3} \bm{I}, \\
            & \bm{\bar{B}}^\ve := \frac{1}{\rho^\ve}\int_{\RR^{3}} \bm{B}(\bm{V}) f^\ve(v)
            %-\mathcal{M}_{\rho^\epsilon, \bm{u}^\epsilon, T^\epsilon} 
            \, dv, && \bm{B}(\bm{V}) = \frac12 \left [ |\bm{V}|^2 - (3+2)\right ] \bm{V}, 
	      \end{aligned}\right.
        \end{equation}
        and where we used the shorthand 
        \[ \bm{V}(v) = \frac{v-\bm{u}}{\sqrt{T}}. \]
        Therefore, depending on the order in $\ve$ of the truncation of the series \eqref{devChapEnsk}, we will obtain different hydrodynamic description of the fluid.

			\medskip
			
			\subsubsection{Zeroth order: compressible Euler system}
			
			  At zeroth order with respect to  $\ve$, we have $f^\ve = \mathcal{M}_{\rho, \bm{u}, T}$. 
			  This distribution is isotropic in $v-\bm{u}$ and its odd moments with respect to $(v-\bm{u})$ are all equal to zero. Since the matrix $\bm{\bar{A}}^\ve$ is traceless, we then have that 
			  \[ 
			    \bm{\bar{A}}_{Euler} = \bm{0}_{M_{3}}.
			  \]
			  Moreover, since $ \bm{\bar{B}}^\ve$ involves odd, centered moments of $f^\ve$, we also obtain
			  \[
			    \bm{\bar{B}}_{Euler} = \bm{0}_{\RR^{3}}.
			  \]
			  Hence, the moments $(\rho, \bm{u}, T)$ are solution to the compressible Euler system \eqref{eqHydroClosedEuler}.
			  We notice in particular that the Maxwellian distribution in \eqref{devChapEnsk} is independent of $\ve$.
			  
			\medskip

			\subsubsection{First order: compressible Navier-Stokes system}
			
			  Going to the next order in $\ve$, we plug-in the expansion \eqref{devChapEnsk} in the Boltzmann equation \eqref{eqCollision}. 
			  Since the Maxwellian distribution is an equilibrium of the collision operator (according to \ref{hypEquilib}), the fluctuation $g^{(1)}$ is given by 
			  \begin{equation}
			    \label{eqCollisionExp}
			    \partial_t \mathcal M_{\rho,\bm{u},T} + v \cdot \nabla_x \mathcal M_{\rho,\bm{u},T} \,=\, \LL_{\mathcal M_{\rho,\bm{u},T}} \, g^{(1)}  + \mathcal{O}(\ve),
			  \end{equation}
			  where $(\rho, \bm{u}, T)$ are the unknowns
                          to be determined  and $\LL_{\mathcal M}$ is the linearized\footnote{Namely the Frechet derivative of the collision operator.} collision operator around the Maxwellian distribution.
					  
        Besides, using a logarithmic derivative, one has\footnote{ The Double-Dot product for matrices is the Frobenius inner product, and it is calculated computing the Hadamard product and adding up all the elements of the resulting matrix.}
			  \begin{multline*}
			    \partial_t \mathcal M_{\rho,\bm{u},T} + v \cdot \nabla_x \mathcal M_{\rho,\bm{u},T} \, = \, \\
			      \mathcal M_{\rho,\bm{u},T} \left [ \frac{1}{\rho}\partial_t \rho +\frac{1}{\rho} v \cdot \nabla_x \rho + \frac{1}{\sqrt{T}} \left ( \bm{V} \cdot \partial_t \bm{u} + \bm{V}\otimes v : \nabla_x \bm{u}\right ) + \frac{1}{2 T} \left (|\bm{V}|^2 - 3\right ) \left (\partial_t T + v\cdot \nabla_x T\right ) \right ].
			  \end{multline*}
			  Then, using the conservation laws
              \eqref{eqHydroNonClosed}, we replace the
              time derivatives by  spatial ones, and drop
              the terms of order $\ve$ in
              \eqref{eqCollisionExp}, hence we find after some computations that \cite{lev-1996, lev-mor-nad-1998}%\cite{fil-rey-2014, lev-mor-nad-1998} 
			  \begin{equation}
			    \label{eqLinearG1}
			    \LL_{\mathcal M_{\rho,\bm{u},T}} g^{(1)} =  \mathcal M_{\rho,\bm{u},T}\left [ 
                \frac{1}{2}
                \bm{A}(\bm{V}):\bm{D}(\bm{u}) + 2 \bm{B}(\bm{V}) \cdot \nabla_x \sqrt{T}\right ],
			  \end{equation}
			  where $\bm{A}$, $\bm{B}$ and $\bm{V}$ are defined in \eqref{eqDefABV} and the traceless \emph{deformation tensor} $\bm{D}$ of $\bm{u}$ is given by \cite{lev-mor-nad-1998}	  
			  \[\bm{D}(\bm{u}) = \nabla_x \bm{u} + \left (\nabla_x \bm{u}\right )^\intercal - \frac{2}{3} \left (\diverg_x \bm{u}\right ) \bm{I}. \]
			  
			  Moreover, using the hypothesis \ref{hypConservations} on the conservation laws of the collision operator, it is possible to show that linear combinations of collisional invariants form exactly the kernel of the linear operator $\LL_{\mathcal M_{\rho,\bm{u},T}}$. In particular, we have the orthonormal family 
			  \[ 
			    \ker \LL_{\mathcal M_{\rho,\bm{u},T}} = \vect \left \{ \frac{1}{\rho},  \frac {\bm{V}}{\rho}, \frac{1}{2\rho} \left (|\bm{V}|^2 - 3 \right ) \right \}.
			  \]
			  Using the orthogonality properties of the moments of a Maxwellian distribution, we have that on $L^2 \left (\mathcal M_{\rho,\bm{u},T} \right )$
			  \[ \bm{A} (\bm{V}), \bm{B}(\bm{V}) \perp \ker \LL_{\mathcal M_{\rho,\bm{u},T}}. \]
			  Since the operator $\LL_{\mathcal M_{\rho,\bm{u},T}}$ is invertible on the orthogonal of its kernel and finally  using \eqref{eqLinearG1}, it yields
			  \begin{equation}
			    \label{eqG1}
			    g^{(1)} =  \frac{1}{2} 
                \LL_{\mathcal M_{\rho,\bm{u},T}}^{-1}\left ( \mathcal M_{\rho,\bm{u},T} \, \bm{A} \right ): \bm{D}(\bm{u}) \,+\, 2 \,\LL_{\mathcal M_{\rho,\bm{u},T}}^{-1}\left ( \mathcal M_{\rho,\bm{u},T} \, \bm{B} \right ) \cdot \nabla_x \sqrt{T}. 
			  \end{equation}
			  We can then plug this expression into the definition \eqref{eqDefABV} to obtain using some classical symmetry properties of the collision operator \cite{lev-mor-nad-1998} that
			  \begin{equation}
			    \label{eqABNS}
			    \left\{\begin{aligned}
                %Before they were wih the minus sign, now I use the plus sign
			      & \bm{\bar{A}}^\ve_{NS} := \frac{1}{\rho} \int_{\RR^{3}} \bm{A}(\bm{V}) \mathcal M_{\rho,\bm{u},T}(v)\left [1 + \ve \, g^{(1)}(v) \right ] dv = - \ve \frac{\mu}{\rho \, T} \bm{D}(\bm{u}), \\ 
	          & \bm{\bar{B}}^\ve_{NS} := \frac{1}{\rho}\int_{\RR^{3}} \bm{B}(\bm{V}) \mathcal M_{\rho,\bm{u},T}(v)\left [1 + \ve \, g^{(1)}(v) \right ] dv = - \ve \frac{\kappa}{\rho \,T^{3/2}} \nabla_x T.
	        \end{aligned}\right.
			  \end{equation}
	      The scalar quantities $\mu$ and $\kappa$ in \eqref{eqABNS}, respectively the \emph{viscosity} and the \emph{thermal conductivity}, are given by \cite{lev-mor-nad-1998, lev-1996} 
	      % \begin{gather*}
	      %   \mu := - \frac{T}{2} \int_{\RR^{3}} \mathcal M_{\rho,\bm{u},T}(v) \bm{A}(\bm{V}):\LL_{\mathcal M_{\rho,\bm{u},T}}^{-1}\left ( \mathcal M_{\rho,\bm{u},T} \, \bm{A} \right )(v) \, dv, \\
	      %   \kappa := - T \int_{\RR^{3}} \mathcal M_{\rho,\bm{u},T}(v) \bm{B}(\bm{V}) \cdot \LL_{\mathcal M_{\rho,\bm{u},T}}^{-1}\left ( \mathcal M_{\rho,\bm{u},T} \, \bm{B} \right )(v) \, dv.
	      % \end{gather*}  
          \begin{gather*}
          %Before they were wih the minus sign, now I use the plus sig
	        \mu := \frac{T}{10} \int_{\RR^{3}}  \bm{A}(\bm{V}): \mathcal M_{\rho,\bm{u},T}(v) \LL_{\mathcal M_{\rho,\bm{u},T}}^{-1}\left ( \mathcal M_{\rho,\bm{u},T} \, \bm{A} \right )(v) \, dv, \\
	        \kappa := \frac{T}{3} \int_{\RR^{3}}  \bm{B}(\bm{V}) \cdot \mathcal M_{\rho,\bm{u},T}(v) \LL_{\mathcal M_{\rho,\bm{u},T}}^{-1}\left ( \mathcal M_{\rho,\bm{u},T} \, \bm{B} \right )(v) \, dv.
	      \end{gather*}	
          
	      They depend on the collision kernel of the model. For example, for the Boltzmann operator in the 3D hard-sphere case, it can be shown \cite{Golse:2005} that there exists some positive constants $\mu_0$, $\kappa_0$ such that
	      \[ \mu = \mu_0 \sqrt T \quad \text{and} \quad \kappa = \kappa_0 \sqrt{T}.\]
	      %In the ES-BGK case, we have \cite{str-2005}
	     % \[ \mu = \frac{1}{1-\beta} \frac{\rho \, T}{\nu} \quad \text{and} \quad \kappa = \frac{5}{2} \frac{\rho \, T}{\nu}.\]
	      
		    Finally, the evolution of the macroscopic quantities at first order with respect to $\ve$ is given by the compressible Navier-Stokes equations 
			  \begin{equation}
			    \label{eqHydroClosedNS}
					\left\{ \begin{aligned}
						& \partial_t \rho + \diverg_x (\rho\, \bm{u} ) \,=\, 0, 
						\\
						& \,
						\\
						& \partial_t(\rho \, \bm{u} ) + \diverg_x \left (\rho \, \bm{u} \otimes \bm{u} +  \rho \, T \, \bm{I} \right ) \,=\, \ve \diverg_x \left ({\mu} \,\bm{D}(\bm{u})\right ), 
						\\
						& \,
						\\
						& \partial_t E + \diverg_x \left ( \bm{u} \left ( E +\rho \, T\right ) \right ) \,=\, \ve \diverg_x \left ( \mu  \, \bm{D}(\bm{u}) \cdot \bm{u} + \kappa  \nabla_x T \right ).
					\end{aligned} \right.
				\end{equation}
				%where the matrix $\bm{\sigma} := -\mu  \, \bm{D}(\bm{u})$ is the so-called \emph{viscosity tensor} and the vector $\bm{q} := -\kappa  \nabla_x T$  the \emph{heat flux}.

			\subsubsection{Second order: Burnett equations}
				
				Pushing the expansion \eqref{devChapEnsk} at second order in $\ve$, we can use the same type of argument that for the compressible Navier-Stokes system to obtain another correction of the compressible Euler equations: the Burnet system.
				Although this system is ill-posed \cite{Golse:2005}, the computation of its coefficients is still possible.
				We have for the BGK case \cite{str-2005, fil-rey-2014}
				\begin{align}
				  \bm{\bar{A}}^\ve_{Burnett} & := \frac{1}{\rho}\int_{\RR^{3}} \bm{A}(\bm{V}) \mathcal M_{\rho,\bm{u},T}(v)\left [1 + \ve \, g^{(1)}(v) + \ve^2  g^{(2)}(v) \right ] dv \notag \\
				    & = - \ve \frac{\mu}{\rho \,T} \bm D( \bm{u} ) 
				        - 2 \ve^2 \frac{\mu^2}{\rho^2 T^{2}}\bigg \{ -\frac{T}{\rho} {\rm Hess}_x(\rho ) + \frac{T}{\rho^2} \nabla_x \rho \otimes \nabla_x \rho - \frac{1}{\rho} \nabla_x T \otimes \nabla_x \rho \notag \\
				    &  \qquad \qquad \qquad \qquad \qquad \qquad \ + \left (\nabla_x \bm{u} \right )  \left (\nabla_x \bm{u} \right )^\intercal - \frac{1}{3}\bm D( \bm{u} ) \diverg_x( \bm{u} ) + \frac{1}{T} \nabla_x T \otimes \nabla_x T\bigg \};  \label{eqAepsBurnett} \\				  
          \bm{\bar{B}}^\ve_{Burnett} & := \frac{1}{\rho}\int_{\RR^{3}} \bm{B}(\bm{V}) \mathcal M_{\rho,\bm{u},T}(v)\left [1 + \ve \, g^{(1)}(v) + \ve^2  g^{(2)}(v) \right ] dv \notag \\
				    & = - \ve \frac{\kappa}{\rho \,T^{3/2}} \nabla_x T 
				        - \ve^2 \frac{\mu^2}{\rho^2 T^{5/2}}\bigg \{ +\frac{25}{6} \left (\diverg_x \bm{u}\right ) \nabla_x T \notag \\
				    &  \qquad \qquad \qquad \qquad \qquad \qquad \quad \ \, - \frac{5}{3} \left [T \diverg_x \left (\nabla_x \bm{u}\right ) + \left (\diverg_x\bm{u} \right ) \nabla_x T + 6 \left (\nabla_x \bm{u} \right )  \nabla_x T\right ] \notag \\
				    &  \qquad \qquad \qquad \qquad \qquad \qquad \quad \ \, + \frac{2}{\rho} \bm D( \bm{u} ) \, \nabla_x \left (\rho \, T \right ) + 2 \, T \diverg_x \left (\bm D( \bm{u} )\right ) + 16 \bm D( \bm{u} ) \,\nabla_x T \bigg \} .
            \label{eqBepsBurnett}
				\end{align}

            \subsection{The moment realizability criterion}
	 	  \label{subFluid2Kin}

            The matrix $\bm{\bar{A}}^\ve$ and the vector $\bm{\bar{B}}^\ve$ will allow us to define our hydrodynamic break down criterion.  Let us set the vector of the reduced collisional invariants for $\bm{V}=(v-\bm{u^\ve})/\sqrt{T^\ve}$,
            \[
            \bm{m} := \left (1, \bm{V}, \left (\frac{2}{5}\right )^{1/2}\left (\frac{|\bm{V}|^2}{2} - \frac{5}{2}\right ) \right ).
            \] 
            Similarly to \cite{lev-mor-nad-1998, tiw-2000}, now we define the so-called \emph{moment realizability matrix} $\bm{M}$ by setting
            \begin{equation}
                \label{eqMomRealiz}
                \bm{M} := \frac{1}{\rho^\ve}\int_{\RR^{3}} \left(\bm{m} \otimes \bm{m} \right) f^\ve(v) \, dv.
            \end{equation}
            Indeed, if the distribution function is nonnegative, then the quantity $\int \phi^2(v) f \, dv $ is also nonnegative \cite{tiw-2000}, for any function $\phi(v)$. The same is true if we consider $\bm{m}=\bm{m}(\bm{v})$ a columns vector of $c$ polynomials, and we set $\bm{\phi}=\bm{a}^\intercal \bm{m}$ for a generic $\bm{a}\in R^c$, and we have
            \begin{equation}
                \bm{a}^\intercal\left(\int_{\RR^{3}}\bm{m}\otimes \bm{m} f \, dv \right)\bm{a}=\int_{\RR^{3}} (\bm{a}^\intercal \bm{m})^2f\, dv \geq 0.
            \end{equation}
            Then if $f$ is nonnegative and not identically zero, then  for every $\bm{m}$ the $c\times c$ matrix
            \begin{equation}\label{eq1005}
                \bm{M}=\frac{1}{\rho}\int_{\RR^{3}} \bm{m}\otimes \bm{m} f\, dv,
            \end{equation}
            is positive semi-definite.\\
            In the Appendix \ref{appendix:mrm} we prove that under the assumption of being sufficiently close enough to the equilibrium distribution function (i.e. the Maxwellian) then the quantity $\bm{M}$ is positive semi-definite if and only if the following quantity is also positive semi-definite:
            \begin{equation}\label{eqMomRealizReduced}
                \mathcal{V}_{\ve} = \bm{I} + \bm{\bar{A}}^\ve - \frac{2}{3}\bm{\bar{B}}^\ve\otimes \bm{\bar{B}}^\ve.
            \end{equation}
            Then the idea is to monitor the evolution of the matrix $\mathcal{V}_{\ve}$ \cite{tiw-1998, lev-mor-nad-1998}: indeed the moment realizability criterion states that the fluid dynamic description has broken down when the perturbation is too large that $\mathcal{V}_{\ve}$, and thus $\bm{M}$, is no longer positive semi-definite. More details on how to build a hierarchical solver are discussed in the next section.

        \subsection{Three layers hierarchical model}
	 	  \label{subSwitching}
        Now we present a hybrid decomposition numerical method able to automatically switch between three different layers: Euler equations (Layer 0), kinetic model with ES-BGK operator (Layer 1) and kinetic model with Boltzmann operator (Layer 2).

    \bigskip
    
        \begin{center}           
            \begin{tabular}{||c|c|c|c||}
                \hline
                LAYER   & MODEL & Computational Cost \& Accuracy \\[0.5 ex]
                \hline \hline
                2 & Kinetic model with Boltzmann operator & Most accurate \& highest computational cost \\
                1 & Kinetic model with ES-BGK operator & \\
                0 & Euler equations & Least accurate \& lowest computational cost\\
                \hline
            \end{tabular}
        \end{center}

    \bigskip
        
        Now from the remarks of the previous section, let us define a criterion to determine the appropriate model -- fluid or kinetic -- to be used. The key idea is the following: we evaluate the quantity $\mathcal{V}_{\ve}$ defined in Equation \ref{eqMomRealizReduced} for three orders in $\epsilon$ of the Chapman-Enskog expansion, and, since the greater the order of expansion, then the greater the accuracy of the model, we associate to each expansion order a different model. In particular, to 0-th order in the Chapman-Enskog expansion we associate the Euler equations (lowest accuracy), to 1-st order we associate the ES-BGK equation, and to 2-nd order we associate the full Boltzmann equation (highest accuracy). 
        Then, we will have three different matrices accordingly to the expansion order in $\epsilon$, that are: $\mathcal{V}_{\ve^0}$ (associated to 0-th expansion order, and to Euler equations), $\mathcal{V}_{\ve^1}$ (associated to 1-st expansion order, ES-BGK equation) and $\mathcal{V}_{\ve^2}$ (associated to 2-nd expansion order, full Boltzmann equation). So, in each cell of our spatial domain we can evaluate these three matrices, and we compare them in order to decide if we have to switch from one level to another one. If in a given cell the matrix $\mathcal{V}_{\ve^i}$ is \emph{too} different from matrix $\mathcal{V}_{\ve^{i-1}}$, then it means that the layer $i-1$ is not accurate, and then we have to use the layer $i$. In order to establish if the previous matrices are \emph{too} different, we will compare their eigenvalues. This is the main idea behind our hybrid domain decomposition numerical scheme.

Let's consider the zeroth order Chapman-Enskog expansion with respect to $\ve$: we get that $\bm{\bar{A}}_{0} = \bm{0}_{M_{3}}$ and $\bm{\bar{B}}_{0} = \bm{0}_{R^{3}}$ and then%. Moreover, we also have in that case  $\bar{C}^\ve = 1$ and then
				\[	
					\mathcal{V}_{_{\ve^0}} = \bm{I}. 
				\]

 On the other hand, consider  the first order model. By cutting the Chapman-Enskog expansion \eqref{devChapEnsk} at the first order with respect to $\ve$, we can compute explicitly the matrix $\mathcal{V}_{\ve^1}$ obtaining
		    %We have in this case using the expressions \eqref{eqABNS} and by symmetry arguments that $\bar C^\ve = 1$ \cite{tiw-2000}, hence
		    \begin{equation}
		      \label{eqMomRealizNS}
		      \mathcal{V}_{\ve^1} = \bm{I} - \ve \frac{\mu}{\rho T} \bm{D}\left (\bm{u}\right ) - \ve^2 \frac{2}{3} \frac{\kappa^2}{\rho^2 T^3} \nabla_x T \otimes \nabla_x T,
		    \end{equation}
		    where $(\rho, \bm{u}, T)$ are solution to the Navier-Stokes equations \eqref{eqHydroClosedNS}. 

     Hence, we claim that the layer 0 (Euler equations) is correct when the matrix $ \mathcal{V}_{\ve^1}$ behaves like the matrix $\mathcal{V}_{\ve^0}={\bm I}$, that is, it is  positive definite and if \emph{its eigenvalues are close to $1$ or not}: the Layer 0 description of the fluid will be considered incorrect if the difference $\lambda_A$ between the eigenvalues $\lambda_{\ve^1}$ of $\mathcal{V}_{\ve^1}$ and the eigenvalues $\lambda_{\ve^0}=1$ of $\mathcal{V}_{\ve^0}$ is \emph{too} big. Then our criterion to switch from Euler description (Layer 0) to ES-BGK (Layer 1) is the following inequality   
    				\begin{equation}
    			  	\label{crit:CompEuler}
    			  	\left| \lambda_{A} \right|=\left |\lambda_{\ve^1}-1\right | > \eta_0, \qquad\forall \lambda_{\ve^1} \in\, {\rm Sp}(\mathcal{V}_{\ve^1}),
    				\end{equation}
    	where $\eta_0$ is a small parameter (and it is a problem dependent quantity); in fact if this conditions it is satisfied then it means that the eigenvalues associated to the first order of the expansions are too far from those associated to the zeroth order, and so the Layer 0 is not accurate. Then we decide to use this criterion in our code to automatically switch between Euler equations (Layer 0) and ES-BGK equation (Layer 1).

        For the switching criterion between the Layer 1 (ES-BGK equation) and Layer 2 (Boltzmann equation), we follow the same strategy as before. So we evaluate $\mathcal{V}_{\ve^1}$ and $\mathcal{V}_{\ve^2}$, we compute the respective eigenvalues, and we compare them. If they are \emph{too} different, then we switch from Layer 1 to Layer 2.
        In the paper \cite{fil-xio-2018}, the authors evaluated the dominant differences between the matrices associated to $\mathcal{V}_{\ve^1}$ and $\mathcal{V}_{\ve^2}$. In their work, the $\mathcal{V}_{\ve^1}$ was associated to Navier-Stokes equations, and $\mathcal{V}_{\ve^2}$ to the ES-BGK equation, while in our case, they are associated respectively to ES-BGK and Boltzmann equations. We report now the main differences between the matrices of the two levels in the Chapman-Enskog expansion, as reported in \cite{fil-xio-2018}

        \begin{equation}\label{eq1006}
            \lambda_B(t, \mathbf{x}):=\epsilon^2
            \left(\frac{|\nabla_{\mathbf{x}}T|^2}{T} + |\nabla_{\mathbf{x}}\mathbf{u}|^2 +
            \sqrt{\left(|\Delta_{\mathbf{x}}\mathbf{u}|^2 + \left|\frac{\Delta_{\mathbf{x}}\rho}{\rho}\right|^2\right)\left(1+T^2\right)}\right).
        \end{equation}

        % For the switching criterion between the kinetic equation with BGK operator (Layer 1) and the kinetic equation with Boltzmann operator (Layer 2), we start noticing that for the main components in $\mathbf{\bar{A}}^\epsilon$ and $\mathbf{\bar{B}}^\epsilon$ in the expression of \eqref{eqMomRealizReduced}, the biggest differences between $\mathbf{\bar{A}}^\epsilon_{NS}$ and  $\mathbf{\bar{A}}^\epsilon_{Burnett}$ as well as $\mathbf{\bar{B}}^\epsilon_{NS}$ and  $\mathbf{\bar{B}}^\epsilon_{Burnett}$ are of order $\mathcal{O}(\epsilon^2)$ \cite{fil-xio-2018}. As indicator, to switch between the two regimes, it is possible to take into account only the main spatial derivatives appeared in those terms and roughly estimate their magnitudes, in order to speed up the tedious computation \cite{fil-xio-2018}. In particular, an indicator for switching from Navier-Stokes to Burnett is defined as
        % \begin{equation}
        %     \lambda_B(t, \mathbf{x}):=\epsilon^2
        %     \left(\frac{|\nabla_{\mathbf{x}}T|^2}{T} + |\nabla_{\mathbf{x}}\mathbf{u}|^2 +
        %     \sqrt{\left(|\Delta_{\mathbf{x}}\mathbf{u}|^2 + |\Delta_{\mathbf{x}}\rho/\rho|^2\right)\left(1+T^2\right)}\right).
        % \end{equation}
        We will use this quantity as indicator to switch from the ES-BGK equation to the full Boltzmann equation. In particular, we will say that the ES-BGK equation is not appropriate at some point $\left(t, \mathbf{x}\right)$ if:
        \begin{equation}
            |\lambda_B(t, \mathbf{x})|>\eta_1,
        \end{equation}
        where $\eta_1$ is a small parameter (and it is a problem dependent quantity). If this inequality is satisfied then we switch to the Boltzmann equation.
				  
      \begin{remark}
        Note that this last two criteria must be also verified when going from Layer 2 to Layer 1, and from Layer 1 to Layer 0, respectively. In particular, for moving from Layer 2 the following condition must be satisfied
        \begin{equation}
            |\lambda_B(t, \mathbf{x})|\leq\eta_1.
        \end{equation}
        In our code this last condition is sufficient to switch from Layer 2 to Layer 1.
        Analogously, for moving from Layer 1 to Layer 0 the following condition must be satisfied
        \begin{equation}
            |\lambda_A|=|\lambda_{\ve^1}(t, \mathbf{x})-1|\leq\eta_0.
        \end{equation}
        This last condition is necessary but not sufficient to switch from Layer 1 to Layer 0, as we explain in the following.
      \end{remark}

	 In fact, knowing the full kinetic description of a gas (Layer 1 or Layer 2), there exists several methods \cite{sai-2009} to decide how far this gas is from the thermal equilibrium, \ie the fluid regime (Layer 0). \\
	    In our code we compare $f^\ve$, solution to the collisional equation \eqref{eqCollision} with the truncated Chapman-Enskog distribution $f^\ve_k$ given by \eqref{devChapEnsk}, whose moments match the one of $f^\ve$, and whose order $k$ corresponds to the order of the macroscopic model considered.
%      If $f$ and $g$ are two function of the variable $v\in \RR^d$, we call \emph{relative entropy} of $f$ with respect to $g$ the quantity
%	    \[
%	      \mathcal H[f|g] := \int_{\RR^{3}} f(v) \log \left ( \frac{f(v)}{g(v)}\right ) dv.
%	    \]
%	    This quantity (which can be computed easily numerically) almost defines a norm between distribution functions and the thermal equilibrium associated. Indeed, we have the classical inequality:
%	    
%	    \begin{theorem} [Csiszar-Kullback-Pinsker inequality, \cite{Csiszar:67}]
%	    
%	      If $\mathcal M$ is the Maxwellian distribution associated to $f$, then
%	      \[
%	        \|f - \mathcal M\|_{L^1}^2 \leq 2 \mathcal H[f|\mathcal M].
%	      \]
%	      
%	    \end{theorem}
	    
%	    From this theorem, it is then easy to derive a \emph{kinetic to fluid} criterion. Let us introduce a small parameter $\delta_0>0$ and take $f^\ve(t,x,\cdot)$ a solution to the kinetic equation  at a given time $t >0$ and position $x \in \Omega$. 
%	    If $(\rho^\ve, \bm{u}^\ve, T^\ve)$ are the associated macroscopic quantities, then we say that the kinetic regime is close to a fluid one if we have
%	    \begin{equation}
%	      \label{eqKintoFluid1}
%	      \mathcal H\left [f^\ve(t,x,\cdot)|\mathcal M_{\rho^\ve,\bm{u}^\ve,T^\ve}(t,x,\cdot)\right ] \leq \delta_0.
%	    \end{equation} 

     Our criterion is then the following: the kinetic description given by the ES-BGK equation (Layer 1) at point $(t,x)$ corresponds to an \emph{hydrodynamic closure of order 0} (Layer 0) if
		  \begin{equation}
		    \label{eqKintoFluid1}
		    \left \| f^\ve(t,x,\cdot) - \mathcal M_{\rho,\bm{u},T}^\ve(t,x,\cdot) \right \|_{L^1_v} \leq \delta_0,
		  \end{equation}
		  where $\delta_0$ is a small parameter (which is a problem dependent quantity) 
		  and if the eigenvalues of the moment realizability matrix $\mathcal V_{\ve^1}$ computed using the moments of $f^\ve$ does not verify the criterion \eqref{crit:CompEuler}.

	\section{Numerical schemes}
	  \label{secNumSim}   

            \subsection{Spatial, velocity and time settings}
            We consider a 2D spatial cartesian grid, of intervals $[0, L_x]$ and $[0, L_y]$, and discretized with $N_x$ and $N_y$ points in the $x$ and $y$ direction respectively.   
            
            The spatial cartesian grid is defined by nodes $(x_i, y_j)=(i\Delta x, j \Delta y)$ and cells
            \begin{equation}
                I_{ij} = \left(x_{i-\frac{1}{2}}, x_{i+\frac{1}{2}} \right) \times \left(y_{j-\frac{1}{2}},y_{j+\frac{1}{2}} \right).
            \end{equation}
            
            The discrete velocity grid which represents the set of admissible velocities in the discrete model. Let $\mathcal{K}\subseteq\ZZ^3$ be a set of $N_v$ integer vectors, and let
            \begin{equation}
                \nu = \{v_k \in \RR^3, k \in \mathcal{K} \},
            \end{equation}
            be a discrete-velocity grid of $N_v$ points indexed by $k=(k_1, k_2, k_3)\in \mathcal{K}$, and defined by
            \begin{equation}
                v_k=(v_{k_1}, v_{k_2}, v_{k_3})=(k_1\Delta v_1, k_2 \Delta v_2, k_3 \Delta v_3),
            \end{equation}
            where $\Delta v_1$, $\Delta v_2$, $\Delta v_3$ are three positive numbers characterising the size of the velocity mesh. To simplify notations we denote by $\Delta v_{\mathcal{K}}=\Delta v_1 \Delta v_2 \Delta v_3$. The velocity distribution $f$ (Layer 1 or Layer 2) is then replaced by a vector $\mathbf{f} = (f_k(t,x,y))_{k\in \mathcal{K}} \in \RR^{N_x}$ where each component $f_k(t,x,y)$ is assumed to be an approximation of $f(t,x,y, v_k)$.

             Each cell $(x_i, y_j)$ is in one of three following regimes: Layer 0 is hydrodynamic (described by Euler equations), Layer 1 is the ES-BGK equation, or Layer 2 is the full Boltzmann equation. We identify with the symbol $r^n_{ij}=\{0, 1, 2\}$ the regime of the cell $(x_i, y_j)$ at the timestep $n$.
             If $r^n_{ij}=2$ then the fluid in the cell $(x_i, y_j)$, at time $t^n$, is described by a distribution function $f^{2}_k(t^n, x_i, y_j)$ which is updated using the spectral method for the Boltzmann operator; if $r^n_{ij}=1$ then the fluid in the cell $(x_i, y_j)$ is described by a distribution function $f^{1}_k(t^n, x_i, y_j)$ which is updated using the ES-BGK operator; if $r^n_{ij}=0$ then the fluid in the cell $(x_i, y_j)$ is described by the moments $U_{t^n, x_i, y_j}=(\rho, \vec{u}, T)_{ij}$ and they are updated using the Euler equations. To each cell we can associate a Maxwellian $\mathcal{M}_{\rho, \vec{u}, T}^n(t^n, x_i,y_j)$, computed using the moments $\vec{m}(t^n, x_i, y_j)=\left[\rho, \vec{u}, T\right](t^n, x_i, y_j)$ (density, velocity and temperature) obtained by integrating the distribution function in the cell itself.
             
            % Now let us discuss the discrete velocity grid which represents the set of admissible velocities in the discrete model. Let $\mathcal{K}\subseteq\ZZ^3$ be a set of $N_v$ integer vectors, and let
            % \begin{equation}
            %     \nu = \{v_k \in \RR^3, k \in \mathcal{K} \}
            % \end{equation}
            % be a discrete-velocity grid of $N_v$ points indexed by $k=(k_1, k_2, k_3)\in \mathcal{K}$, and defined by
            % \begin{equation}
            %     v_k=(v_{k_1}, v_{k_2}, v_{k_3})=(k_1\Delta v_1, k_2 \Delta v_2, k_3 \Delta v_3)
            % \end{equation}
            % where $\Delta v_1$, $\Delta v_2$, $\Delta v_3$ are three positive numbers characterising the size of the mesh. To simplify notations we denote by $\Delta v_{\mathcal{K}}=\Delta v_1 \Delta v_2 \Delta v_3$. The velocity distribution $f$ (Layer 1 or Layer 2) is then replaced by a vector $\mathbf{f} = (f_k(t,x,y))_{k\in \mathcal{K}} \in \RR^{N_x}$ where each components $f_k(t,x,y)$ is assumed to be an approximation of $f(t,x,y, v_k)$.\\

            Our implementation is characterised by one single time step $\Delta t = 0.1 dx$, used to update both the hydrodynamic solution, and the kinetic ones (ES-BGK and Boltzmann operator), where $dx$ is the spatial discretization (in our simulations we use the same discretization in both spatial directions, thus $dx=dy$). The time integration, discussed in further detail in the next sections, is performed using a TVD Runge Kutta 2 \cite{got-shu-1996, shu-osh-1989} for the hydrodynamic regime and the so-called IMEX PR(2,2,2) \cite{par-rus-2001} numerical scheme for the kinetic levels.

		  \subsection{Euler equations - Layer 0}		  
			  In this subsection, we shall focus on the space and time discretization of the system of $N$ conservation laws
			  \begin{equation}
				  \label{sysConservLaws}
				  \left \{ \begin{aligned}
				   & \frac{\partial U}{\partial t} + \diverg_x F(U) = 0, \ \forall \, (t,x) \in \RR_+ \times \Omega, \\
				   & \, \\
				   & U(0,x) = U_0(x), 
				   \end{aligned} \right.
				\end{equation}
				for a smooth function $F : \RR^N \to M_{N\times d_x}(\RR)$ and a Lipschitz-continuous domain $\Omega \subset \RR^{d_x}$. 

                Here we apply a finite volume scheme, where the reconstruction of the variables at the interfaces is performed using a CWENO scheme \cite{lev-pup-rus-2000} and the flux is reconstructed using the Lax-Friedrich flux \cite{pup-sem-vis-2022}.
                In particular, we use the 1D reconstruction in both the x and y directions\footnote{Such an approach is also used to compute the linear free transport part of the kinetic layers.}.

                The Euler equations in a 2D space can be expressed using the following conservation laws
                \begin{equation}
                    \frac{\partial U}{\partial t} + \frac{\partial \mathbf{h}}{\partial x} +  \frac{\partial \mathbf{g}}{\partial y} =0,
                \end{equation}
                where
                \begin{equation}
                    U=
                    \begin{bmatrix}
                    \rho \\
                    \rho u_x\\
                    \rho u_y\\
                    E\\
                    \end{bmatrix}, \ \ 
                     \mathbf{h}=
                    \begin{bmatrix}
                    \rho u_x \\
                    \rho u_x^2 + P\\
                    \rho u_x u_y\\
                    u_x (E+P)\\
                    \end{bmatrix}, \ \ 
                    \mathbf{g}=
                    \begin{bmatrix}
                    \rho u_y \\
                    \rho u_y u_x\\
                    \rho u_y^2 + P\\
                    u_y (E+P)\\
                    \end{bmatrix},
                \end{equation}
                and $P=\left(\Xi - 1 \right)\left(E-\frac{1}{2}\rho\left(u_x^2 + u_y^2\right)\right)$ is the pressure, $\rho$, $u_x$, $u_y$, $E$ and $\Xi$ are, respectively, the density, velocity along x-direction, velocity along y-direction, energy and specific heat ratio.\\

                The time integration is performed using the second order TVD Runge Kutta scheme presented in \cite{got-shu-1996, shu-osh-1989}, which is the following one
                \begin{align}
                    U^{(1)} &= U^n + \Delta t L( U^n),\\
                    U^{n+1} &= \frac{1}{2} U^n + \frac{1}{2} U^{(1)} + \frac{1}{2} \Delta t L(U^{(1)}),
                \end{align}
                where, in our case, $L(U)= - \frac{\partial \mathbf{h}}{\partial x} - \frac{\partial \mathbf{g}}{\partial y}$.

			\medskip

\subsection{ES-BGK model - Layer 1}

  			  We now focus briefly on the evolution of the ES-BGK equation
			  \begin{equation} 
					\label{eqES-BGK}
					\left\{ \begin{aligned}
					  & \frac{\partial f^\ve}{\partial t} + v \cdot \nabla_x f^\ve \,=\, \frac{\nu}{\ve} \left( \G[f] \,-\, f \right), 
					  \\
				  	&  \;
					  \\
					  & f^\ve(0, x, v) = f_{0}(x,v),
					\end{aligned} \right.
				\end{equation}
                Here we will consider $f=f^1$, i.e. the distribution function associated to the Layer 1, which is updated using the ES-BGK operator. In the following we denote by $f_{ij,k}^n=f_k(t^n, x_i, y_j), \forall k \in \mathcal{K}$ and by $M_{ij,k}^n[\mathbf{f}]=M_k[\mathbf{f}](t^n, x_i, y_j)$ the discrete Maxwellian equilibrium. We further denote with $\nu_{ij}^n=\nu(\rho_{ij}^n, T_{ij}^n)$ the collision frequency.\\

                We numerically integrate the  ES-BGK equation using the same IMEX strategy presented in \cite{fil-jin-2011}. The advantage of such a time discretization is that it is uniformly stable with respect to small Knudsen number, thus removing the stiffness of the relaxation term, yet the implicit relaxation term can be solved explicitly, thanks to the special structure of the relaxation term. In addition, the authors have proved the asymptotic-preserving property of this scheme, which allows it to capture the fluid dynamic behaviour without resolving numerically the small Knudsen number. Here we now present the IMEX strategy as illustrated in \cite{fil-jin-2011}, but in our code we integrate in time using the PR(2,2,2) scheme \cite{par-rus-2001}, that we will discuss at the end of this section. Applying a fist order implicit-explicit IMEX scheme to the ES-BGK equation one obtains \cite{fil-jin-2011}
                \begin{equation}\label{eq1000}
                    \begin{cases}
                        \displaystyle \frac{f^{n+1}-f^n}{\Delta t} + v \cdot \nabla_x f^n = \frac{\nu^{n+1}}{\epsilon}\left( \mathcal{G}\left[f^{n+1}\right] - f^{n+1} \right),\\
                        \displaystyle f^0(x, v) = f_0 (x, v).
                    \end{cases}
                \end{equation}
                After some manipulations, one can obtain the following scheme \cite{fil-jin-2011}
                \begin{equation}
                    \begin{cases}
                       \displaystyle  U^{n+1} = \int_{\RR^{d_v}}\phi(v)
                        \left(f^n - \Delta t v \cdot \nabla_x f^n \right) dv, \\
                        \begin{aligned}
                            \displaystyle \Sigma^{n+1} = \frac{\epsilon}{\epsilon + (1 - \beta)\nu^{n+1}\Delta t }
                            \left(\Sigma^n - \Delta t \int_{\RR^{d_v}}v\otimes v \, v \cdot \nabla_x f^n dv\right)+\\
                            \frac{(1-\beta)\nu^{n+1}\Delta t}{\epsilon + (1-\beta)\nu^{n+1}\Delta t}\rho^{n+1}
                        \left(T^{n+1}I+ u^{n+1}\otimes u^{n+1} \right),
                        \end{aligned}
                        \\
                       \displaystyle f^{n+1} = \frac{\epsilon}{\epsilon + \nu^{n+1}\Delta t}\left[f^n - \Delta t v\cdot \nabla_x f^n \right] +
                        \frac{\nu^{n+1}\Delta t}{\epsilon + \nu^{n+1}\Delta t}\mathcal{G}\left[f^{n+1}\right],
                    \end{cases}
                \end{equation}
                where 
                \begin{equation}
                    \displaystyle \Sigma^{n+1}=\int_{\RR^{d_v}}v\otimes v f^{n+1}dv = \rho^{n+1}(\Theta^{n+1}+u^{n+1}\otimes u^{n+1}),
                \end{equation}
                and
                \begin{equation}
                    \displaystyle \Theta(t,x)\, =\, \frac{1}{\rho} \int_{{\RR}^{3}}  (v-\bm{u})\otimes (v-\bm{u})\,f(t,x,v) \,dv.
                \end{equation}
                Thus, even if \eqref{eq1000} is nonlinearly implicit, it can be solved explicitly.

                The system \eqref{eq1000} is a first order IMEX scheme, however we used the same strategy to implement the PR(2,2,2) second order IMEX \cite{par-rus-2001}, which is adoptable for any stiff systems of differential equation in the form
                \begin{equation}\label{eq1001}
                    y'=q(y) + \frac{1}{\epsilon} g(y),
                \end{equation}
                where $y=y(t)\in \RR^N$, $q, g: \RR^N\to \RR^N$. In our case $y$ is the distribution function, when considering only the temporal dependence, $q(y)$ is the transport part of the Boltzmann equation, and $g(y)$ is the ES-BGK operator.
                
                An Implicit-Explicit IMEX Runge Kutta scheme for system \eqref{eq1001} is of the form
                \begin{equation}
                    \begin{cases}
                        \displaystyle Y_i = y_n + \Delta t \sum_{j=1}^{i-1} \tilde{a}_{ij}q(t_n + \tilde{c}_j \Delta t, Y_j) + \Delta t \sum_{j=1}^\Lambda a_{ij}\frac{1}{\epsilon} g(t_n + c_j \Delta t, Y_j),\\
                        \displaystyle y_{n+1}=y_n + \Delta t \sum_{i=1}^\Lambda \tilde{\omega}_{i}q(t_n+\tilde{c}_i \Delta t, Y_i)+\Delta t \sum_{i=1}^\Lambda \omega_i \frac{1}{\epsilon} g(t_n + c_i \Delta t, Y_i).
                    \end{cases}                    
                \end{equation}
                The matrices $\tilde{A}=(\tilde{a}_{ij})$, $\tilde{a}_{ij}=0$ for $j\geq i$ and $A=(a_{ij})$ are $\Lambda \time \Lambda$ matrices such that the resulting scheme is explicit in $q$, and implicit in $g$. An IMEX Runge-Kutta scheme is characterized by these two matrices and the coefficient vectors $\tilde{c}=(\tilde{c}_1, ..., \tilde{c}_\Lambda)^\intercal$, $\tilde{\omega}=(\tilde{\omega}_1, ..., \tilde{\omega}_\Lambda)^\intercal$, $c=(c_1, ..., c_\Lambda)^\intercal$, $\omega = (\omega_1, ..., \omega_\Lambda)^\intercal$. They can be represented by a double \emph{tableau} in the usual Butcher notation,
                \begin{center}
                    \begin{tabular}{c|c}
                        $\tilde{c}$ & $\tilde{A}$\\                
                        \hline
                        \\
                                    & $\tilde{\omega}^\intercal$\\
                    \end{tabular}, \, \,  
                    \begin{tabular}{c|c}
                        $c$ & $A$\\                
                        \hline
                        \\
                            & $\omega^\intercal$\\
                    \end{tabular}.
                \end{center}
                The PR(2,2,2) \cite{par-rus-2001} is characterized by the following Butcher tableaus
                \begin{center}
                    \begin{tabular}{c|c c}
                        $0$ & $0$   & $0$\\   
                        $1$ & $1$   & $0$\\ 
                        \hline
                            & $1/2$ & $1/2$\\
                    \end{tabular},  \, \,
                    \begin{tabular}{c|c c}
                        $1-C$ & $1-C$      & $0$\\   
                        $C$   & $C-\delta$ & $\delta$\\ 
                        \hline
                            & $1/2$ & $1/2$\\
                    \end{tabular},
                \end{center}
                where $\delta = 1 - 1/(2C)$. In our scheme we use $C=1/\sqrt{2}$.\\

                Now it remains only to discuss the the spatial integration of the transport part in the Boltzmann equation, which is performed using the same CWENO scheme \cite{lev-pup-rus-2000, pup-sem-vis-2022} used for the hydrodynamic setting.

            \subsection{Boltzmann equation - Layer 2}

  			  We now focus on the evolution of the kinetic equation where we use the Boltzmann operator
			  \begin{equation} 
					\label{eqBoltzmann}
					\left\{ \begin{aligned}
					  & \frac{\partial f^\ve}{\partial t} + v \cdot \nabla_x f^\ve \,=\, \frac{1}{\epsilon}\Q(f^\ve), 
					  \\
				  	&  \;
					  \\
					  & f^\ve(0, x, v) = f_{0}(x,v).
					\end{aligned} \right.
				\end{equation}
                Here we consider $f=f^2$, i.e. the distribution function associated to the Layer 2. We shall adopt the strategy proposed in \cite{fil-jin-2010}. When the collision operator $\mathcal{Q}$ is the BGK operator, it is well known that even an implicit collision term can be solved explictly, using the property that $\mathcal{Q}$ preserves mass, momentum and energy. The main idea in \cite{fil-jin-2010} is to use this property penalizing the Boltzmann collision operator $\mathcal{Q}$ by the BGK operator as follows
                \begin{equation}
                    \mathcal{Q} = \left[\mathcal{Q} - \lambda \left(\mathcal{M}-f\right)\right] + \lambda \left[\mathcal{M}-f\right],
                \end{equation}
                where $\mathcal{M}$ is the Maxwellian associated to the distribution function $f$ and $\lambda$ is the spectral radius of the linearized collision operator of $\mathcal{Q}$ around $\mathcal{M}$.\\
                In this last equation, the first term on the right-hand side is less stiff (or not stiff at all) compared to the second one. So the first term can be discretized \emph{explicitly}, thus avoiding to invert the nonlinear Boltzmann operator $\mathcal{Q}$, and the second term can be treated explicitly.
                
                In \cite{fil-jin-2010} the authors propose the following split procedure in the equation \ref{eqBoltzmann}
                \begin{equation}
                    \frac{\mathcal{Q}(f)}{\epsilon} = 
                    \frac{\mathcal{Q}(f) - P(f)}{\epsilon}+ \frac{P(f)}{\epsilon},
                \end{equation}
                where the first term on the RHS is less stiff than $\frac{P(f)}{\epsilon}$. The quantity $P(f)$ is a linear operator asymptotically close to the source term $\mathcal{Q}(f)$, and that preserves the steady state $P(\mathcal{M})=0$.
                
                It is possible to perform a Taylor expansion leading to
                \begin{equation}
                    \mathcal{Q}(f) = \mathcal{Q}(\mathcal{M}) + 
                    \nabla \mathcal{Q}(\mathcal{M})(f-\mathcal{M}) + O(||f-\mathcal{M}||_H^2),
                \end{equation}
                and thus we can choose 
                \begin{equation}
                    P(f)\coloneqq \nabla \mathcal{Q}(\mathcal{M})(f-\mathcal{M}).
                \end{equation}
                However, since it is not always possible to compute $\nabla \mathcal{Q}(\mathcal{M})$ we can introduce an upper bound (or some approximation) $\beta$ of $||\nabla \mathcal{Q}(\mathcal{M})||$.
                
                At this point we have splitted the Boltzmann collision operator $\mathcal{Q}$ in a stiff and in a non-stiff term. Thus, we can adopt an IMEX approach like we have done in the previous section for the ES-BGK operator, where in this case, the explicit part is constituted by
                \begin{equation}
                    -\nabla_x f +\frac{\mathcal{Q}(f) - P(f)}{\epsilon},
                \end{equation}
                where $ P(f)= \beta(f-\mathcal{M})$;
                while the implicit term is
                \begin{equation}
                    \frac{P(f)}{\epsilon}.
                \end{equation}
                A first order IMEX scheme as presented in \cite{fil-jin-2010}  takes the following form
                \begin{equation}\label{eq1003}
                    \begin{cases}
                        \displaystyle \frac{f^{n+1}-f^n}{\Delta t} + v \cdot \nabla_x f^n = \frac{\mathcal{Q}(f^n) - P(f^n)}{\epsilon}+ \frac{P(f^{n+1})}{\epsilon}\\
                        \displaystyle f^0 (x, v)=f_0(x,v).
                    \end{cases}
                \end{equation}
                This system can be rewritten as
                \begin{equation}\label{eq1002}
                    f^{n+1} = \frac{\epsilon}{\epsilon + \beta^{n+1}\Delta t}\left[f^n - \Delta t v \nabla_x f^n\right] + \Delta t \frac{\mathcal{Q}(f^n) - P(f^n)}{\epsilon + \beta^{n+1}\Delta t} + \frac{\beta^{n+1}\Delta t}{\epsilon + \beta^{n+1}\Delta t}\mathcal{M}^{n+1},
                \end{equation}
                where $\beta^{n+1}=\beta(\rho^{n+1}, T^{n+1})$ and $\mathcal{M}^{n+1}$ is the Maxwellian associated to $f^{n+1}$. In our simulations we considered $\beta = 2\pi\rho$ \cite{fil-jin-2010}.
                
                The equation \eqref{eq1002} appears to be nonlinearly implicit, but it can be computed explicitly \cite{fil-jin-2010}. Indeed, multiplying it by $\phi(v)\coloneqq(1, v, |v|^2)$ (collisional invariants) and introducing the macroscopic quantity $U\coloneqq(\rho, \rho u , T)$ one gets
                \begin{equation}
                    U^{n+1} = \int \phi(v) (f^n - \Delta t v \cdot \nabla_x f^n)dv.
                \end{equation}
                From this last relation, the quantity $U^{n+1}$ can be obtained explicitly from $f^n$, and consequently it is possible to compute $\mathcal{M}^{n+1}$.\\
                Summing up, it is possible to compute explicitly $f^{n+1}$ using the first order IMEX scheme \eqref{eq1003}.
                
                At this point it is possible to generalize this approach to higher order IMEX scheme, and in particular, in our code we have adopted it for the PR(2,2,2) scheme, already discussed in the previous section.
                The spatial part of the Boltzmann equation has been computed using the already discussed CWENO scheme \cite{lev-pup-rus-2000, pup-sem-vis-2022}.

               % Let us discuss  the numerical approximation of the Botzmann operator $\mathcal{Q}(f)$.

                %\subsubsection{From kinetic equations to Discrete Velocity Models}\label{sub_Discretization}
                \subsubsection{Fourier-Galerkin spectral methods for the Boltzmann equation}\label{spectralMethods}
                The aim is to build a numerical method to compute the Boltzmann operator \eqref{eqQBoltz}, which depends only on the velocity variable $v$. We shall then neglects the time and space variables in this section.

                The starting point is the truncation of the integration domain of the Boltzmann integral \ref{eqQBoltz}, then we suppose that the distribution function $f$ has compact support on the ball $B_0(R)$, of radius $R$ centred in the origin. We have the classical result from \cite{PaRu:SINUM:2000}
                \begin{proposition}
                    Let the distribution function $f$ be compactly supported on the ball $B_0(R)$ of radius $R$ centered in the origin, then 
                    \begin{equation}
                        \text{Supp}_v(\Q(f,f)) \subset B_0 (\sqrt{2}R).
                    \end{equation}
                \end{proposition}
                It is possible to introduce a collision operator with cut-off in the following way                 
                \begin{equation}
                    \Q^R(f,g)(v)=\int \int_{B_0(R)\times \SSS^{d-1}}\left[f(v')g(v_*')-f(v)g(v_*)\right]B(v-v_*, \sigma)dv_* d\sigma, \ \ \ \forall v \in \RR^d
                \end{equation}
                It is important to have a spectral approximation without superposition of periods. In order to reach this goal, it is sufficient to have a distribution function $f(v)$ restricted on the cube $\left[-T, T\right]^d$, with $T\geq (2+\sqrt{2})R$. Then, we assume $f(v)=0$ on $\left[-T, T\right]^d \backslash B_0(R)$ and extend $f(v)$ to a periodic function on the set $\left[-T, T\right]^d$. We want to remark that the lower bound for $T$ can be improved: in fact, the choice $T=(3+\sqrt{2})R/2$ guarantees the absence of intersection between periods where the distribution function $f$ is different from 0. Unfortunately, the support of $f$ will increase in time, so we can just minimise the errors due to aliasing \cite{rey-2023}.
                
                In the following we will assume $T=\pi$ and hence $R=2\pi/(3+\sqrt{2})$ in order to simplify the notation.
                Projecting the collision operator onto the space of trigonometric polynomials we obtain
                \begin{equation}\label{eqColOpTrigPol}
                    \hat{Q}_k = \int_{\left[-\pi, \pi\right]^d} \Q^R(f_N, f_N)e^{-ik\cdot v}dv, \ \ \ k = -N, ..., N
                \end{equation}
                where $f_N$ is the truncated Fourier series of the distribution function $f$,
                \begin{equation}
                    \begin{aligned}                            
                        f_N(v)=\sum_{k=-N}^{N} \hat{f}_k e^{ik\cdot v},\\
                        \hat{f}_k = \frac{1}{(2\pi)^d}\int_{D_t} f(v) e^{-ik\cdot v}dv,
                    \end{aligned}
                \end{equation}
                where we have used just one index to denote the \emph{d}-dimension sums with respect to the vector $k=(k_1, ..., k_d)\in \ZZ^d$, and so we have
                \begin{equation}
                    \sum_{k=-N}^N \coloneq \sum_{k_1, ..., k_d = -N}^N.
                \end{equation}
                Now we can substitute the truncated Fourier series $f_N$ inside the expression \eqref{eqColOpTrigPol}, and after some computations we obtain \cite{PaRu:SINUM:2000}
                \begin{equation}\label{eqQk}
                \hat{Q}_k = \sum_{\substack{l,m=-N\\l+m=k}}^N \hat{f}_l \hat{f}_m \beta(l, m), \ \ \ k = -N, ..., N,
                \end{equation}
                where $\beta(l,m)=\hat{B}(l,m) - \hat{B}(m,m)$, and
                \begin{equation}
                    \hat{B}(l,m)=\int_{B_0(2\lambda \pi)}\int_{\SSS^{d-1}} |q|\sigma(|q|, \cos{(\theta)})e^{-i(l\cdot q^+ + m \cdot q^-)}d\omega \ dq,
                \end{equation}
                with
                \begin{equation}
                    \lambda = 2/(3+\sqrt{2}), \ \ \ q^+ = \frac{1}{2}(q+|q|\sigma), \ \ \ q^-=\frac{1}{2}(q-|q|\sigma), \ \ \ q = v - v_*.
                \end{equation}
                We would like to remark that the naive evaluation of \eqref{eqQk} requires $\mathcal{O}(n^2)$ operations, where $n=N^d$, which means that this spectral method is computationally expensive. The computation can be speed up by reducing the number of operations needed to evaluate the collision integral, introducing a FFT-based convolution approach in our spectral method.

                It is possible to find a new representation, called \emph{Carleman-like representation}, of the collision operator which conserves more symmetries when approximated on a bounded domain.
                \begin{lemma}
                    Introducing the change of variables $x=r\sigma/2$ and $y = v_* - v - x$ the collision operator can be rewritten as \cite{MoPa:2006}
                    \begin{equation}
                        \Q(f,f)(v)=\int_{x\in \RR^d}\int_{y\in \RR^d}\tilde{B}(x,y)\delta(x\cdot y)\left[f(v+y)f(v+x)-f(v+x+y)f(v)\right]dx \ dy,
                    \end{equation}
                    where
                    \begin{equation}\label{BetaTildeCarleman}
                        \tilde{B}(x,y) = \tilde{B}(|x|, |y|) = 2^{d-1}B\left(\sqrt{|x|^2+|y|^2},\frac{|x|}{\sqrt{|x|^2+|y|^2}}  \right)(|x|^2+|y|^2)^{-\frac{d-2}{2}},
                    \end{equation}
                where the last expression of $\tilde{B}(x,y)$ is valid on the manifold defined by $x\cdot y = 0$.
                \end{lemma}
                
                Let us notice that, if we assume the validity of the VHS model \eqref{eqVHS}, then we obtain $\tilde{B}(x,y)$ constant if $\gamma=0$ and $d=2$ or if $\gamma=1$ and $d=3$.
                
                Now let's discuss two possible truncation strategies \cite{rey-2023}. We start considering a bounded domain $\mathcal{D}_T=\left[-T,T\right]^d \ (0<T<+\infty)$. A first approach is removing the collisions related with some points out of the box. This is the natural first stage for deriving conservative schemes based on the discretization of the velocity. In this case it is not necessary a truncation on the modulus of $x$ and $y$ since we are imposing them to stay in the box. So we obtain 
                \begin{multline}
                    \Q^{tr}(f,f)(v) = \int \int_{\{x,\ y\in \ \RR^d \ | \ v+x,\ v+y,\ v+x+y\ \in \ \mathcal{D}_T\}}\tilde{B}(x,y)\delta(x\cdot y)\\
                    \left[f(v+y)f(v+x)-f(v+x+y)f(v)\right]dx \ dy,
                \end{multline}
                where $v\in \mathcal{D}_T$.\\
                This last operator satisfies the following weak form
                \begin{multline}
                    \int \Q^{tr}(f,f)(v) \phi(v) dv =
                    \frac{1}{4} \int \int \int_{\{v,\ x,\ y\in \ \RR^d \ | v,\ v+x,\ v+y,\ v+x+y\ \in \ \mathcal{D}_T\}}\tilde{B}(x,y)\delta(x\cdot y)\\
                    f(v+x+y)f(v)\left[\phi(v+y)+\phi(v+x)-\phi(v+x+y)-\phi(v)\right] dv \ dx \ dy,
                \end{multline}
                from which we can deduce the conservation of mass, momentum, energy and the validity of the H-theorem for the entropy. Doing this derivation we have changed the collision kernel by adding artificial dependencies on $v, v_*, v', v_*'$, and so the convolution-like properties are lost. This is the problem of this kind of truncation.

                A second possible approach is truncating the integration with respect to variables $x$ and $y$ by setting them to be in $B_0(R)$. In order to obtain all possible collisions we consider $R=2S$ for a compactly supported function $f$ whose support is $B_0(S)$.\\
                Our goal is to implement the FFT algorithm to evaluate the quadrature approximation, so we will make use of periodic distribution functions, and then we have to properly take into account the aliasing effect due to periods superposition in the Fourier space. A geometrical argument shows that using the periodicity of the function it is enough to take $T\geq (3+\sqrt{2})S/2$ to prevent intersections of the region where $f$ is different from zero. The operator now can be written as 
                \begin{equation}\label{eqQRapSpe}
                    \mathcal{Q}^R(f,f)(v)=
                    \int_{x\in \mathcal{B}_0(R)}
                    \int_{y\in \mathcal{B}_0(R)}
                    \tilde{B}(x,y)\delta(x\cdot y)
                    \left[f(v+y)f(v+x)-f(v+x+y)f(v)\right] \ dx \ dy,
                \end{equation}
                for $v\in \mathcal{D}_T$. The main interest of this representation is preserving the real collision kernel and its properties.

                The main idea to develop a more efficient algorithm is to use the representation \eqref{eqQRapSpe}, where taking $\phi(v)=e^{ik\cdot v}$ we obtain the following spectral quadrature formula
                \begin{equation}
                    \hat{Q}_k = \sum_{\substack{l,m=-N\\l+m=k}}^N \hat{f}_l \hat{f}_m \hat{\beta}_F(l, m), \ \ \ k = -N, ..., N,
                \end{equation}
                where $\hat{\beta}_F (l,m)=\hat{B}_F(l,m) - \hat{B}_F(m,m)$ are now given by
                \begin{equation}
                    \hat{B}_F(l,m) =
                    \int_{\mathcal{B}_0(R)}
                    \int_{\mathcal{B}_0(R)}
                    \tilde{B}(x,y)\delta(x\cdot y)
                    e^{i(l\cdot x + m \cdot y)} \ dx \ dy.
                \end{equation}
                Finding a convolution structure in the previous equation  reduces the number of operations needed to evaluate it. In particular, we want to approximate each $\hat{\beta}_F(l,m)$ as follows
                \begin{equation}\label{BSpectral}
                    \hat{\beta}_F(l,m) \simeq \sum_{p=1}^A \alpha_p(l)\alpha_p'(m),
                \end{equation}
                where $A$ represents the number of finite possible directions of collisions. In the end, we will obtain a sum of $A$ discrete convolutions, which can be computed in $\mathcal{O}(AN\log_2{N})$ operations using the standard FFT technique.\\
                Now we make the decoupling assumption
                \begin{equation}
                    \tilde{B}(x,y)=a(|x|)b(|y|),
                \end{equation}
                in order to obtain the convolution form. We would like to remark that this assumption is satisfied when $\tilde{B}$ is constant.\\
                This is the case of the Maxwellian molecules in dimension two, and hard-spheres in dimension three, the two most used cases. Indeed, putting the variable hard-sphere kernel inside the equation \eqref{BetaTildeCarleman} we obtain
                \begin{equation}
                    \tilde{B}(x,y)=2^{d-1}C_\gamma (|x|^2 + |y|^2)^{-\frac{d-\gamma-2}{2}},
                \end{equation}
                and so we obtain that $\tilde{B}$ is constant if $d=2$, $\gamma=0$ (Maxwellian molecules in dimension two) or if $d=3$, $\gamma=1$ (hard-sphere molecules in dimension three).

                \paragraph{\textit{The 2d case.}}
                We start with Maxwelliam molecules in dimension two, assuming $\tilde{B}=1$. We move to spherical coordinate $x=\rho e$ and $x'=\rho' e'$ inside \eqref{BSpectral}, and we get
                \begin{equation}
                    \tilde{B}_F(l,m)=\frac{1}{4}\int_{\SSS^1}\int_{\SSS^1}
                    \delta(e\cdot e')
                    \left[\int_{-R}^Re^{i\rho(l\cdot e)}d\rho\right]
                    \left[\int_{-R}^Re^{i\rho'(m\cdot e')}d\rho'\right]
                    de \ de'.
                \end{equation}
                Now we introduce $\phi_R^2(s)=\int_{-R}^R e^{i\rho s} d\rho$, for $s\in \RR$, and we obtain the explicit formula
                \begin{equation}
                    \phi_R^2(s)=2R\sinc{(Rs)},
                \end{equation}
                where $\sinc(x)=(\sin(x)/x)$.
                Using the parity properties of the $\sinc$ function inside the expression of $\hat{B}_F(l,m)$ we get 
                \begin{equation}
                    \hat{B}_F(l,m)=
                    \int_0^\pi \phi_R^2(l\cdot e_\theta) 
                    \phi_R^2(m\cdot e_{\theta + \pi/2}) \ d\theta.
                \end{equation}
                Finally a regular discretization of $A$ equally spaced points $\theta_p=\pi p/A$ of the unit sphere, which is spectrally accurate \cite{can-hus-qua-zan-1991} because of the periodicity of the function, gives
                \begin{equation}
                    \hat{B}_F(l,m)=\frac{\pi}{A}
                    \sum_{p=1}^A \alpha_p(l)\alpha_p'(m),
                \end{equation}
                with
                \begin{equation}
                    \alpha_p(l)=\phi_R^2(l\cdot e_{\theta_p}), \ \ \ \ \ \ \ \alpha_p'(m)=\phi_R^2(m\cdot e_{\theta_{p}+\pi/2}).
                \end{equation}

                \paragraph{\textit{The 3d case.}}
                Now we deal with the hard-sphere collision kernel in the three dimensional case. First of all we consider a spherical parametrization $(\theta, \phi)$ of the vector $e \in \SSS^2_+$ and a uniform grid of respective size $A_1$ and $A_2$ for the variables $\theta$ and $\phi$ (spectrally accurate thanks to the periodicity of the function) which leads to the following quadrature for $\tilde{B}_F(l,m)$ \cite{rey-2023}
                \begin{equation}
                    \tilde{B}_F(l,m)=\frac{\pi^2}{A_1 A_2}
                    \sum_{p,q=0}^{A_1, A_2}
                    \alpha_{p,q}(l)\alpha_{p,q}'(m),
                \end{equation}
                where
                \begin{equation}
                    \alpha_{p,q}(l)=\phi_{R,a}^3\left(l\cdot e_{(\theta_p, \phi_q)}\right), \ \ \ \ \ 
                    \alpha_{p,q}'(m)=\psi_{R,b}^3\left(\Pi_{e_{\theta_p, \phi_q}^\perp}(m)\right),
                \end{equation}
                \begin{equation}
                    \phi_{R,a}^3(s)=\int_{-R}^R \rho a(\rho)e^{i\rho s}d\rho, \ \ \ \ \ 
                    \psi_{R,b}^3(s)=\int_{0}^\pi \sin{(\theta)}\phi_{R,b}^3(s \cos{(\theta)})d\theta,
                \end{equation}
                and for all $p$ and $q$
                \begin{equation}
                    (\theta_p, \phi_q)=\left(\frac{p\pi}{A_1}, \frac{q\pi}{A_2}\right).
                \end{equation}

                % The numerical method used is the same as the one implemented in the previous section for the BGK operator, with the only difference that now we have that the collision operator is numerically integrated using the spectral method discussed in Section \ref{spectralMethods}.\\
			  
    \subsection{Evolving the solution and coupling the models}
      \label{subEvoCouple}  
      
		  We are now interested in evolving in time the domain decomposition scheme. At a given time $t^n$, we denote by $K_{ij}$ a control volume, the space domain $\Omega = \Omega_0^n \sqcup \Omega_1^n \sqcup \Omega_2^n$ is decomposed in
		  \begin{itemize}
		    \item \emph{Fluid cells} $K_{ij} \subset \Omega_0^n$, described by the hydrodynamic fields
		      \[
		        U_{ij}^n := \left (\rho_{ij}^n, \bm{u}_{ij}^n, T_{ij}^n\right ) \simeq \left (\rho(t^n, x_i, y_j), \bm{u}(t^n, x_i, y_j), T(t^n, x_i, y_j)\right ) ; 
		      \]
		    \item \emph{Kinetic cells evolved with ES-BGK equation} $K_{ij} \subset \Omega_1^n$, described by the particle distribution function $f^1$ evolved using the ES-BGK equation
		    \[
		      f_{ij}^{1,n}(v) \simeq f^1(t^n, x_i,y_j,v), \quad \forall v \in \RR^{3}.
		    \]
              \item \emph{Kinetic cells evolved with Boltzmann equation} $K_{ij} \subset \Omega_2^n$, described by the particle distribution function $f^2$ evolved using the Boltzmann operator by mean of the spectral method
        		    \[
        		      f_{ij}^{2,n}(v) \simeq f^2(t^n,x_i,y_j,v), \quad \forall v \in \RR^{3}.
        		    \]
		  \end{itemize}	

             Now let's introduce some important quantities. Since we are in a two dimensional space, the matrix $\mathcal{V}_{\ve^1}$ is equal to
            \begin{equation}
                \mathcal{V}_{\ve^1}=
                \begin{pmatrix}
                    1 - \xi_1 a_1 - \xi_2 a_2 &   - \xi_1 b_1 - \xi_2 b_2 & 0,            \\
                      - \xi_1 b_1 - \xi_2 b_2 & 1 - \xi_1 c_1 - \xi_2 c_2 & 0,            \\
                    0                         & 0                         & 1 - \xi_1 h_1,\\
                \end{pmatrix}
            \end{equation}
            where $a_1=\frac{4}{3}\frac{\partial u_x}{\partial x} - \frac{2}{3}\frac{\partial u_y}{\partial y}$, $a_2=(\frac{\partial T}{\partial x})^2$, $b_1 = \frac{\partial u_x}{\partial y}+\frac{\partial u_y}{\partial x}$, $b_2=(\frac{\partial T}{\partial x}\frac{\partial T}{\partial y})$, $c_1 = \frac{4}{3}\frac{\partial u_y}{\partial y} - \frac{2}{3}\frac{\partial u_x}{\partial x}$, $c_2=(\frac{\partial T}{\partial y})^2$, $h_1=-\frac{2}{3}\left(\frac{\partial u_x}{\partial x}+\frac{\partial u_y}{\partial y}\right)$, $\xi_1 = \frac{\epsilon \mu}{\rho T}$ and $\xi_2 = \epsilon^2\frac{2k^2}{3\rho^2 T^3}$.\\
            The three eigenvalues are then:
            \begin{equation}
                \begin{cases}
                    \lambda_{\ve^1}^a = \frac{1}{2}\left(\sqrt{\left[\xi_1\left(a_1-c_1\right)+\xi_2\left(a_2-c_2\right)\right]^2 + 4 d^2} - a_1\xi_1-a_2\xi_2-c_1\xi_1-c_2\xi_2+2 \right),\\
                    \lambda_{\ve^1}^b = \frac{1}{2}\left(-\sqrt{\left[\xi_1\left(a_1-c_1\right)+\xi_2\left(a_2-c_2\right)\right]^2 + 4 d^2} - a_1\xi_1-a_2\xi_2-c_1\xi_1-c_2\xi_2+2 \right),\\
                    \lambda_{\ve^1}^c = 1 - \xi_1 h_1
                \end{cases}                
            \end{equation}
            where $d=\xi_1 b_1 + \xi_2 b_2$.\\

\begin{figure}
            \newcommand{\varA}{0.9}
                \begin{tikzpicture}
                    \draw(-\varA,-\varA) rectangle (\varA,\varA);
                    
                    \draw(-\varA,\varA) rectangle (\varA,3*\varA);
                    \draw(-\varA,3*\varA) rectangle (\varA,5*\varA);

                    \draw(-\varA,-\varA) rectangle (\varA,-3*\varA);
                    \draw(-\varA,-3*\varA) rectangle (\varA,-5*\varA);

                    \draw (\varA,-\varA) rectangle (3*\varA,\varA);
                    \draw (3*\varA,-\varA) rectangle (5*\varA,\varA);

                    \draw (-\varA,-\varA) rectangle (-3*\varA,\varA);
                    \draw (-3*\varA,-\varA) rectangle (-5*\varA,\varA);

                    \filldraw[black] (0,0) circle (0pt) node[anchor=south]{$(i,j)$};
                    \filldraw[black] (0,0) circle (0pt) node[anchor=north]{Layer 0};
                    
                    \filldraw[black] (2*\varA,0) circle (0pt) node[anchor=south]{$(i+1,j)$};
                    \filldraw[black] (2*\varA,0) circle (0pt) node[anchor=north]{Layer 1};
                    
                    \filldraw[black] (4*\varA,0) circle (0pt) node[anchor=south]{$(i+2,j)$};
                    \filldraw[black] (4*\varA,0) circle (0pt) node[anchor=north]{Layer 1};
                    
                    \filldraw[black] (-2*\varA,0) circle (0pt) node[anchor=south]{$(i-1,j)$};
                    \filldraw[black] (-2*\varA,0) circle (0pt) node[anchor=north]{Layer 2};
                    
                    \filldraw[black] (-4*\varA,0) circle (0pt) node[anchor=south]{$(i-2,j)$};
                    \filldraw[black] (-4*\varA,0) circle (0pt) node[anchor=north]{Layer 0};
                    
                    \filldraw[black] (0,2*\varA) circle (0pt) node[anchor=south]{$(i,j+1)$};
                    \filldraw[black] (0,2*\varA) circle (0pt) node[anchor=north]{Layer 0};
                    
                    \filldraw[black] (0,4*\varA) circle (0pt) node[anchor=south]{$(i,j+2)$};
                    \filldraw[black] (0,4*\varA) circle (0pt) node[anchor=north]{Layer 2};
                    
                    \filldraw[black] (0,-2*\varA) circle (0pt) node[anchor=south]{$(i,j-1)$};
                    \filldraw[black] (0,-2*\varA) circle (0pt) node[anchor=north]{Layer 2};
                    
                    \filldraw[black] (0,-4*\varA) circle (0pt) node[anchor=south]{$(i,j-2)$};
                    \filldraw[black] (0,-4*\varA) circle (0pt) node[anchor=north]{Layer 1};
                \end{tikzpicture} 
                \caption{Example of schematic view of the stencil for the computation of the first and second spatial derivatives, and for the application of a finite volume method.}
                \label{fig2}
            \end{figure}
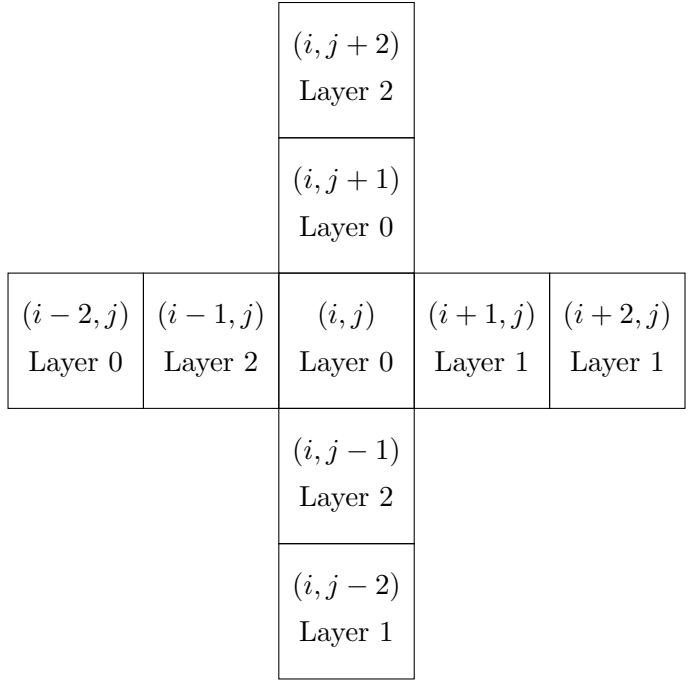

            Now we briefly discuss the computation of the first and second spatial derivative that appear inside the expressions of $\lambda_{\ve^1}^a$, $\lambda_{\ve^1}^b$ and $\lambda_B$ (Equation \eqref{eq1006}), taking into account that in general there are cells close to each other that are of different type (Layer 0, Layer 1 and Layer 2). For computing these quantities it is necessary to consider a stencil of one cell in positive and negative direction for both $x$ and $y$. A schematic example representation is illustrated in Figure \ref{fig2}. We are interested in computing the first and second spatial derivative in each direction ($x$ and $y$) of the following quantities: $\{\rho, u_x, u_y, T\}:=\xi$; and $\nabla_x T$, $\Delta_\mathbf{x} \mathbf{u}$ and $\Delta_\mathbf{x} \rho$.
            Then these quantities have been computed as follows
            \begin{equation}
                \begin{cases}
                    \displaystyle \frac{\partial \xi}{\partial x} = \displaystyle \frac{\xi_{i+1,j}-\xi_{i-1, j}}{2\Delta x},\\[10pt]
                    \displaystyle \frac{\partial \xi}{\partial y} = \displaystyle \frac{\xi_{i,j+1}-\xi_{i,j-1}}{2\Delta y},\\[10pt]
                    \Delta_\mathbf{x} \mathbf{u}=\Delta_\mathbf{x} u_x + \Delta_\mathbf{x} u_y,\\[10pt]
                    \Delta_\mathbf{x} \rho = \displaystyle \frac{\rho_{i-1, j}+\rho_{i+1, j}+\rho_{i, j-1}+\rho_{i, j+1}-4\rho_{i, j}}{\Delta x \Delta y},\\[10pt]
                    \nabla_\mathbf{x} T = \displaystyle \mathbf{\hat{x}}\frac{\partial T}{\partial x} + \mathbf{\hat{y}}\frac{\partial T}{\partial y}.             
                \end{cases}
            \end{equation}
            where
            \begin{equation}
                \begin{cases}
                    \Delta_x u_x= \displaystyle \frac{u_{x, i-1, j}+u_{x, i+1, j}+u_{x, i, j-1}+u_{x, i, j+1}-4u_{x, i, j}}{\Delta x \Delta y},\\[10pt]
                    \Delta_x u_y=\displaystyle \frac{u_{y, i-1, j}+u_{y, i+1, j}+u_{y, i, j-1}+u_{y, i, j+1}-4u_{y, i, j}}{\Delta x \Delta y}.\\[10pt]
                \end{cases}
            \end{equation}
            It is evident that for applying these differential operators to one of the $\xi$ in the cell $(x_i, y_j)$ it is necessary to know the value of $\xi$ in the cells $(x_{i}, y_{j-1})$, $(x_{i}, y_{j+1})$, $(x_{i-1}, y_{j})$ and $(x_{i+1}, y_{j})$. The code properly estimate the value of $\xi$ in each of these cells, according to its Layer. For example, in Figure \ref{fig2} the cell $(x_{i+1}, y_{j})$ is of type Layer 1, so the code estimate $\xi_{i+1,j}$ as a moment associated to the distribution function $f^1(i+1,j)$; while, in the same Figure, the cell $(x_{i}, y_{j+1})$ is of type Layer 0, so the quantity $\xi_{i, j+1}$ is obtained directly from the hydrodynamic solution $U(i,j+1)$.

            The Algorithm \ref{alg1} summarises the switching procedure between different regimes, and the following quantities are used: $\lambda_{\ve^1}^k(i,j)=\lambda_{\ve^1}^k(x_i, y_j)$ for $k=\{a,b,c\}$, indicate the three eigenvalues associated to the matrix $\mathcal{V}_{\ve^1}$; and $\lambda_B(i,j)=\lambda_B(x_i, y_j)$.\\

            % in each cell $(x_i, y_j)$ we can compute the two eigenvalues associated to the matrix $\nu_{NS}$ that we will indicate as $\lambda_{NS}^k(i,j)=\lambda_{NS}^k(x_i, y_j)$ for $k=\{a,b\}$. Analogously we can also compute $\lambda_B(i,j)=\lambda_B(x_i, y_j)$.\\
    
		  \begin{algorithm}
            \caption{Regime update}
            \algsetup{indent=2em}
            \begin{algorithmic}[1]
                \STATE Actual time: $t^n$
                \IF{$r^n_{ij}==2$}
                    \STATE compute $\lambda_B(i,j)$
                    \IF{$|\lambda_B(i,j)| > \eta_1$} 
                        \STATE $r^{n+1}_{ij}=2$
                        \ELSE
                        \STATE $r^{n+1}_{ij}=1$
                    \ENDIF
                \ELSIF{$r^n_{ij}==1$}
                    \STATE compute $\lambda_B(i,j)$, $\lambda_{\ve^1}^a(i,j)$, $\lambda_{\ve^1}^b(i,j)$ and
                    $\lambda_{\ve^1}^c(i,j)$
                    \IF{$|\lambda_B(i,j)| > \eta_1$}
                        \STATE $r^{n+1}_{ij}=2$
                    \ELSE 
                        \IF{$|\lambda_{\ve^1}^a - 1|\leq\eta_0$ and $|\lambda_{\ve^1}^b - 1|\leq\eta_0$ and $|\lambda_{\ve^1}^c - 1|\leq\eta_0$}
                            \STATE estimate $\tilde{U}_{ij}=(\rho, \mathbf{u}, T)$ as moments of $f^{1}(i,j)$
                            \STATE compute $\xi=||f^{1}(i,j)-\mathcal{M}_{\tilde{U}_{ij}}^n||_{L^1_v}$
                            \IF{$\xi\leq \delta_0$}
                                \STATE $r^{n+1}_{ij}=0$
                            \ELSE
                                \STATE $r^{n+1}_{ij}=1$
                            \ENDIF
                        \ELSE
                            \STATE $r^{n+1}_{ij}=1$
                        \ENDIF
                    \ENDIF                   
                \ELSE
                    \STATE compute $\lambda_{\ve^1}^a(i,j)$ and $\lambda_{\ve^1}^b(i,j)$ and $\lambda_{\ve^1}^c(i,j)$
                        \IF{$|\lambda_{\ve^1}^a - 1|>\eta_0$ or $|\lambda_{\ve^1}^b - 1|>\eta_0$ or $|\lambda_{\ve^1}^c - 1|>\eta_0$}
                            \STATE $r^{n+1}_{ij}=1$
                        \ELSE
                            \STATE $r^{n+1}_{ij}=0$
                        \ENDIF          
                \ENDIF
            \end{algorithmic}
            \label{alg1}
            \end{algorithm}

			Now we explain what happens between two cells of different types when performing the spatial integration. Let us consider the situation at time $t^n$ depicted in Figure \ref{fig2}. Since the cell $(x_i, y_j)$ is hydrodynamic (Layer 0), the code uses the CWENO TVD RK2 scheme to integrate it. This method requires a stencil in both $x$ and $y$ direction. As we can see from the same Figure, not all the cells of the stencil are of the correct type (hydrodynamic) Layer 0, so in order to use the numerical method we estimate the value of the hydrodynamic solution $U$ in the other cells using the distribution function $f^2$ for the cells of Layer 2 and the distribution function $f^1$ for the cells of Layer 1. For example, in the cells $(i-1, j)$, $(i+1, j)$ and $(i, j-1)$, the hydrodynamic solution is computed using, respectively the moments of the distribution function $f^2(i-1,j)$, $f^1(i+1,j)$ and $f^2(i, j-1)$. We perform the same operation for all cells constituting the stencil.\\

            If now we consider a small variation of the Figure \ref{fig2}, using the Layer 2 in the cell $(i,j)$, then a consideration similar to previous one is still valid. In particular the method that the code uses to evolve in time the distribution functions $f^1$ and $f^2$, uses as well a stencil in both $x$ and $y$ directions. As we can see from the same Figure, not all the cells of the stencil are of the correct type Boltzmann operator (Layer 2), so in order to use the numerical method we estimate the value of the distribution function $f^2$ in the cells $(i, j+1)$, $(i-2, j)$, using the Chapman-Enskog expansion up to the zeroth order computed with the hydrodynamic solutions $U(i,j+1)$ and $U(i-2,j)$; while in the cells $(i+1,j)$, $(i+2,j)$, $(i,j-2)$, we replace the distribution function as follows $f^2(i+1, j)=f^1(i+1,j)$, $f^2(i+2,j)=f^1(i+2,j)$ and $f^2(i, j-2)=f^1(i, j-2)$. We perform the same operation for all cells constituting the stencil. So after these substitutions we execute the numerical scheme.\\

            It now remains to discuss what happens to the solution in a specific cell of coordinate $(x_i, y_j)$ when its regime changes. We have four possible situations:
            \begin{itemize}
                \item the first situation is a cell that switches from Layer 0 (Euler equations) to Layer 1 (ES-BGK equation). It means that in the previous timestep the cell has been updated using Euler equations, but in the actual timestep it must be updated using the kinetic solver, so we have to associate a distribution function to the cell. In this situations we decide that the distribution function associated to the cell is the Maxwellian (corresponding to 0-th order in the Chapman-Enskog expansion) computed using the moments of the solutions of the Euler equations obtained in the previous timestep using the CWENO TVD RK2 scheme. Better approximation of the distribution function could be obtained using higher order in the Chapman-Enskog expansion (see \cite{fil-rey-2014});
                \item cell switching from Layer 1 (ES-BGK equation) to Layer 0 (Euler equations). In this situation the solution computed in the previous timestep was the distribution function updated with the ES-BGK equation, but now we have to update it using the Euler equations. So we compute the moments of the distribution function and we use them as input for the CWENO TVD RK2 scheme;
                \item cell switching from Layer 1 (ES-BGK equation) to Layer 2 (Boltzmann equation). In this situation we decide to use for the Layer 2 the same distribution function that we obtained using the ES-BGK equation in the previous timestep, but of course, in the actual timestep the distribution function is updated using the Boltzmann operator;
                \item the last situation is a cell switching from Layer 2 (Boltzmann equation) to Layer 1 (ES-BGK equation). Similarly to the previous case, in this configuration we decide to use for the Layer 1 the same distribution function that we updated using the Boltzmann operator in the previous timestep, but of course, in the actual timestep the distribution function is updated using the ES-BGK equation.
            \end{itemize}
			  
\section{Numerical simulations}
	    \label{subNumSim}
	   The following simulations show the stability of the method implemented and the speed-up offered by the hybrid domain decomposition method compared to using only the Boltzmann operator. The solution computed with the hybrid scheme is compared with the two solutions obtained using, respectively, only the Euler equations and only the Boltzmann operator. The hybrid code (and the full Euler one) is written using C++20 ISO, it is executed only on CPU, and it takes advantage of the multithreads capabilities: in particular in order to speed up the execution, the space domain is divided in subdomains, each of which is assigned to a specific thread to apply the numerical integration method, and all threads are executed in parallel. All the simulations have been executed on a non-dedicated computer equipped with 
       %with a \emph{12th Gen Intel(R) Core (TM) i7-12700H 2.70 GHz} CPU, \emph{16.0 GB} of RAM and a \emph{NVIDIA(R) GeForce RTX-3060 Laptop} GPU for the first device, and
       a \emph{Intel (R) Core (TM) Ultra 7 155H} CPU, \emph{16.0 GB} of RAM and a \emph{NVIDIA(R) RTX 1000 Ada Generation Laptop} GPU. This shows the great scalability, and numerical optimization of the code. The CPU FFT and IFFT numerical operations required by the Spectral Methods for the Boltzmann operator have been performed using the FFTW library. Since some operations of this library are not multithreads safe, the \emph{semaphores} \cite{frigo1999fftw} have been implemented in order to control the access to shared resources, while the respective GPU operations have been implemented using the cuFFT library. 
       
       The solutions computed using only the Boltzmann operator have been obtained using CUDA Programming Language, by \emph{massively parallelizing} the software. This allows us to highly speed up the execution taking advantage of the intrinsic multithreads capability of the GPUs.
       Of course, the computation of the kinetic solution using GPUs is faster than the computation of the hybrid solution using CPUs. So in order to show the speed-up offered by the hybrid-CPU solution compared to the full kinetic-CPU solution, we have measured the execution time of a few timesteps of the full kinetic-CPU solution, and then we have estimated the total duration of this integration by multiplying this quantity by the ratio between the total number of time steps and the already computed time steps.
       
       In all the simulations we used the simplest boundary condition implementing the ``ghost cells": on the four rows and four columns closest to the border, we copy the value of the solution and the domain typology of the closest physical cell. In all simulations the velocity space of the distribution functions has 3 dimensions.\\

       In all simulations, except for the last one (\hyperref[subNonUniformKnudsen]{\textbf{Test 5}}), the regime is initialized everywhere as Euler equations (Layer 0), and then we evolve the system according to our scheme.\\
       
       In all our simulations we used a collision frequency for the ES-BGK operator equals to $\nu_{ij}^n=\nu(\rho_{ij}^n, T_{ij}^n)=\frac{\pi}{2} \rho_{ij}^n$.

        \subsection{Test 1: Sod shock tube}
        \label{subTestSod}
\begin{figure}
    \centering
    \begin{subfigure}[b]{0.32\textwidth}
        \centering
        \includegraphics[width=\textwidth, trim={0.9cm 0cm 0.9cm 0cm}]{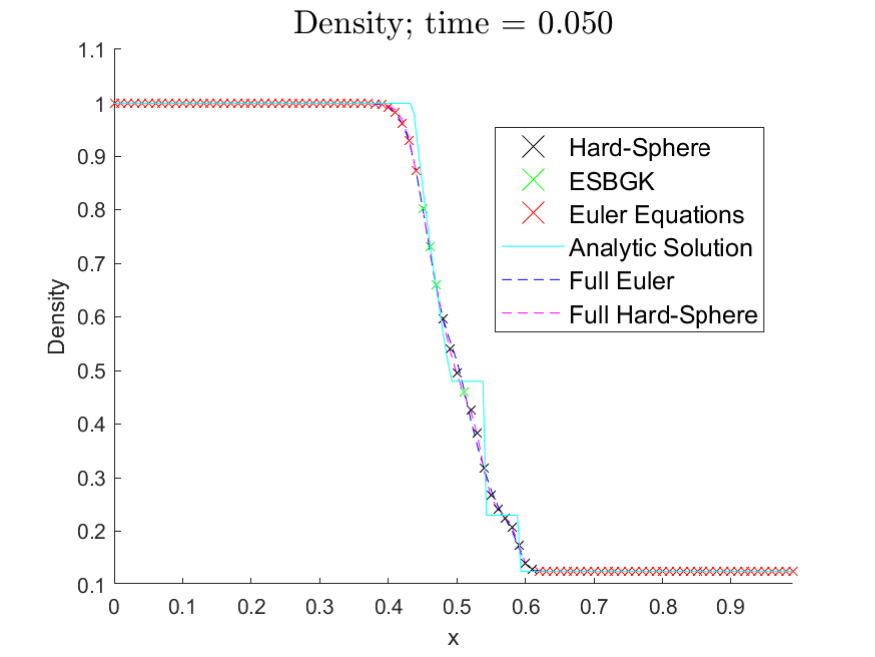}
    \end{subfigure}
    \hfill
    \begin{subfigure}[b]{0.32\textwidth}
        \centering
        \includegraphics[width=\textwidth, trim={0.9cm 0cm 0.9cm 0cm}]{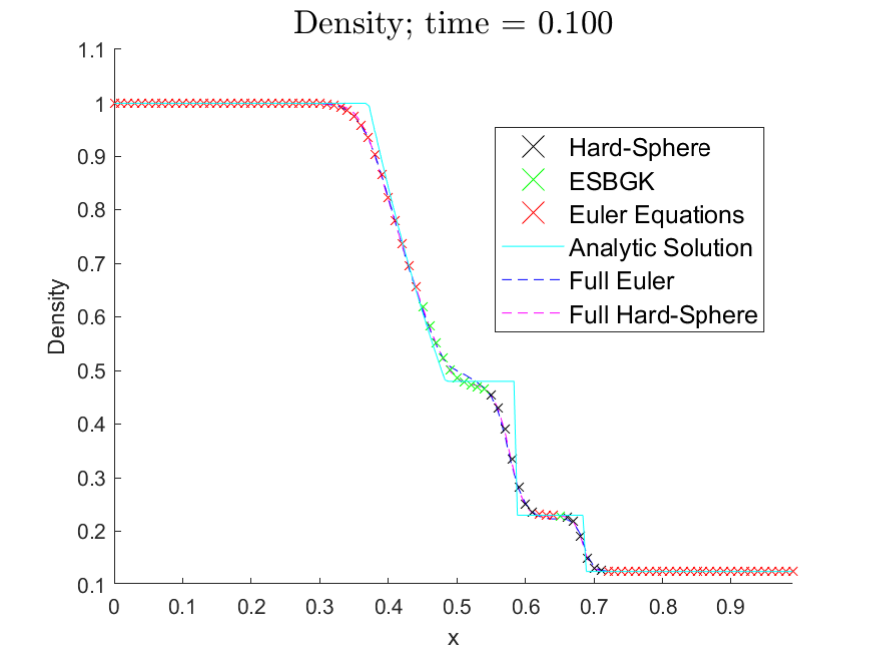}
    \end{subfigure}
    \hfill
    \begin{subfigure}[b]{0.32\textwidth}
        \centering
        \includegraphics[width=\textwidth, trim={0.9cm 0cm 0.9cm 0cm}]{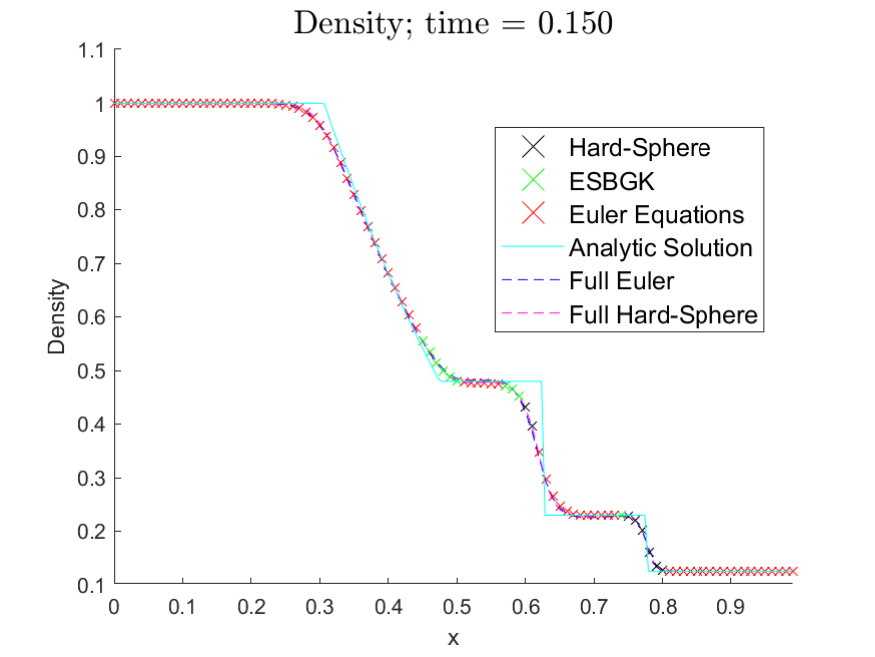}
    \end{subfigure}

    \begin{subfigure}[b]{0.32\textwidth}
        \centering
        \includegraphics[width=\textwidth, trim={0.9cm 0cm 0.9cm 0cm}]{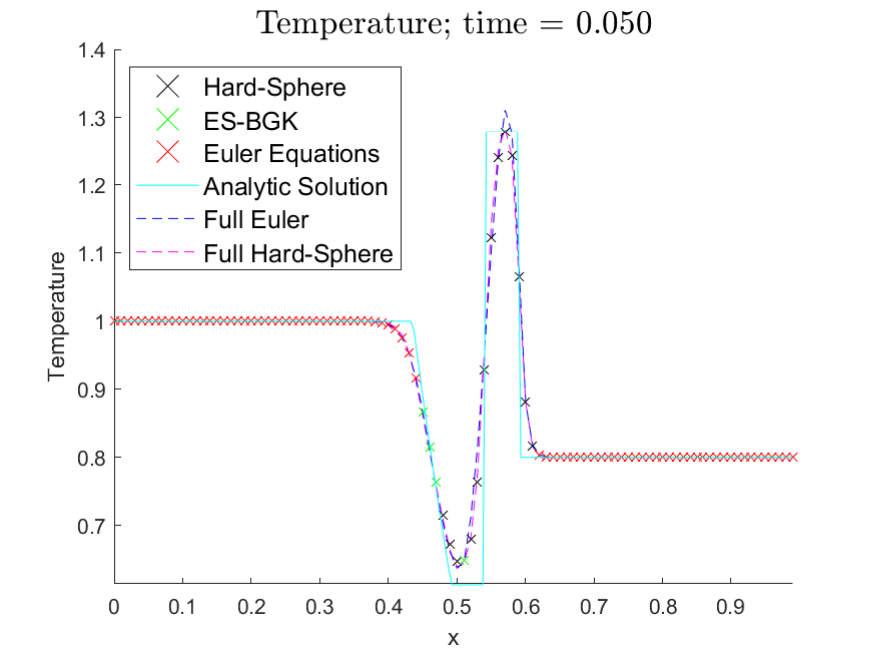}
    \end{subfigure}
    \hfill
    \begin{subfigure}[b]{0.32\textwidth}
        \centering
        \includegraphics[width=\textwidth, trim={0.9cm 0cm 0.9cm 0cm}]{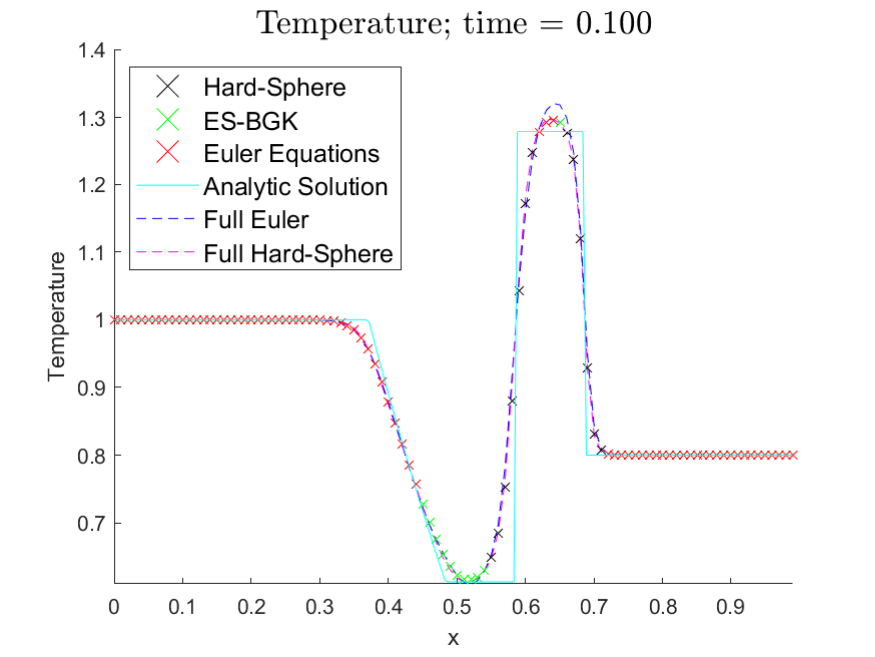}
    \end{subfigure}
    \hfill
    \begin{subfigure}[b]{0.32\textwidth}
        \centering
        \includegraphics[width=\textwidth, trim={0.9cm 0cm 0.9cm 0cm}]{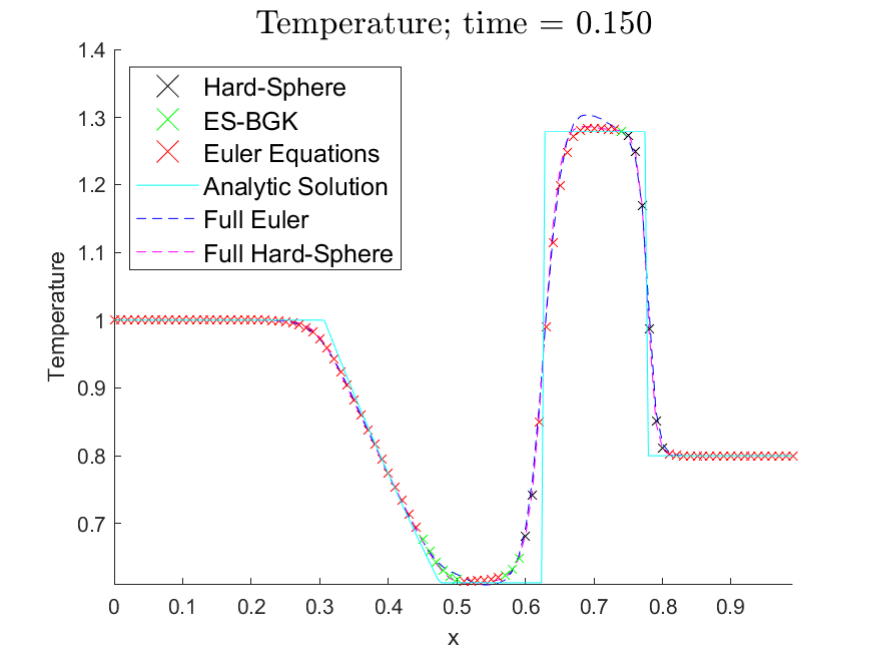}
    \end{subfigure}

    \begin{subfigure}[b]{0.32\textwidth}
        \centering
        \includegraphics[width=\textwidth, trim={0.9cm 0cm 0.9cm 0cm}]{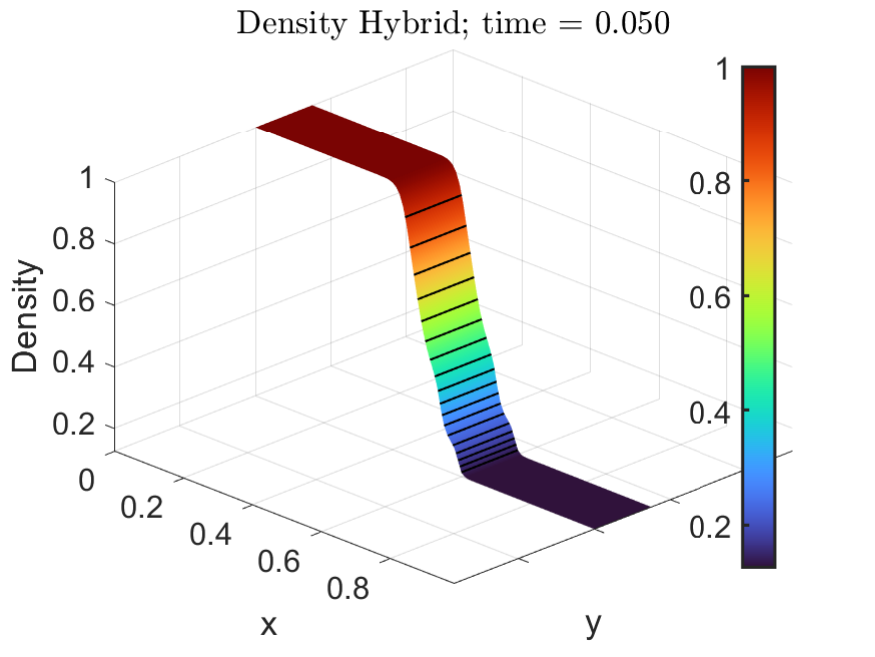}
    \end{subfigure}
    \hfill
    \begin{subfigure}[b]{0.32\textwidth}
        \centering
        \includegraphics[width=\textwidth, trim={0.9cm 0cm 0.9cm 0cm}]{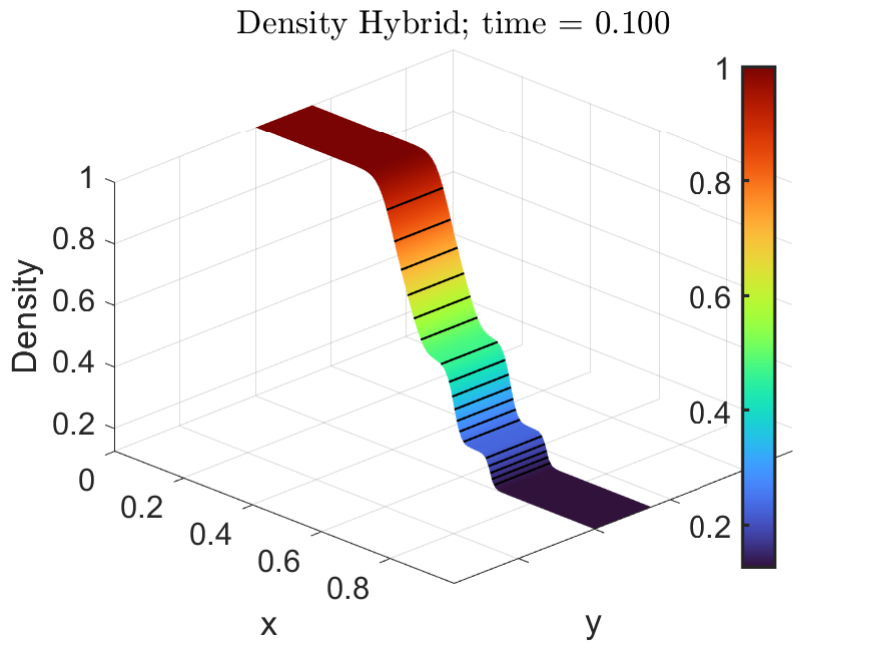}
    \end{subfigure}
    \hfill
    \begin{subfigure}[b]{0.32\textwidth}
        \centering
        \includegraphics[width=\textwidth, trim={0.9cm 0cm 0.9cm 0cm}]{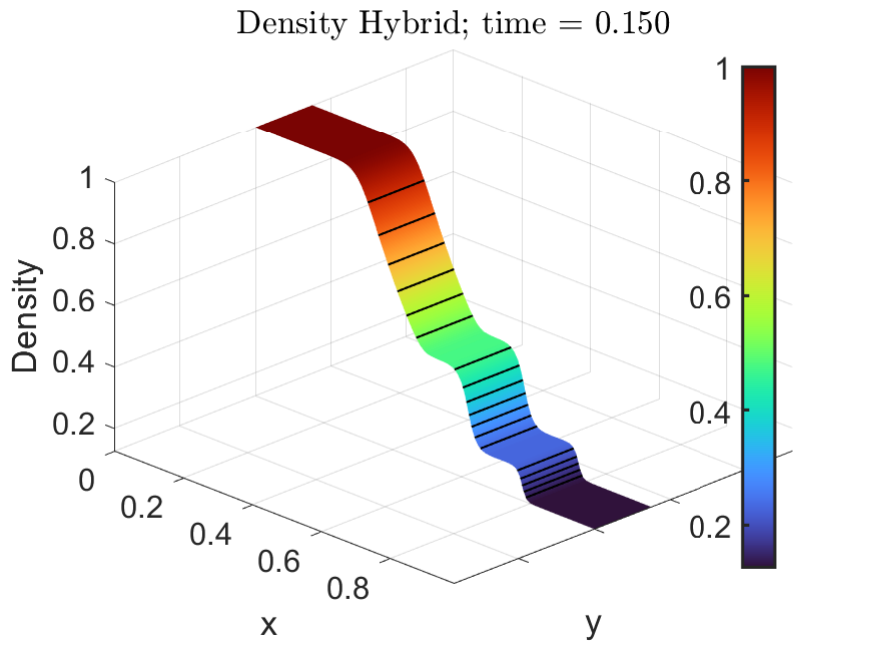}
    \end{subfigure}
    
    \caption{Time evolution of density for \hyperref[subTestSod]{\textbf{Test 1}}. In the first row the regime adaptation and the density at three different times are displayed. Red, green and black dots, identify, respectively, cells updated using Euler equations, Boltzmann equations with ES-BGK operator and Boltzmann equations with Boltzmann (Hard-Sphere) operator, using the new hybrid scheme. The cyan line is the analytical solution, the magenta one is the solution obtained using only the Euler equations, and the blue line is computed using only the Boltzmann equation. The same plots, but for the temperature, are showed in the second row. The time evolution of the density for the hybrid scheme in the 2D space (solution constant along the $y-$direction) is displayed in the last row.}
    \label{fig_density_sim1}
\end{figure}

        A classical benchmark test for the 1D Euler equations is the so-called \emph{Sod shock tube test problem} \cite{tor-2009}: its solution consists of a left rarefaction, a contact and a right shock. In the 1D case, on a spatial domain $[0, 1]$ the initial condition is constituted by a left state $U_L=(\rho_L, u_{xL}, P_L)=(1.0, 0.0, 1.0)$ for $x\leq 0.5$ and by a right state $U_R=(\rho_R, u_{xR}, P_R)=(0.125, 0.0, 0.1)$ for $x>0.5$. In this test we simulate these initial conditions, considering also the $y-$direction, along which the solution is always uniform in space. Also the velocity along the $y-$direction is null. In Figure \ref{fig_density_sim1} we compare the density and the temperature obtained using the hybrid scheme, with the Euler solution, with the Hard-Sphere solution and the analytic one (associated to the Euler equations). We discretize the spatial domain with $100$ points in the $x-$ direction, and $16$ in the $y-$ direction, imposing $\Delta x = \Delta y$. The velocity domain is discretized using 32 points in the three directions, and we considered a cut-off domain $\Omega_v=[-8, 8]^3$. The value of $R$ has been taken approximately equal to 3.3, and both numbers of points $A_1$ and $A_2$ discretizing the spherical integration of the Boltzmann operator (see Section \ref{subBoltzOp}) have been taken equal to 4.
        
        The distribution function is initialized everywhere as a Maxwellian, whose moments, correspond to the ones of the Euler equations (i.e. $U_L$ and $U_R$), and null velocity in the other directions.
        We consider a Knudsen number $\epsilon=10^{-6}$ and a time step $\Delta t=0.1 dx$. It is remarkable that the usage of an IMEX approach for the integration of the ES-BGK equation and of the full Boltzmann one allowed us to use the timestep associated to the transport part of the Boltzmann equation, instead of the one constrained by the collision operator, which in that case, would be much smaller than the one implemented thanks to the IMEX approach. This property significantly accelerates the code execution.
        
        The thresholds used are: $\eta_0=1\times 10^{-5}$, $\eta_1=3.5\times 10^{-10}$, $\delta_0=10^{-3}$.
        
        This test allows us to check that the coupling and the transitions between different regimes are executed correctly and that all the numerical schemes are working properly and converging towards the right solution. Indeed, since this test is performed in a very stiff regime (i.e. $\epsilon << 1$) then the solution of the kinetic levels must converge towards the hydrodynamic limit (Euler equations).
        It is clear that the numerical solution is strongly consistent with the analytical one.
        For this specific test case, the hybrid solution is approximately 1.8 times faster then the full kinetic one.

\subsection{Test 2: Fluid flux towards a fixed rectangular obstacle}
\label{subTestFluidFixedRectObs}
One of the aims of this test is to show the capability of the scheme to handle more complex geometries and real world phenomena. In particular we use our code to simulate the flow around a fixed object with rectangular shape, which is embedded in a computational domain of size $L_x=L_y=1.0$ discretized with 100 points in each direction. We have discretized the velocity domain using 16 points in the three directions, and we considered a cut-off domain $\Omega_v=[-8, 8]^3$. The value of $R$ has been taken approximately equal to 3.3, and both numbers of points $A_1$ and $A_2$ discretizing the spherical integration of the Boltzmann operator (see Section \ref{subBoltzOp}) have been taken equal to 4.

\begin{figure}
    \centering
    \begin{subfigure}[b]{0.32\textwidth}
        \centering
        \includegraphics[width=\textwidth, trim={0.9cm 0cm 0.9cm 0cm}]{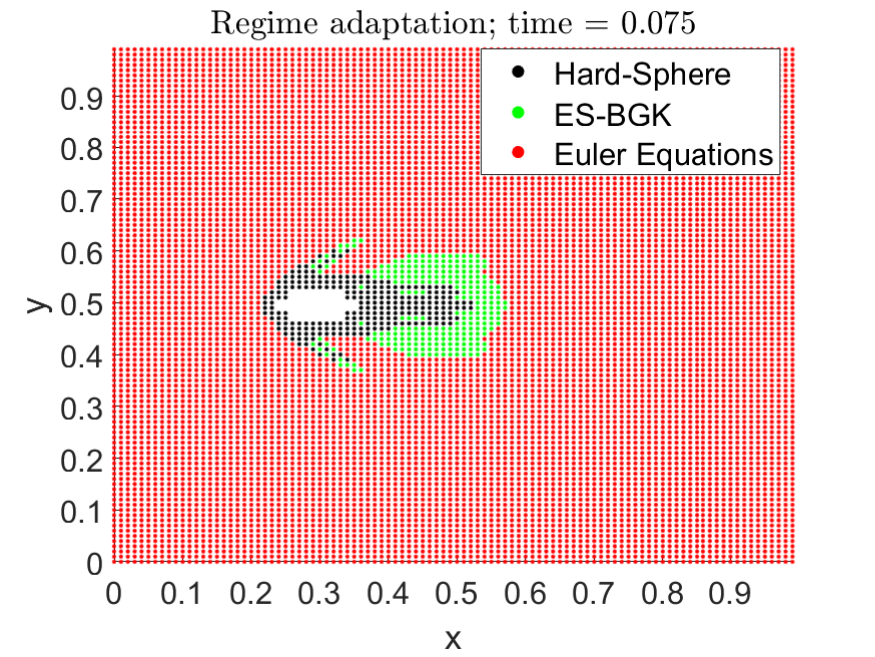}
    \end{subfigure}
    \hfill
    \begin{subfigure}[b]{0.32\textwidth}
        \centering
        \includegraphics[width=\textwidth, trim={0.9cm 0cm 0.9cm 0cm}]{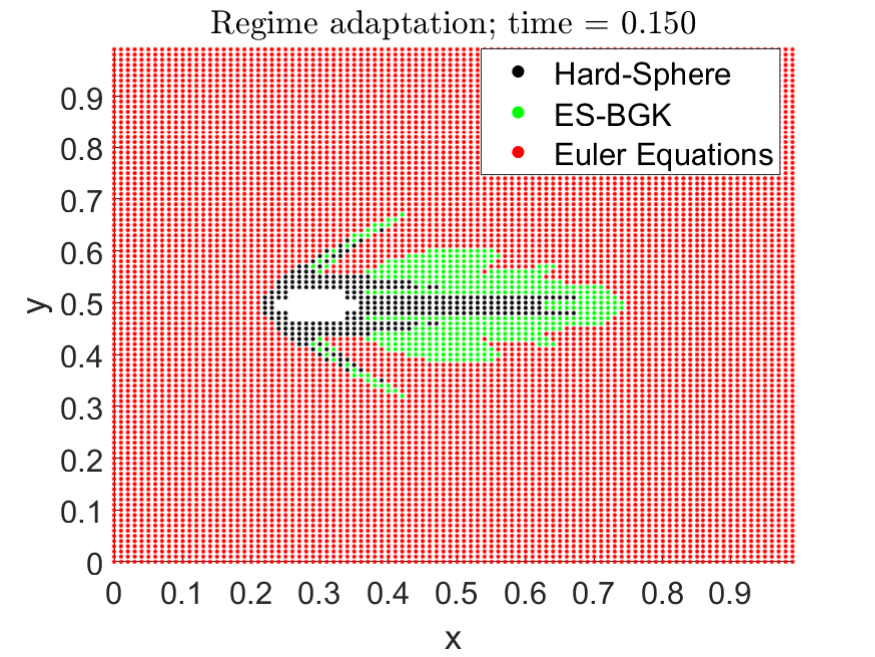}
    \end{subfigure}
    \hfill
    \begin{subfigure}[b]{0.32\textwidth}
        \centering
        \includegraphics[width=\textwidth, trim={0.9cm 0cm 0.9cm 0cm}]{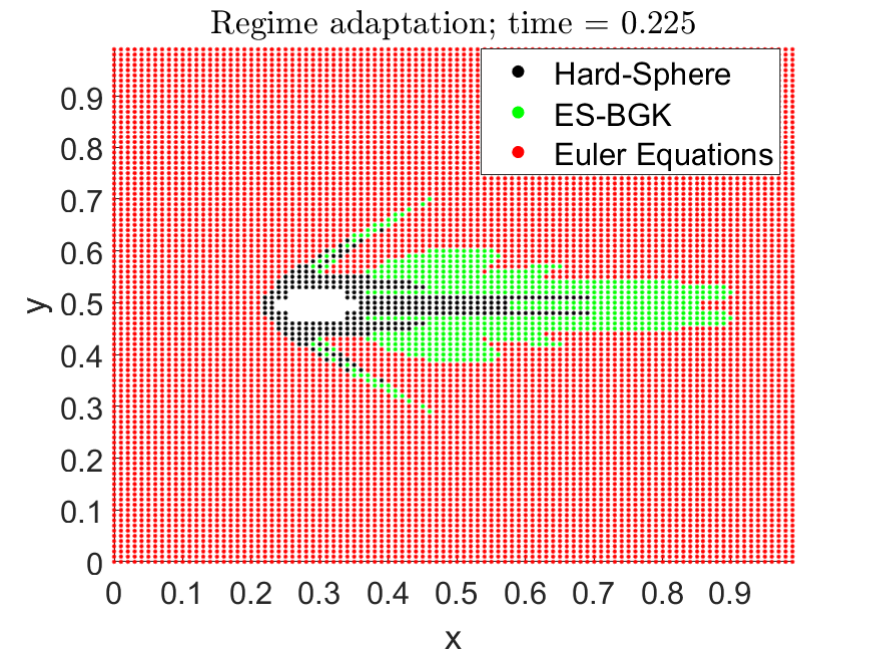}
    \end{subfigure}

    \begin{subfigure}[b]{0.32\textwidth}
        \centering
        \includegraphics[width=\textwidth, trim={0.9cm 0cm 0.9cm 0cm}]{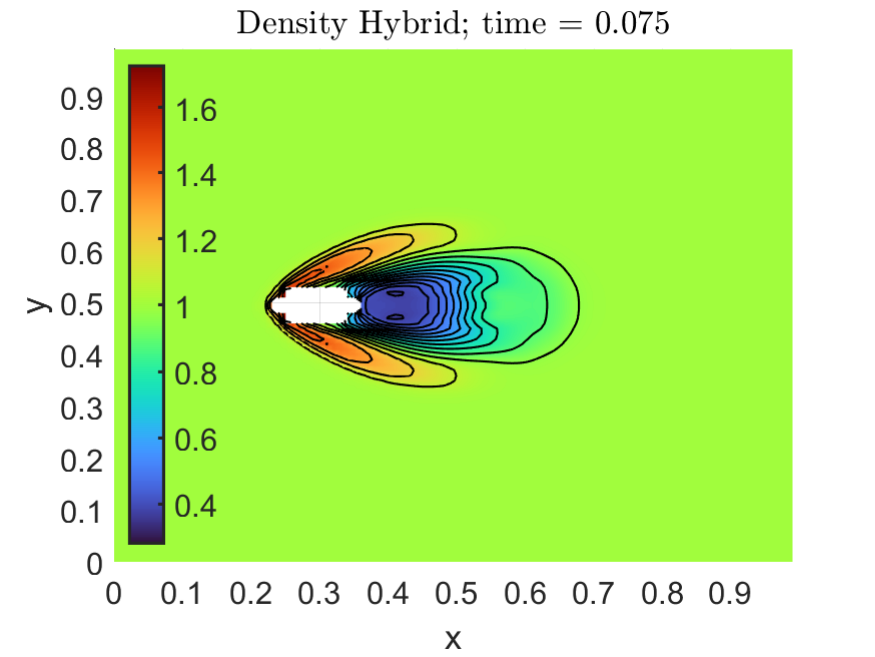}
    \end{subfigure}
    \hfill
    \begin{subfigure}[b]{0.32\textwidth}
        \centering
        \includegraphics[width=\textwidth, trim={0.9cm 0cm 0.9cm 0cm}]{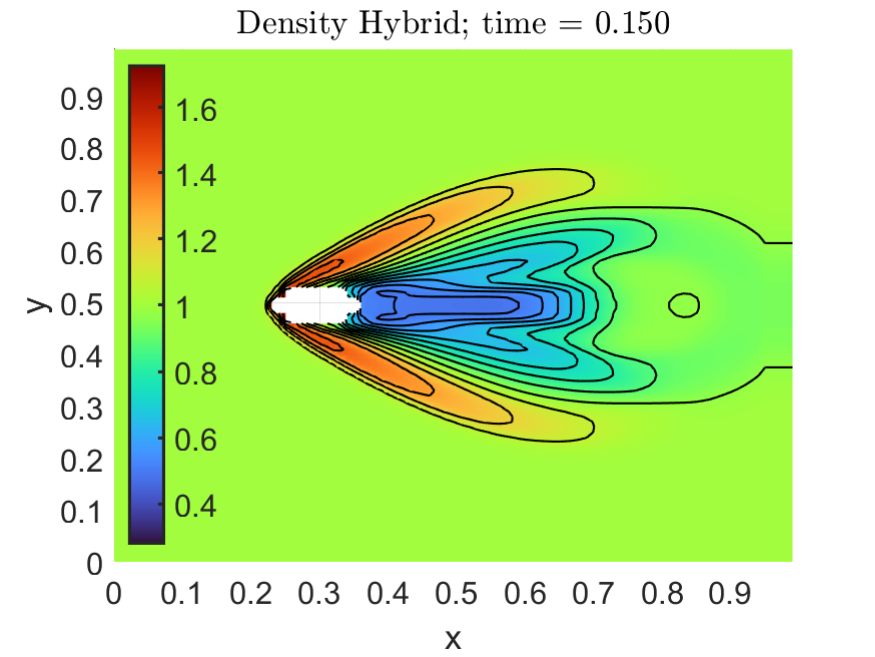}
    \end{subfigure}
    \hfill
    \begin{subfigure}[b]{0.32\textwidth}
        \centering
        \includegraphics[width=\textwidth, trim={0.9cm 0cm 0.9cm 0cm}]{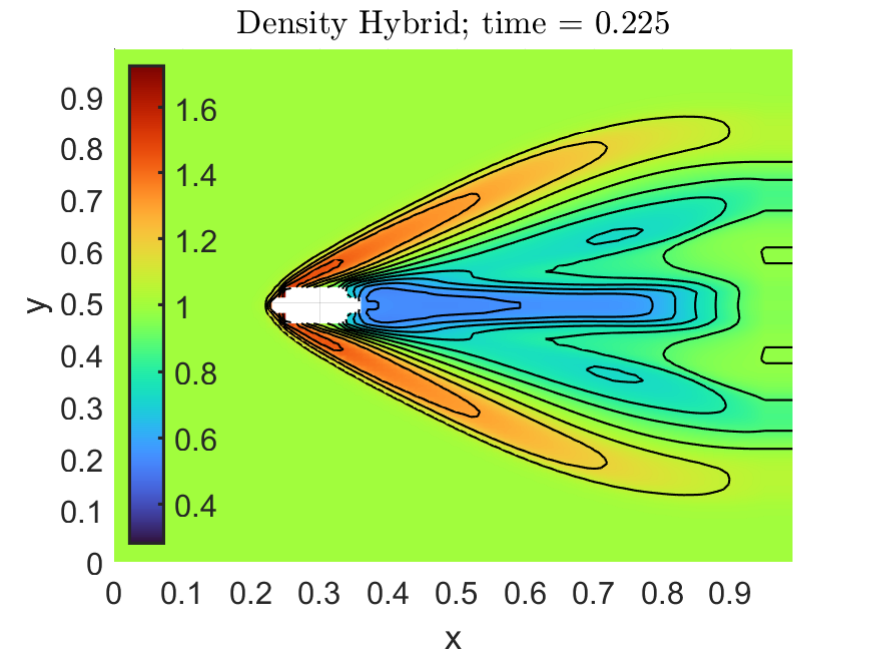}
    \end{subfigure}

    \begin{subfigure}[b]{0.32\textwidth}
        \centering
        \includegraphics[width=\textwidth, trim={0.9cm 0cm 0.9cm 0cm}]{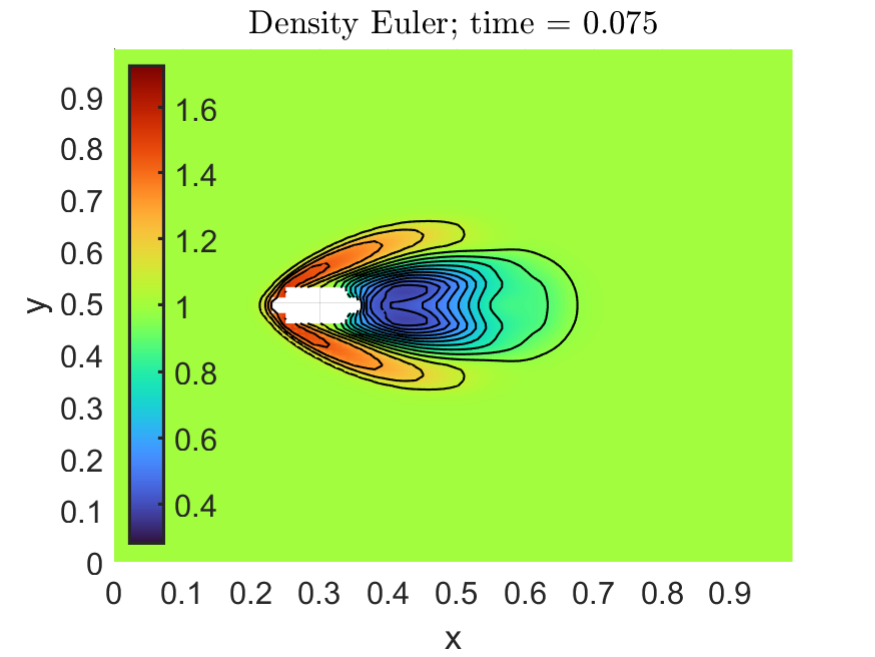}
    \end{subfigure}
    \hfill
    \begin{subfigure}[b]{0.32\textwidth}
        \centering
        \includegraphics[width=\textwidth, trim={0.9cm 0cm 0.9cm 0cm}]{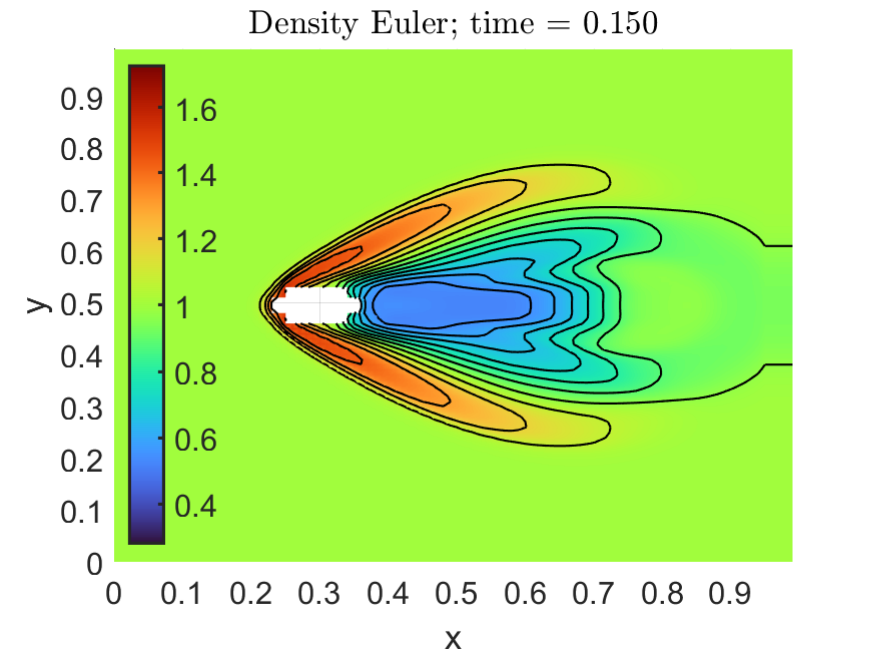}
    \end{subfigure}
    \hfill
    \begin{subfigure}[b]{0.32\textwidth}
        \centering
        \includegraphics[width=\textwidth, trim={0.9cm 0cm 0.9cm 0cm}]{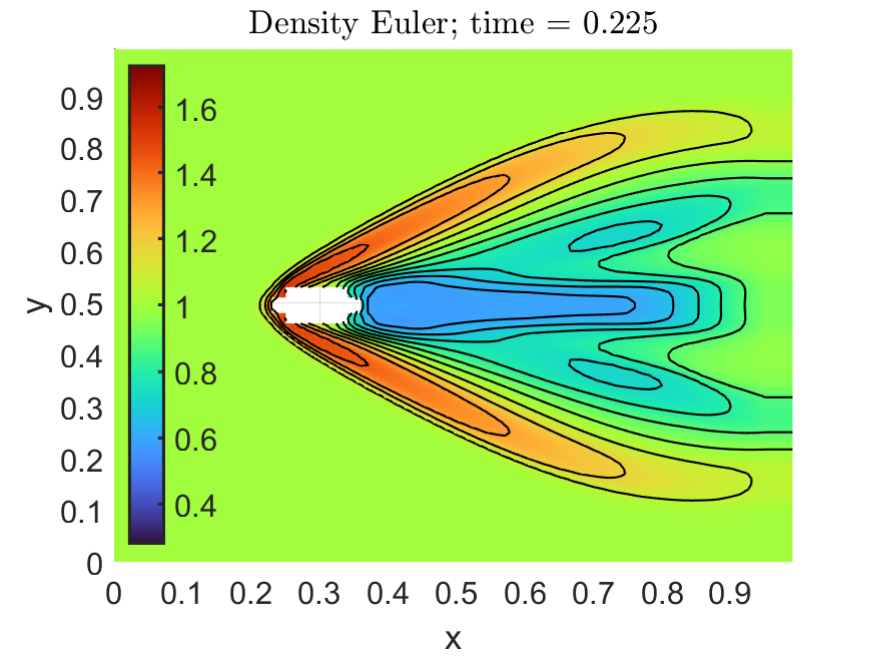}
    \end{subfigure}

    \begin{subfigure}[b]{0.32\textwidth}
        \centering
        \includegraphics[width=\textwidth, trim={0.9cm 0cm 0.9cm 0cm}]{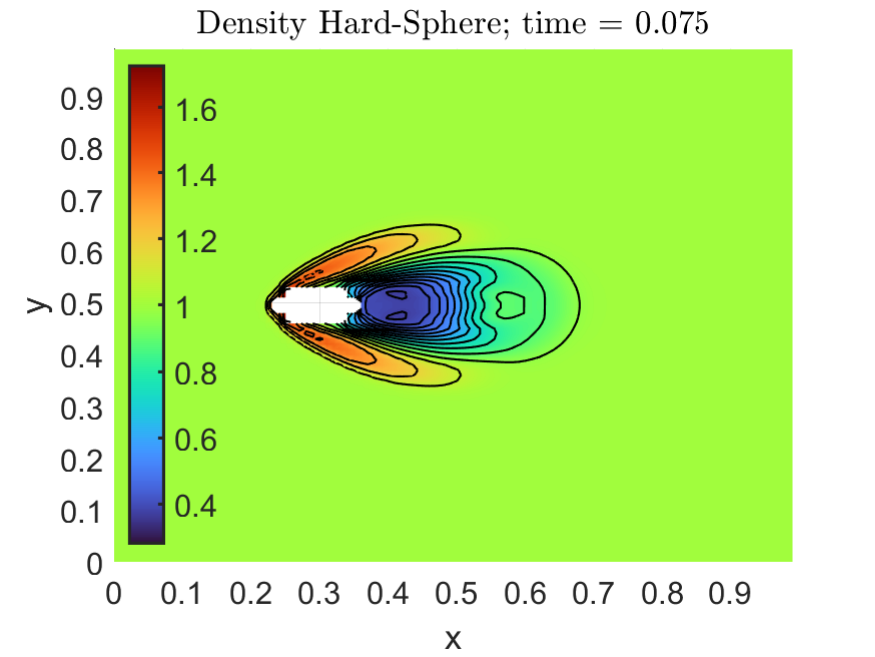}
    \end{subfigure}
    \hfill
    \begin{subfigure}[b]{0.32\textwidth}
        \centering
        \includegraphics[width=\textwidth, trim={0.9cm 0cm 0.9cm 0cm}]{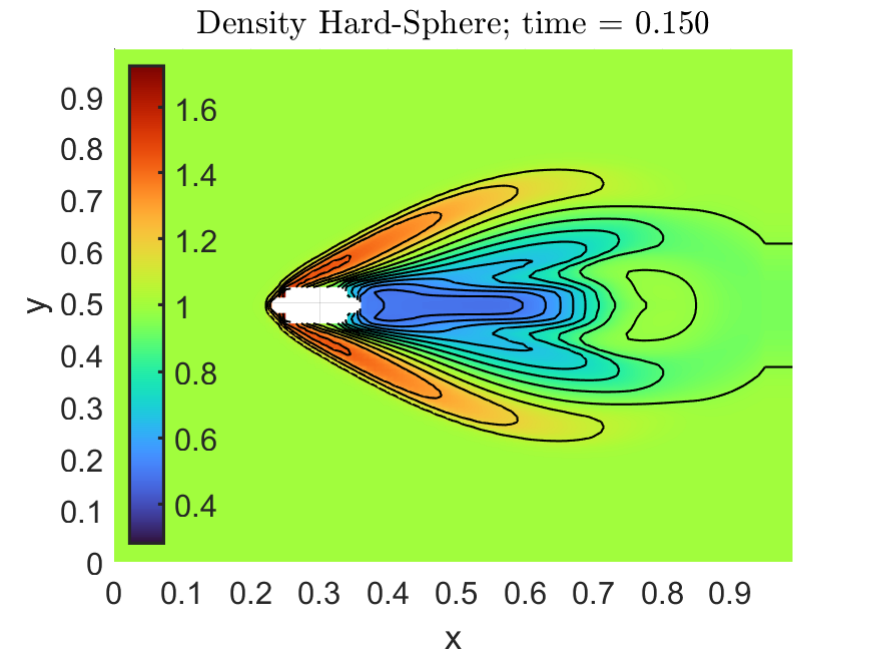}
    \end{subfigure}
    \hfill
    \begin{subfigure}[b]{0.32\textwidth}
        \centering
        \includegraphics[width=\textwidth, trim={0.9cm 0cm 0.9cm 0cm}]{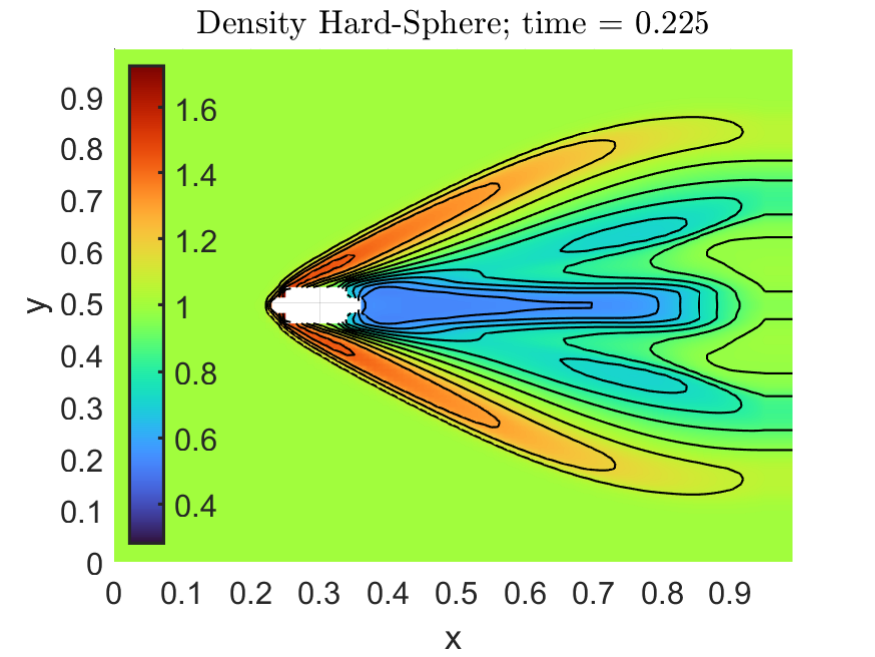}
    \end{subfigure}
    
    \caption{Time evolution of density for \hyperref[subTestFluidFixedRectObs]{\textbf{Test 2}}.  The first row displays the regime adaptation at three different times. Red, green and black dots, identify, respectively, cells updated using Euler equations, ES-BGK equation and Boltzmann equation. In the second, third and fourth row are displayed, respectively, the time evolution of density computed using hybrid scheme, Euler equations and Boltzmann equation.}
    \label{fig_density_sim2}
\end{figure}

\begin{figure}
    \begin{subfigure}[b]{0.32\textwidth}
        \centering
        \includegraphics[width=\textwidth, trim={0.9cm 0cm 0.9cm 0cm}]{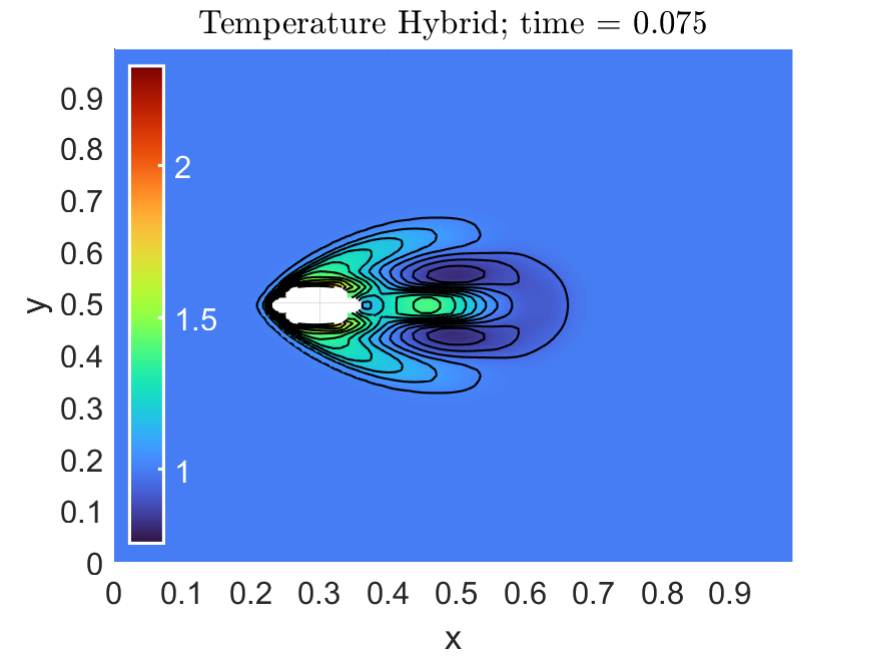}
    \end{subfigure}
    \hfill
    \begin{subfigure}[b]{0.32\textwidth}
        \centering
        \includegraphics[width=\textwidth, trim={0.9cm 0cm 0.9cm 0cm}]{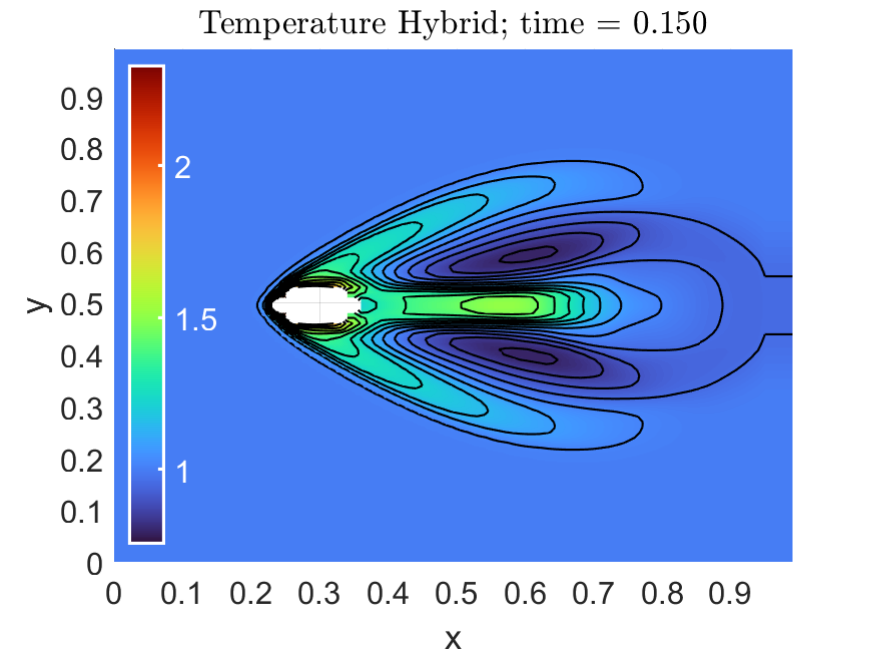}
    \end{subfigure}
    \hfill
    \begin{subfigure}[b]{0.32\textwidth}
        \centering
        \includegraphics[width=\textwidth, trim={0.9cm 0cm 0.9cm 0cm}]{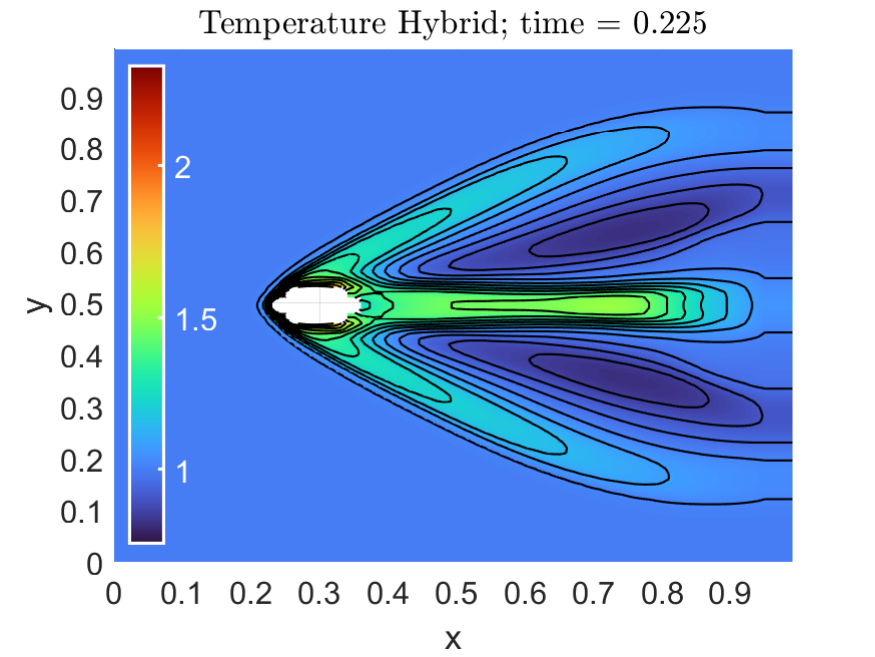}
    \end{subfigure}

    \begin{subfigure}[b]{0.32\textwidth}
        \centering
        \includegraphics[width=\textwidth, trim={0.9cm 0cm 0.9cm 0cm}]{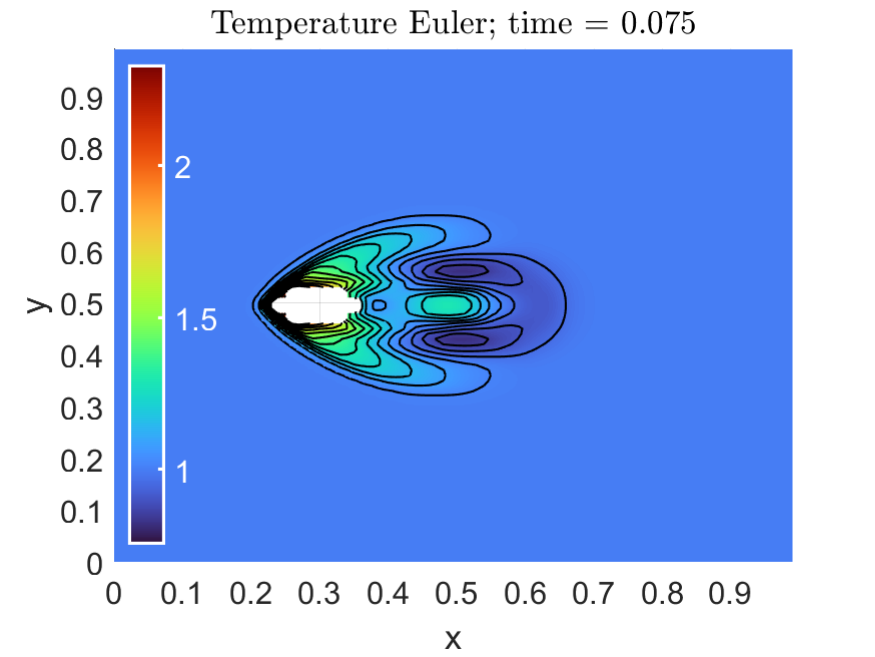}
    \end{subfigure}
    \hfill
    \begin{subfigure}[b]{0.32\textwidth}
        \centering
        \includegraphics[width=\textwidth, trim={0.9cm 0cm 0.9cm 0cm}]{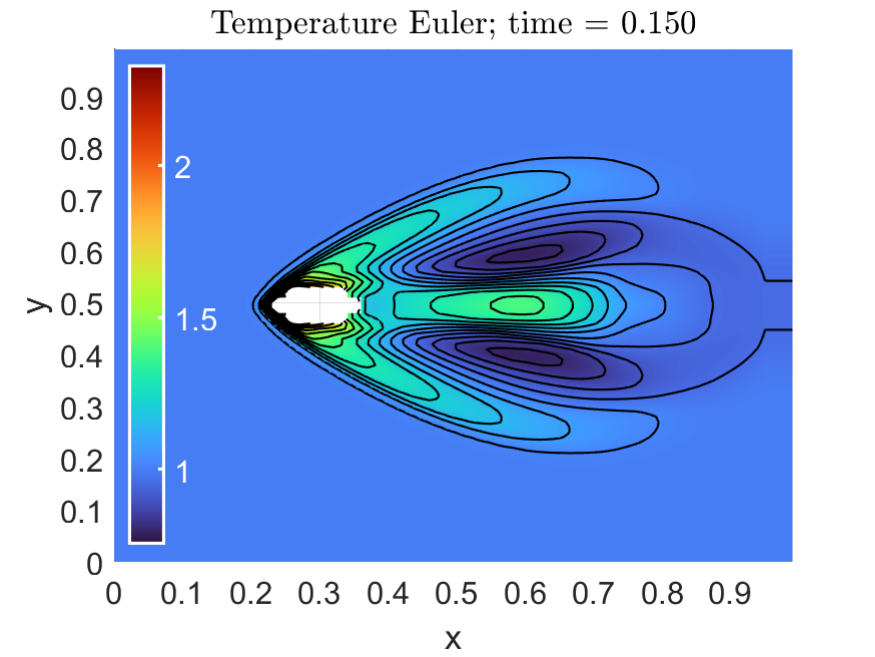}
    \end{subfigure}
    \hfill
    \begin{subfigure}[b]{0.32\textwidth}
        \centering
        \includegraphics[width=\textwidth, trim={0.9cm 0cm 0.9cm 0cm}]{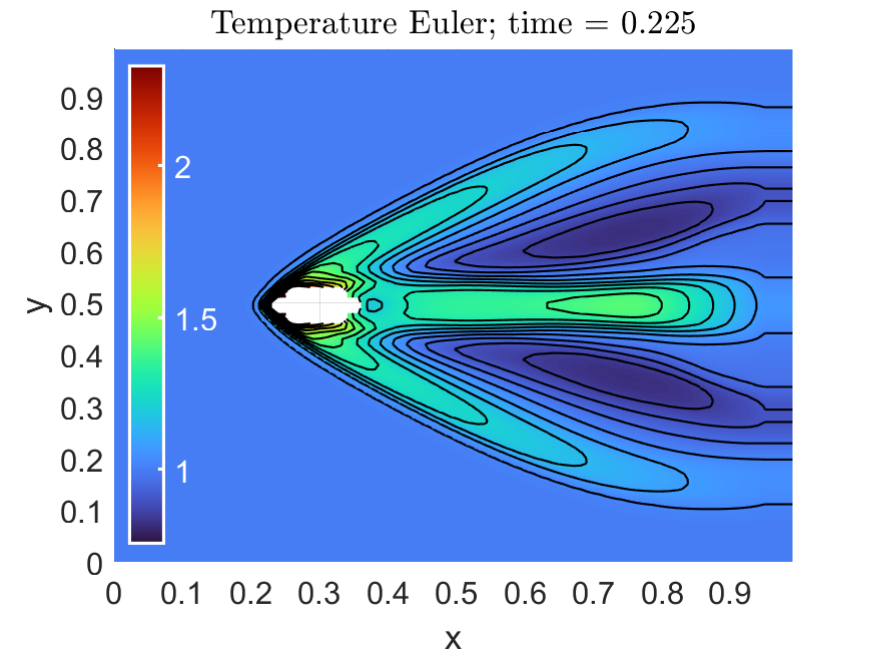}
    \end{subfigure}
    \begin{subfigure}[b]{0.32\textwidth}
        \centering
        \includegraphics[width=\textwidth, trim={0.9cm 0cm 0.9cm 0cm}]{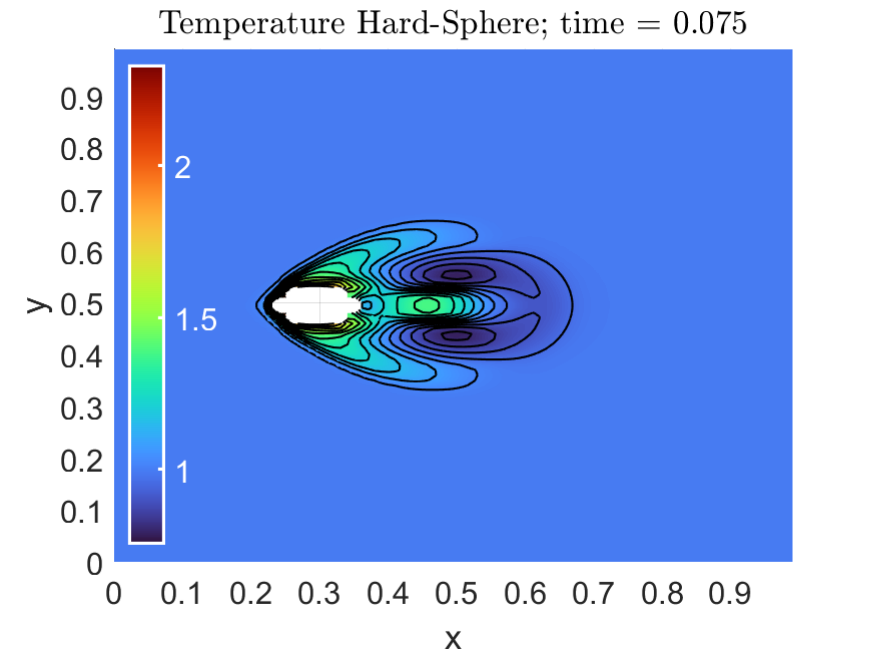}
    \end{subfigure}
    \hfill
    \begin{subfigure}[b]{0.32\textwidth}
        \centering
        \includegraphics[width=\textwidth, trim={0.9cm 0cm 0.9cm 0cm}]{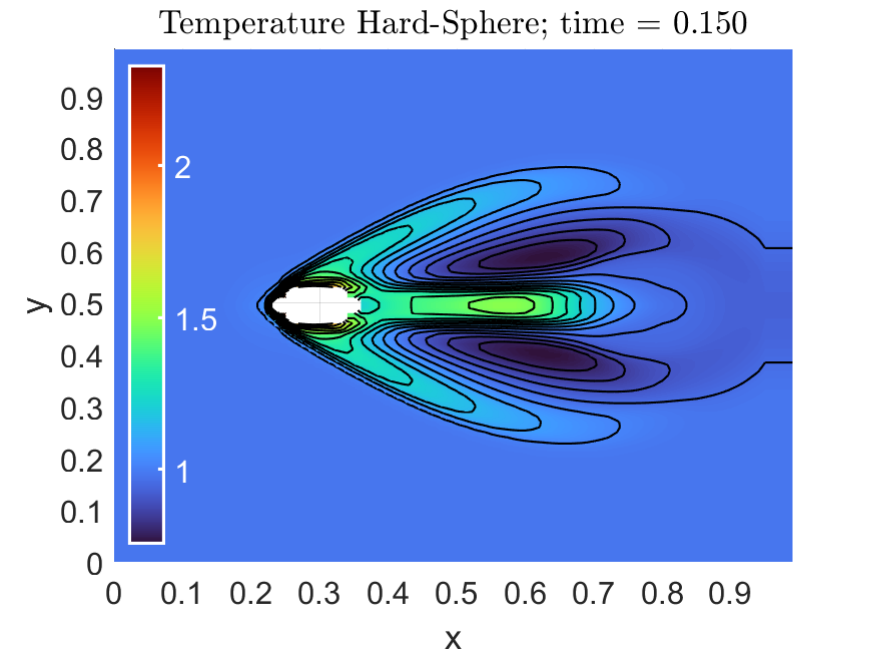}
    \end{subfigure}
    \hfill
    \begin{subfigure}[b]{0.32\textwidth}
        \centering
        \includegraphics[width=\textwidth, trim={0.9cm 0cm 0.9cm 0cm}]{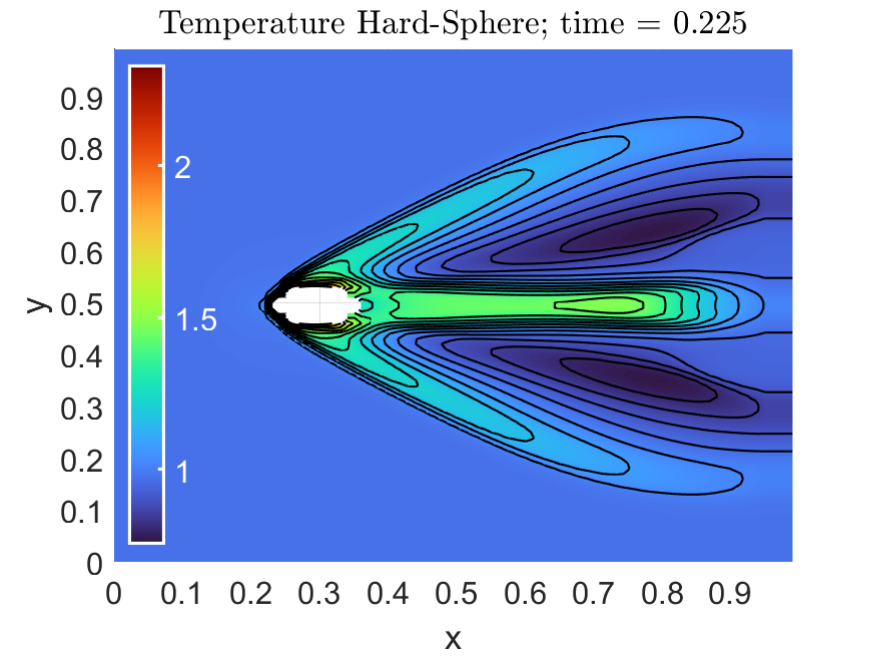}
    \end{subfigure}

    \caption{Time evolution of temperature for \hyperref[subTestFluidFixedRectObs]{\textbf{Test 2}}.  In the first, second and third row are displayed, respectively, the time evolution of temperature computed using hybrid scheme, Euler equations and Boltzmann equation.}
    \label{fig_temperature_sim2}
\end{figure}

To simulate the presence of a rectangular object, we have defined a rectangular domain (white cells in Figure \ref{fig_density_sim2} and Figure \ref{fig_temperature_sim2}) in which the integrators are never executed. These cells, are part of the stencils of the closest cells of the integration domain, thus we have assigned to them a constant and uniform value of temperature, pressure, density and velocity, given by: $P_*=1.5, \rho_*=1.0, T_*=P_* / \rho_*$ and $\vec{u}_*=\vec{0}$. This rectangular shape domain, is surrounded by a small region (2 cells thick) in which the integrator is always the one associated to the full Boltzmann equation with Boltzmann operator.\\
The fluid in the integration domain is initialized with the following parameters: $P=1.0, \rho=1.0, \vec{u}=(3.0, 0.0, 0.0)$. The distribution functions are initialized everywhere as a Maxwellian whose moments corresponds to the initial values of the Euler equations.

% In Figure \ref{fig4} is reported a schematic view of the test. The white cells constitute the fixed object, and numerical integration is never performed in them. Since these cells constitute the stencil cells of the green rhomboid cells, it is necessary to assign to them a temperature, pressure, density and velocity value, which are given respectively by $P_*=2.5, \rho_*=1.0, T_*=P_* / \rho_*$ and $\vec{u}_*=\vec{0}$. These four values are always constant over time. Inside the green cells, the hard-sphere operator is always integrated. In all other (blue dashed) cells, Algorithm \ref{alg1} is executed, and so the integration in each of this cells is performed using the respective regime integrator. The red arrows schematize the direction of the fluid flow. All dashed cells are initialised with the following density, pressure and velocity values, $\tilde{\rho}=1.0$, $\tilde{P}=1.0$ and $\vec{u}=(3.0, 0,0)$. The dimensions of the object in Figure \ref{fig4} correspond to those actually implemented in the code.\\
In Figures \ref{fig_density_sim2} and \ref{fig_temperature_sim2}, are showed, respectively, the density and temperature profiles, with $\epsilon = 10^{-4}$, at time $t=0.075,0.15, 0.225$ for the hybrid scheme, the full Euler integrator and for the Boltzmann operator (computed using CUDA GPU parallelization), as well as the domain indicators for the hybrid scheme (red, green and black dots indicate cells computed using, respectively, Euler equations, ES-BGK operator and Boltzmann operator).
%It is interesting to note that some cells that were initially hydrodynamic (red) become kinetic, and then they return hydrodynamic again at the end of the simulations, since the distribution function converges to equilibrium (Maxwellian).
In this simulation we considered $\Xi=5/3$ as specific heat ratio. We estimated, as explained in Section \ref{subNumSim}, the total duration of the numerical integration performed using just the full Boltzmann operator on CPU over the entire domain, and the ratio between the execution time for the full Boltzmann and the hybrid scheme is equal to 3.
%2.4. 
The thresholds used are: $\eta_0=1\times 10^{-3}$, $\eta_1=2.8\times10^{-5}$, $\delta_0=10^{-4}$. The time step is equal to $0.1 dx$.

It is important to remark that the kinetic regions tend to be confined in time, and in particular the one computed with the Boltzmann operator, which is the most computationally expensive. Indeed the Boltzmann operator is triggered only really close to the object and along the frontal shocks. This feature show the robustness and effectiveness of the implemented criterion.
% \begin{figure}
%     \centering
%     \includegraphics[width=0.9\linewidth]{Figure2 - Rotate.pdf}
%     \caption{Schematization of the geometric implementation of \hyperref[subTestFluidFixedRectObs]{\textbf{Test 2}}. In the dashed (blue) regions, Algorithm \ref{alg1} is executed at each time step. The white cells constitute the fixed object, and numerical integration is never performed in them. Since some of these cells constitute the stencil cells of the green rhomboid cells, we have arbitrarily assigned fixed temperature, pressure, density and velocity values. These same values are also assigned at each iteration to the green cells, into which the Boltzmann operator is always integrated. The red arrows schematize the direction of the fluid flow.}
%     \label{fig4}
% \end{figure}

\subsection{Test 3: Moving object in a fluid at rest}
\label{subTestFluidMovingRest}

The aim of this test is to show the capability of the scheme to handle complex moving geometries, real world phenomena and moving boundaries inside the integration domain. In particular we use our code to simulate the flow generated by a moving object in a fluid at rest, which is embedded in a computational domain of size $L_x=L_y=1.0$ discretized with 100 points in each direction. We have discretized the velocity domain using 16 points in the three directions, and we considered a cut-off domain $\Omega_v=[-8, 8]^3$. The value of $R$ has been taken approximately equal to 3.3, and both numbers of points $A_1$ and $A_2$ discretizing the spherical integration of the Boltzmann operator (see Section \ref{subBoltzOp}) have been taken equal to 4.

\begin{figure}
    \centering
    \begin{subfigure}[b]{0.32\textwidth}
        \centering
        \includegraphics[width=\textwidth, trim={0.9cm 0cm 0.9cm 0cm}]{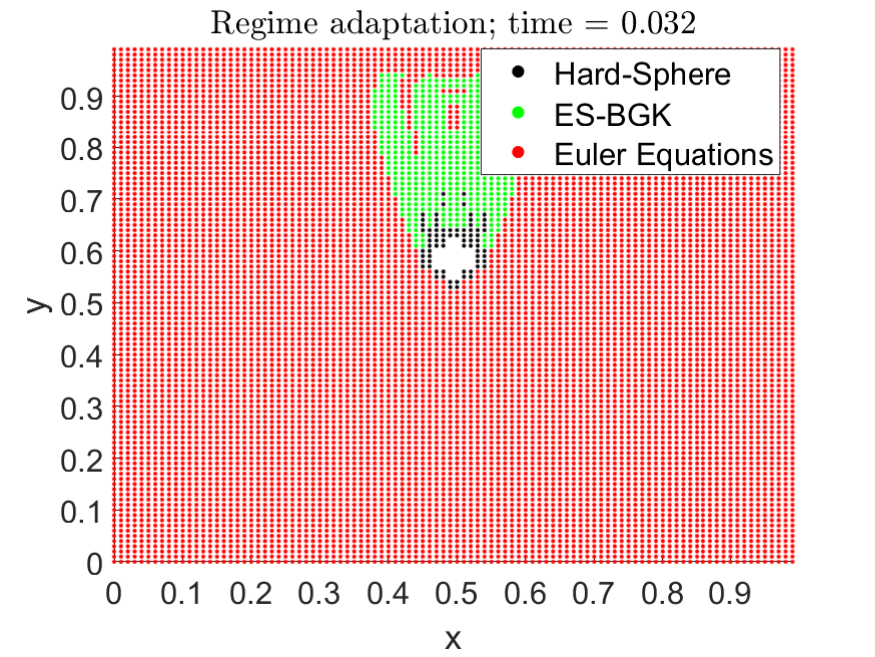}
    \end{subfigure}
    \hfill
    \begin{subfigure}[b]{0.32\textwidth}
        \centering
        \includegraphics[width=\textwidth, trim={0.9cm 0cm 0.9cm 0cm}]{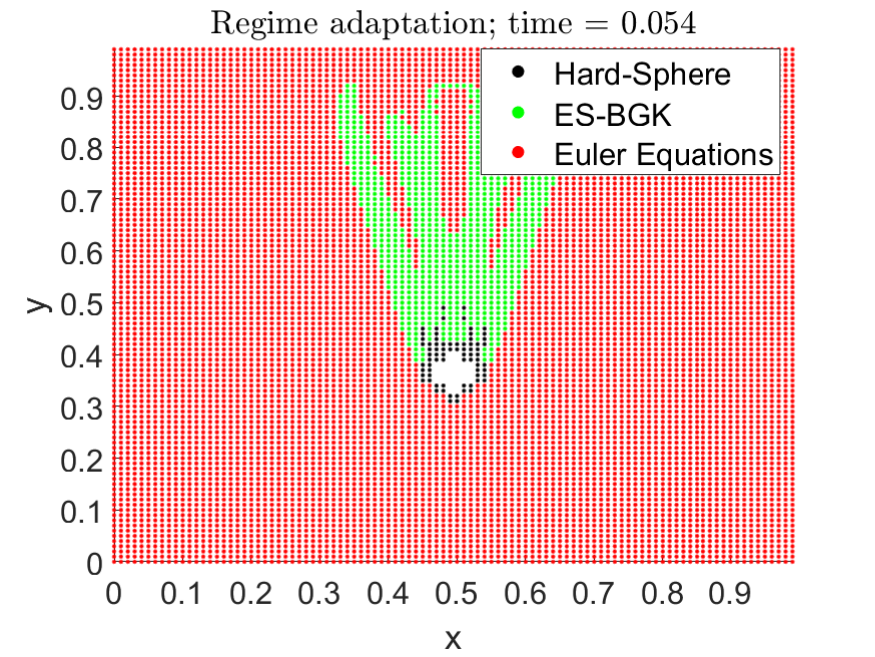}
    \end{subfigure}
    \hfill
    \begin{subfigure}[b]{0.32\textwidth}
        \centering
        \includegraphics[width=\textwidth, trim={0.9cm 0cm 0.9cm 0cm}]{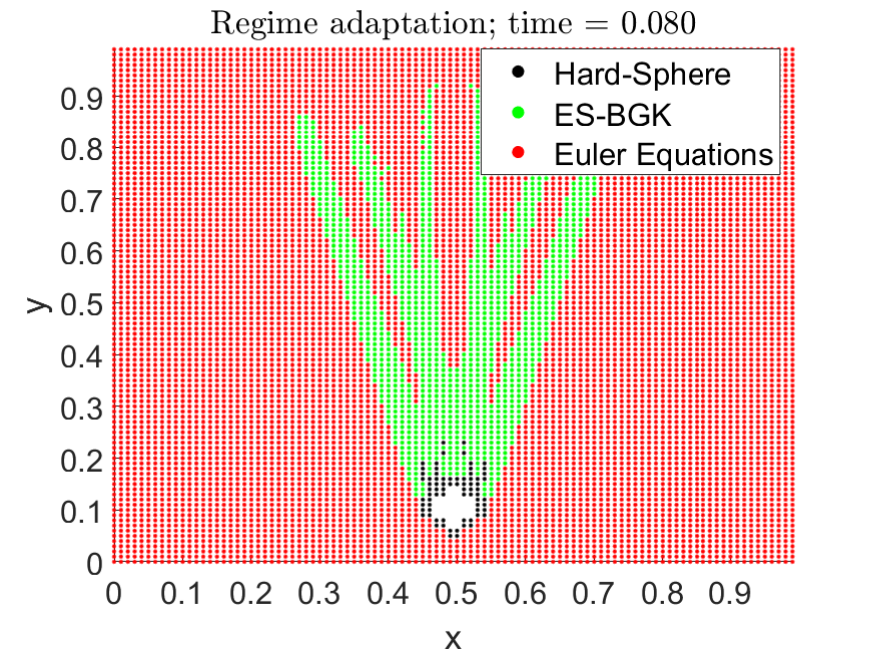}
    \end{subfigure}

    \begin{subfigure}[b]{0.32\textwidth}
        \centering
        \includegraphics[width=\textwidth, trim={0.9cm 0cm 0.9cm 0cm}]{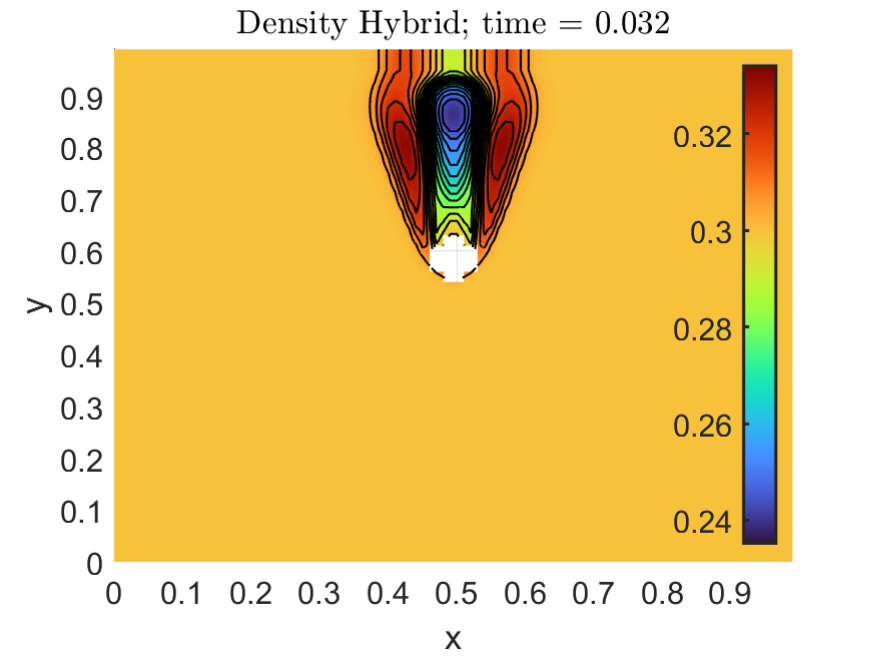}
    \end{subfigure}
    \hfill
    \begin{subfigure}[b]{0.32\textwidth}
        \centering
        \includegraphics[width=\textwidth, trim={0.9cm 0cm 0.9cm 0cm}]{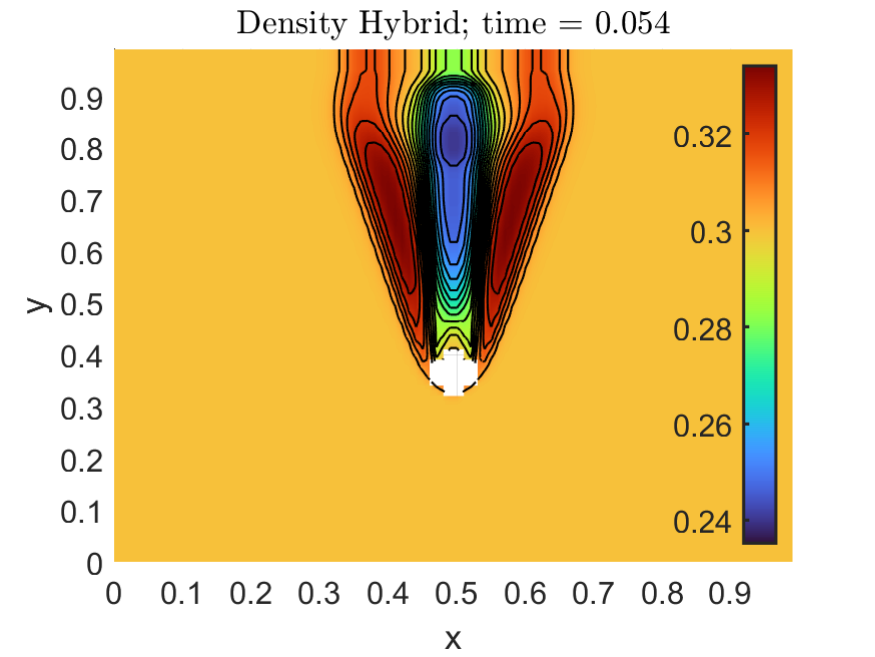}
    \end{subfigure}
    \hfill
    \begin{subfigure}[b]{0.32\textwidth}
        \centering
        \includegraphics[width=\textwidth, trim={0.9cm 0cm 0.9cm 0cm}]{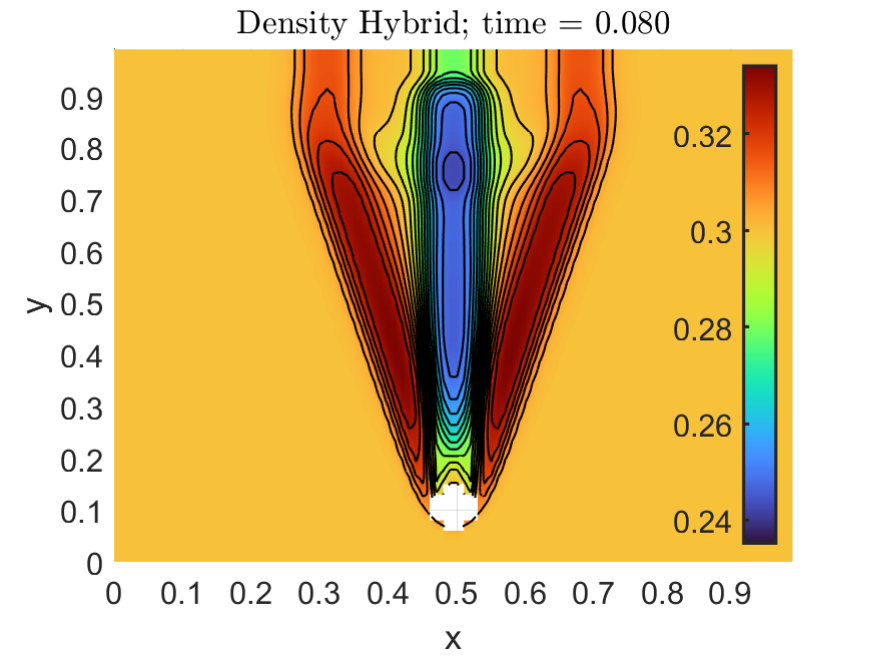}
    \end{subfigure}

    \begin{subfigure}[b]{0.32\textwidth}
        \centering
        \includegraphics[width=\textwidth, trim={0.9cm 0cm 0.9cm 0cm}]{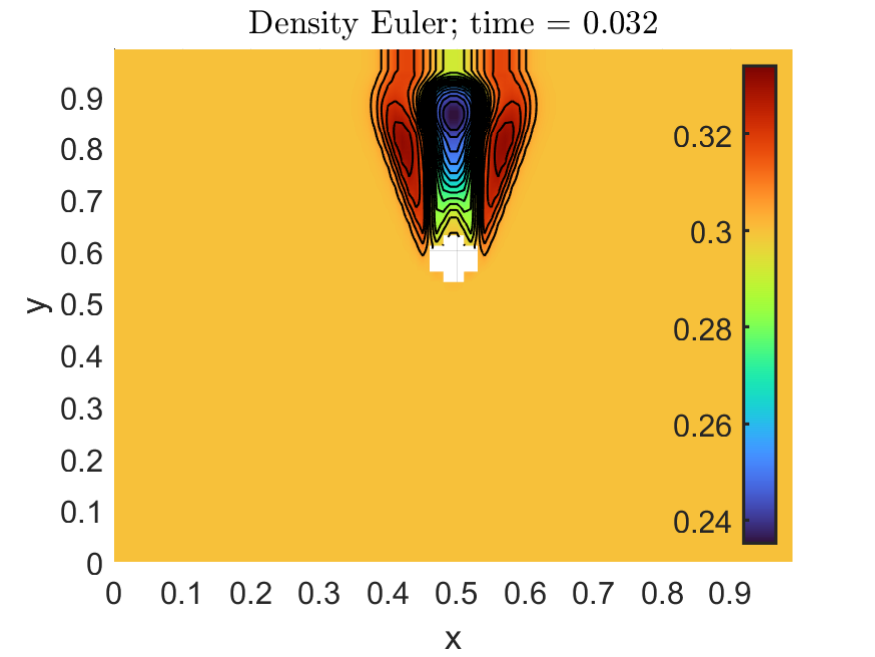}
    \end{subfigure}
    \hfill
    \begin{subfigure}[b]{0.32\textwidth}
        \centering
        \includegraphics[width=\textwidth, trim={0.9cm 0cm 0.9cm 0cm}]{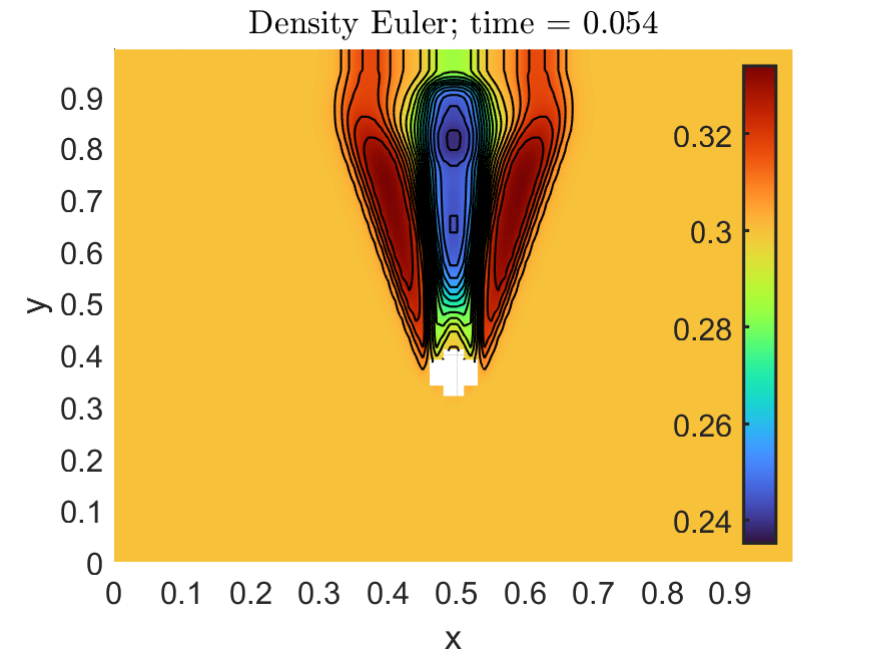}
    \end{subfigure}
    \hfill
    \begin{subfigure}[b]{0.32\textwidth}
        \centering
        \includegraphics[width=\textwidth, trim={0.9cm 0cm 0.9cm 0cm}]{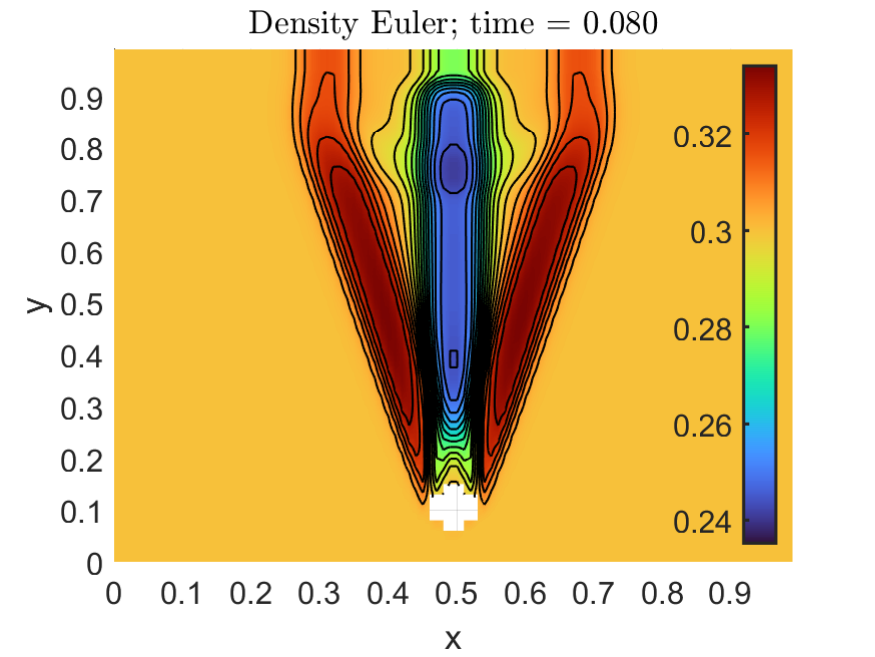}
    \end{subfigure}
    \begin{subfigure}[b]{0.32\textwidth}
        \centering
        \includegraphics[width=\textwidth, trim={0.9cm 0cm 0.9cm 0cm}]{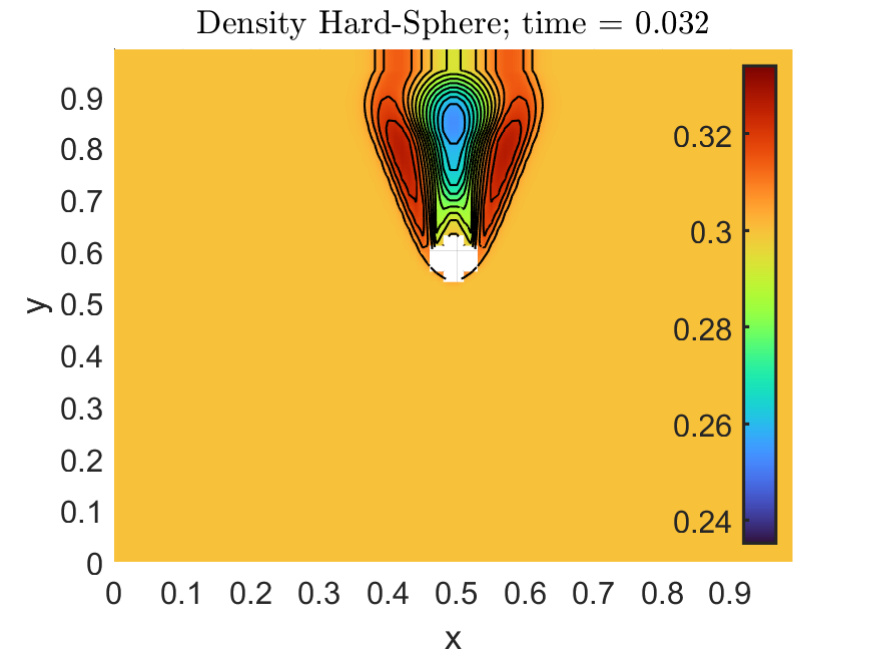}
    \end{subfigure}
    \hfill
    \begin{subfigure}[b]{0.32\textwidth}
        \centering
        \includegraphics[width=\textwidth, trim={0.9cm 0cm 0.9cm 0cm}]{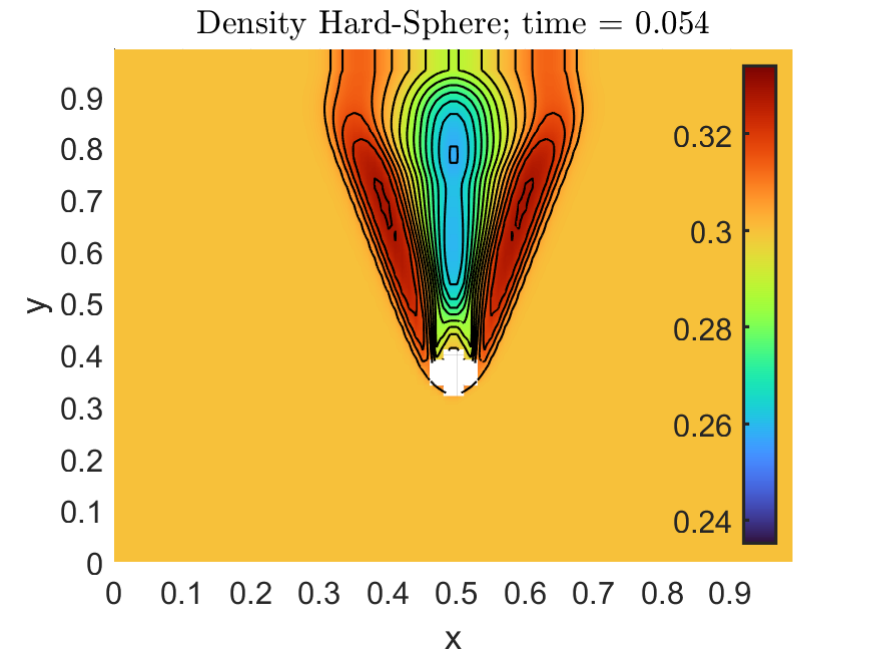}
    \end{subfigure}
    \hfill
    \begin{subfigure}[b]{0.32\textwidth}
        \centering
        \includegraphics[width=\textwidth, trim={0.9cm 0cm 0.9cm 0cm}]{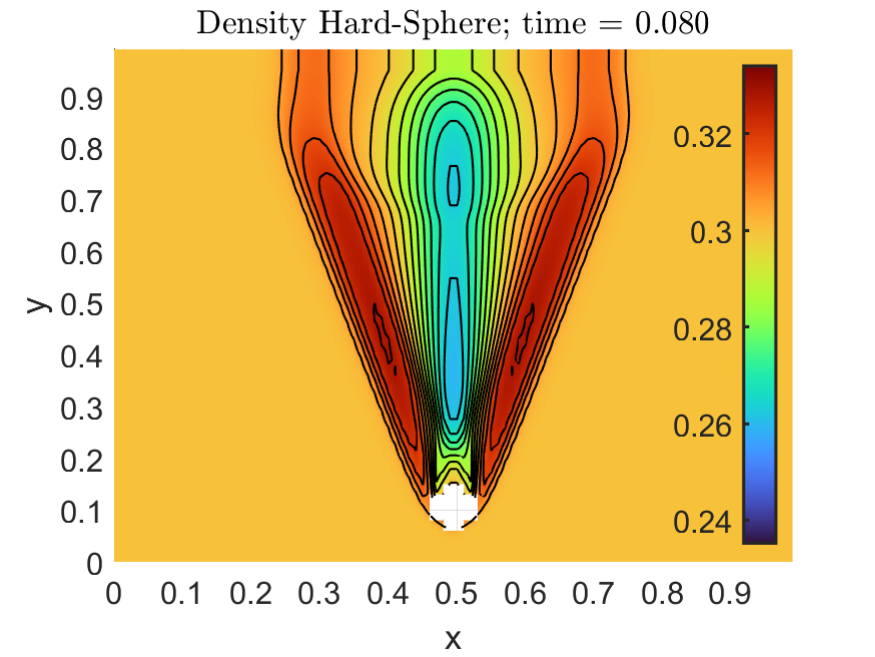}
    \end{subfigure}
    
    \caption{Time evolution of density for \hyperref[subTestFluidMovingRest]{\textbf{Test 3}}. The first row  displays the regime adaptation at three different times. Red, green and black dots, identify, respectively, cells updated using Euler equations, ES-BGK equation and Boltzmann equation. In the second, third and fourth row are displayed, respectively, the time evolution of density computed using hybrid scheme, Euler equations and Boltzmann equation.
    }
    \label{fig_density_sim3}
\end{figure}

\begin{figure}
    \begin{subfigure}[b]{0.32\textwidth}
        \centering
        \includegraphics[width=\textwidth, trim={0.9cm 0cm 0.9cm 0cm}]{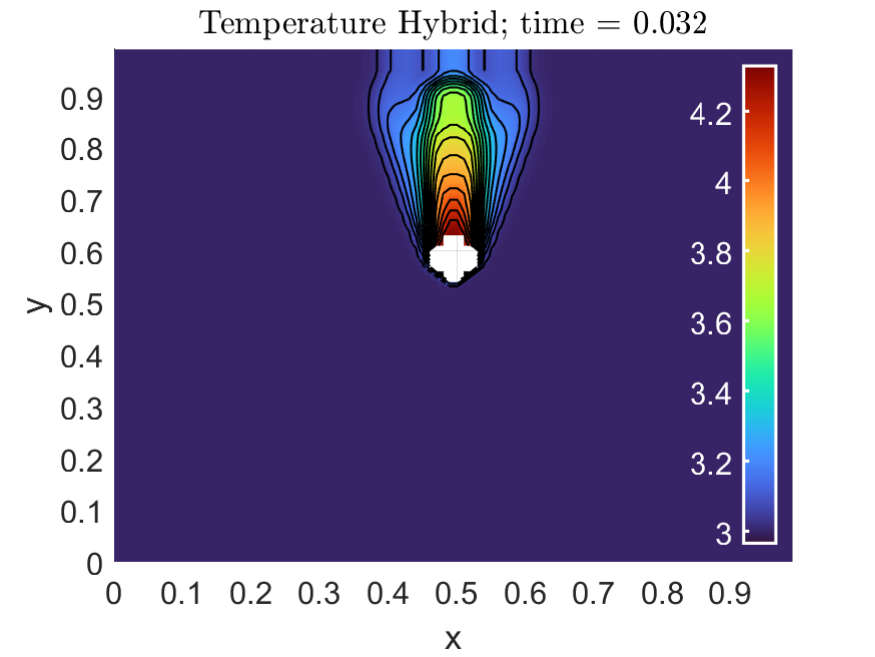}
    \end{subfigure}
    \hfill
    \begin{subfigure}[b]{0.32\textwidth}
        \centering
        \includegraphics[width=\textwidth, trim={0.9cm 0cm 0.9cm 0cm}]{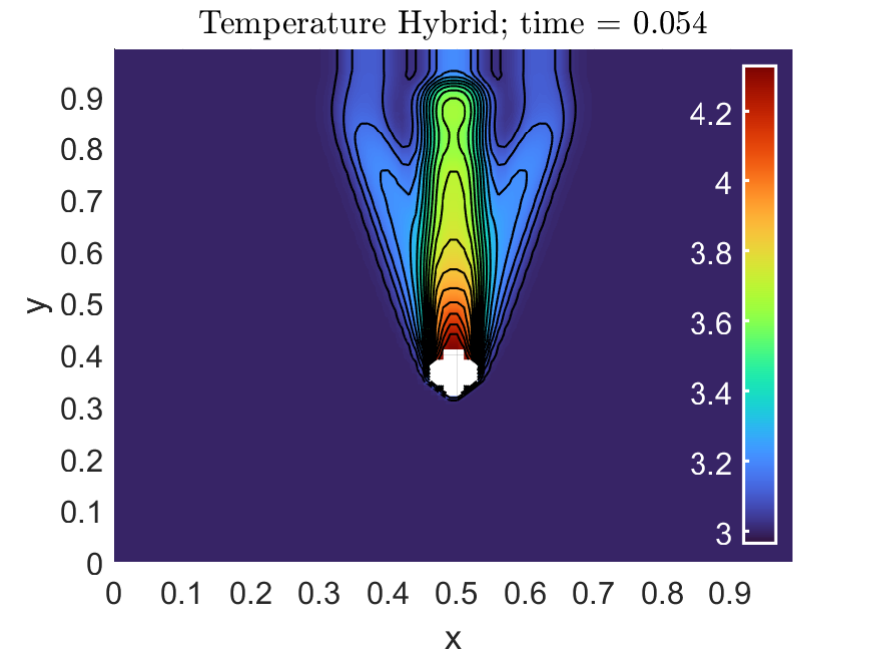}
    \end{subfigure}
    \hfill
    \begin{subfigure}[b]{0.32\textwidth}
        \centering
        \includegraphics[width=\textwidth, trim={0.9cm 0cm 0.9cm 0cm}]{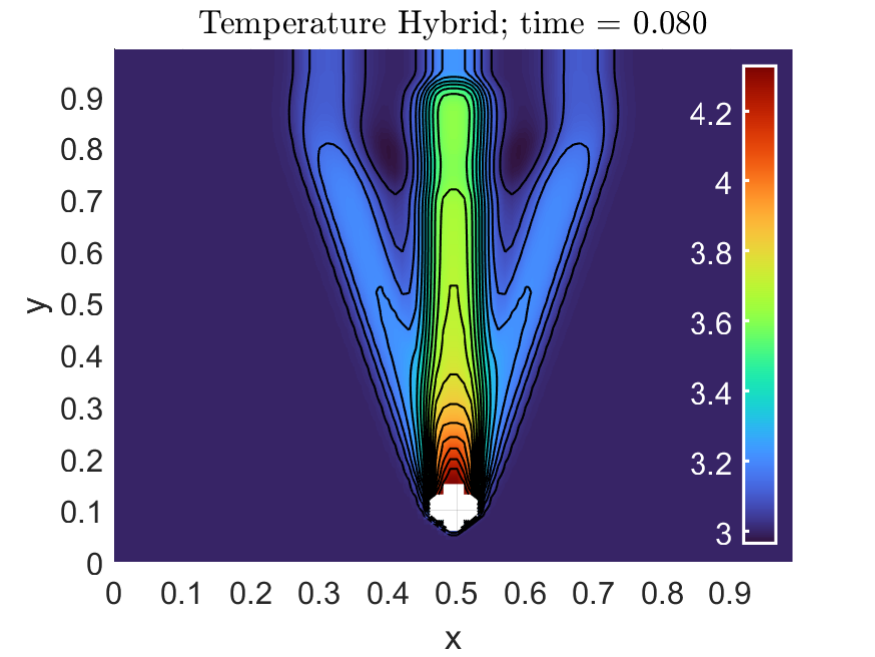}
    \end{subfigure}

    \begin{subfigure}[b]{0.32\textwidth}
        \centering
        \includegraphics[width=\textwidth, trim={0.9cm 0cm 0.9cm 0cm}]{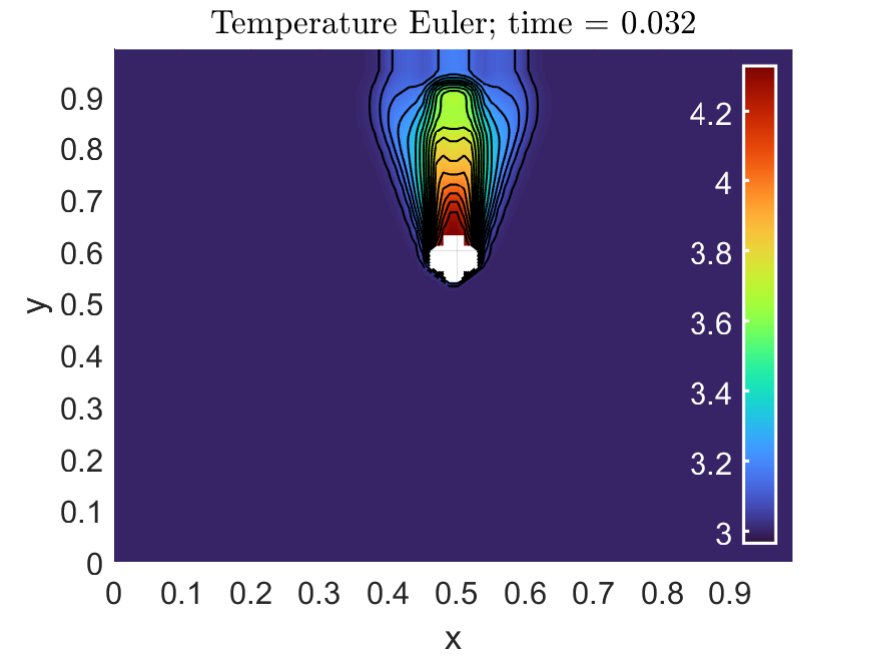}
    \end{subfigure}
    \hfill
    \begin{subfigure}[b]{0.32\textwidth}
        \centering
        \includegraphics[width=\textwidth, trim={0.9cm 0cm 0.9cm 0cm}]{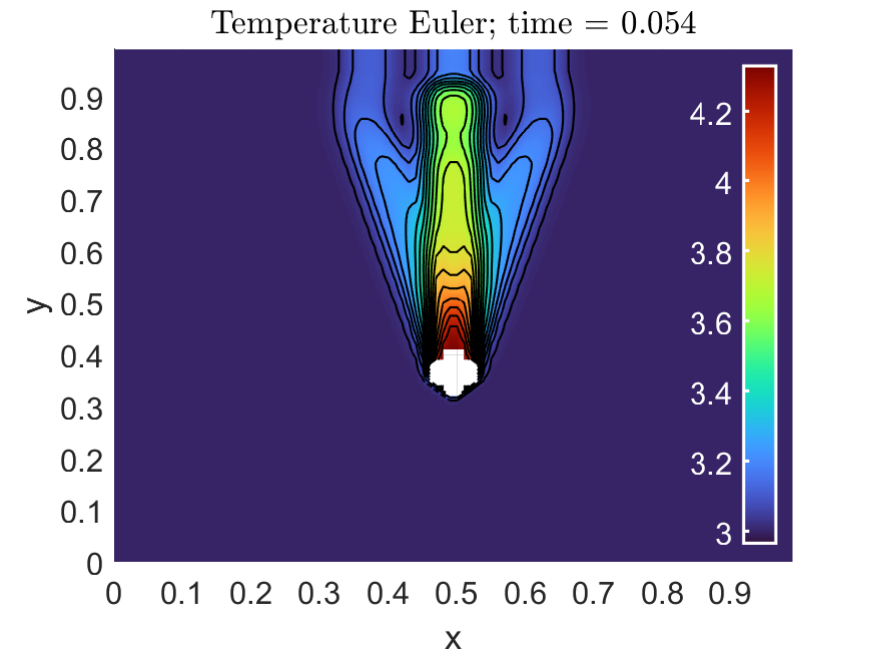}
    \end{subfigure}
    \hfill
    \begin{subfigure}[b]{0.32\textwidth}
        \centering
        \includegraphics[width=\textwidth, trim={0.9cm 0cm 0.9cm 0cm}]{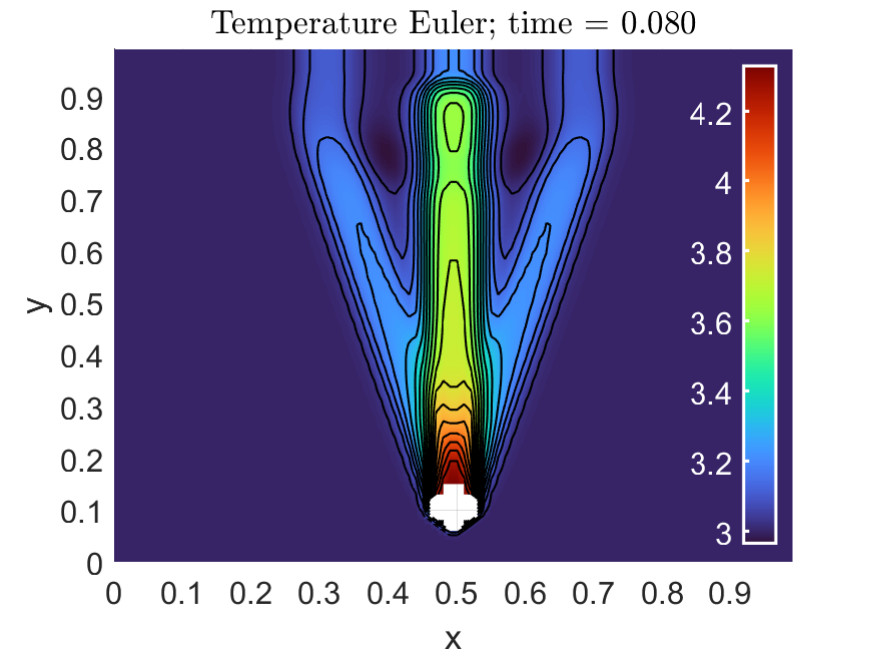}
    \end{subfigure}

    \begin{subfigure}[b]{0.32\textwidth}
        \centering
        \includegraphics[width=\textwidth, trim={0.9cm 0cm 0.9cm 0cm}]{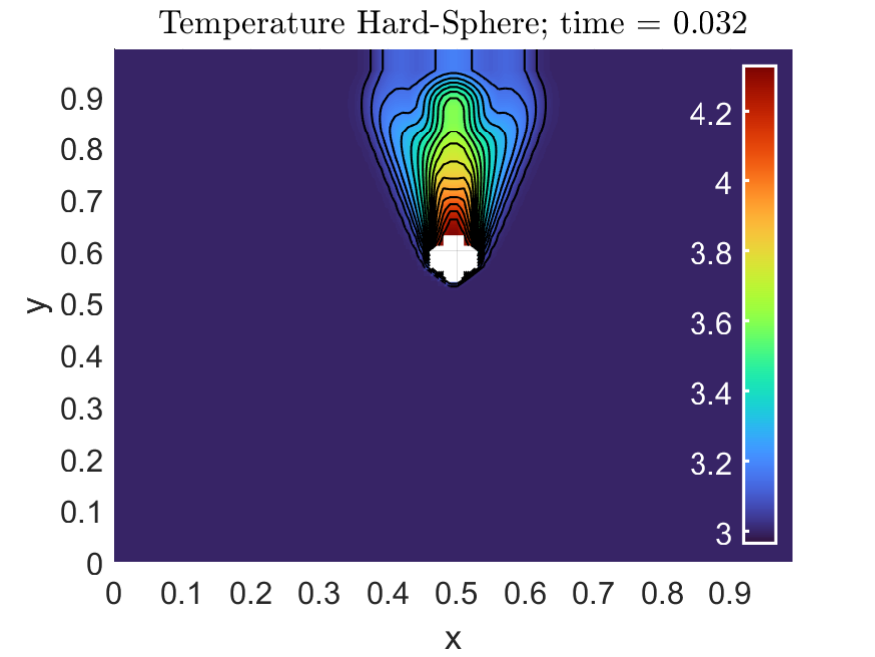}
    \end{subfigure}
    \hfill
    \begin{subfigure}[b]{0.32\textwidth}
        \centering
        \includegraphics[width=\textwidth, trim={0.9cm 0cm 0.9cm 0cm}]{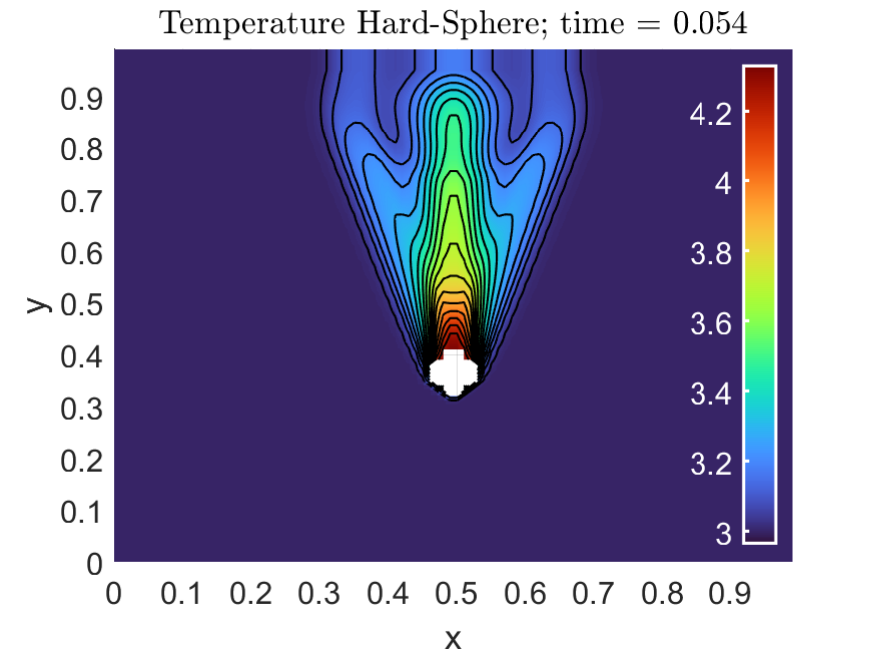}
    \end{subfigure}
    \hfill
    \begin{subfigure}[b]{0.32\textwidth}
        \centering
        \includegraphics[width=\textwidth, trim={0.9cm 0cm 0.9cm 0cm}]{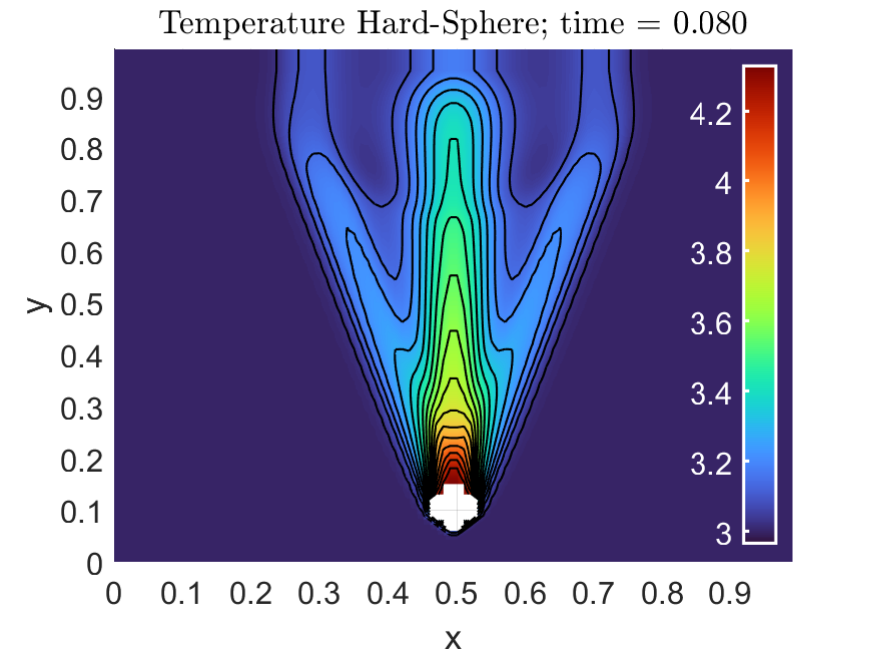}
    \end{subfigure}

    \caption{Time evolution of temperature for \hyperref[subTestFluidMovingRest]{\textbf{Test 3}}. In the first, second and third row are displayed, respectively, the time evolution of temperature computed using hybrid scheme, Euler equations and Boltzmann equation.}
    \label{fig_temperature_sim3}
\end{figure}

To simulate the presence of a rectangular object, we have defined a rectangular domain (white cells in Figures \ref{fig_density_sim3} and \ref{fig_temperature_sim3}) in which the integrators are never executed. These cells, are part of the stencils of the closest cells of the integration domain, thus we have assigned to them a constant and uniform value of pressure, density, temperature and velocity, given by: $P_*=1.3, \rho_*=0.3, T_*=P_* / \rho_*$ and $\vec{u}_*=\vec{0}$. This rectangular shape domain, is surrounded by a small region (2 cells thick) in which the integrator is always the one associated to the full Boltzmann equation with Boltzmann operator. The objects moves in the negative direction of the y-axis, by one cell, at each timestep.
The fluid in the integration domain is initialized with the following values of pressure, density, temperature and velocity: $P=0.9, \rho=0.3, T=P / \rho$ and $\vec{u}=\vec{0}$. The distribution functions are initialized everywhere as a Maxwellian, whose moments correspond to the initial values of the Euler equations.

% In Figure \ref{fig6} is reported a schematic view of the test. In the dotted region, our algorithm is executed at each timestep. The white cells constitute the moving object, and numerical integration is never performed in them. The same Figure shows the position of the same object at two different time instants, at time $t_0$ and time $t_1>t_0$. Since the white cells constitute the stencil cells of the green rhomboid cells, it is necessary to assign to them a temperature, pressure, density and velocity value, which are given respectively by $P_*=2.2, \rho_*=0.3, T_*=P_* / \rho_*$ and $\vec{u}_*$. These four values are always constant over time. Inside the rhomboid green cells the full Boltzmann equation is always integrated. In all other (dashed) cells, Algorithm \ref{alg1} is executed, and then the respective integrator domain is computed. The green arrows schematize the direction of velocity of the moving object. The object moves in the negative direction of the y-axis, by one cell, at each timestep of the hydrodynamic integrator. All dashed cells are initialized with the following density, pressure and velocity values, $\tilde{\rho}=0.3$, $\tilde{P}=0.3$ and $\vec{u}=(0.0, 0.0, 0.0)$. The dimensions of the object in Figure \ref{fig6} correspond to those actually implemented in the code.
In Figure \ref{fig_density_sim3} we show the density profile with $\epsilon = 10^{-4}$, at time $t=0.032, 0.054, 0.080$ for the hybrid scheme, for the full Euler equations and for the Boltzmann (Hard-Sphere) operator (computed using CUDA GPU parallelization), as well as the domain indicators for the hybrid scheme (red, green and black dots indicate cells computed using, respectively, Euler equations, ES-BGK equation and Boltzmann equation). It is interesting to note that some cells that were initially hydrodynamic (red) become kinetic, and then they return hydrodynamic again at the end of the simulations, since the distribution function converges to equilibrium (Maxwellian). In Figure \ref{fig_temperature_sim3} is reported the temperature profile for the same instants in time, for the hybrid scheme, for the Euler and Boltzmann equations. In this simulation we considered $\Xi=5/3$ as specific heat ratio. We estimated, as explained in Section \ref{subNumSim}, the total duration of the numerical integration performed using just the full Boltzmann operator on CPU over the entire domain, and the ratio between the execution time for the full Boltzmann and the hybrid scheme (Euler equations, ES-BGK operator and Boltzmann operator) is equal to 2.4.

%3.1. 
The thresholds used are: $\eta_0=1\times 10^{-3}$, $\eta_1=2.8\times10^{-5}$, $\delta_0=10^{-4}$. The timestep $t_H$ associated to the hydrodynamic domain of the hybrid scheme has been updated in such a way to guarantee $CFL=0.95$, and the timestep is $0.1 dx$. It is important to remark that the kinetic regions tend to be confined in time, and in particular the one computed with the Boltzmann operator, which is the most computationally expensive. Indeed the Boltzmann operator is triggered only really close to the object and along the frontal shocks. This feature show the robustness and effectiveness of the implemented criterion.

% \begin{figure}
%     \centering
%     \includegraphics[scale=0.7, trim={0cm 8cm 0cm 0cm}]{Figure6.pdf}
%     \caption{Schematic view of the geometry of \hyperref[subTestFluidMovingRest]{\textbf{Test 3}}, at two given times. In the dotted region, Algorithm \ref{alg1} is executed at each timestep. The white cells constitute the moving object, where numerical integration is never performed. Some of the white cells constitute the stencil cells of the green rhomboid cells, where temperature, pressure, density and velocity values are constant over time. Inside the rhomboid green cells the Boltzmann operator is always integrated. The green arrows schematize the direction of motion of the moving object: it moves in the negative direction of the y-axis, by one cell, at each hydrodynamic timestep.}
%     \label{fig6}
% \end{figure}

\subsection{Test 4: Fluid flux towards a fixed round obstacle}
\label{subTestFluidObstacle}
One of the aim of this test is to show the capability of the scheme to handle round geometries. In particular we use our code to simulate the flow around a fixed round object, which is embedded in a computational domain of size $L_x=L_y=1.0$ discretized with 100 points in each direction. We have discretized the velocity domain using 16 points in the three directions, and we considered a cut-off domain $\Omega_v=[-8, 8]^3$. The value of $R$ has been taken approximately equal to 3.3, and both numbers of points $A_1$ and $A_2$ discretizing the spherical integration of the Boltzmann operator (see Section \ref{subBoltzOp}) have been taken equal to 4.

\begin{figure}
    \centering
    \begin{subfigure}[b]{0.32\textwidth}
        \centering
        \includegraphics[width=\textwidth, trim={0.9cm 0cm 0.9cm 0cm}]{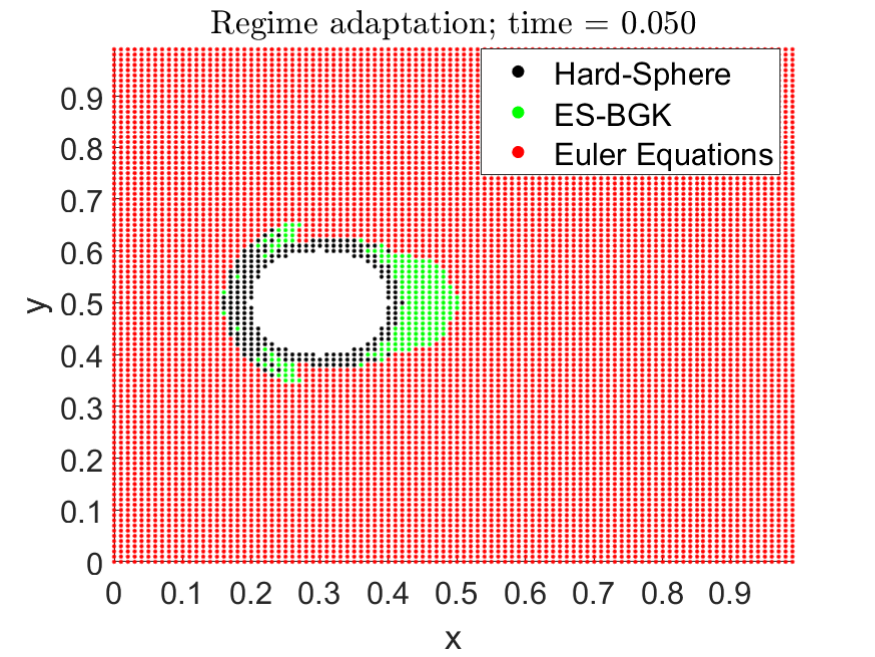}
    \end{subfigure}
    \hfill
    \begin{subfigure}[b]{0.32\textwidth}
        \centering
        \includegraphics[width=\textwidth, trim={0.9cm 0cm 0.9cm 0cm}]{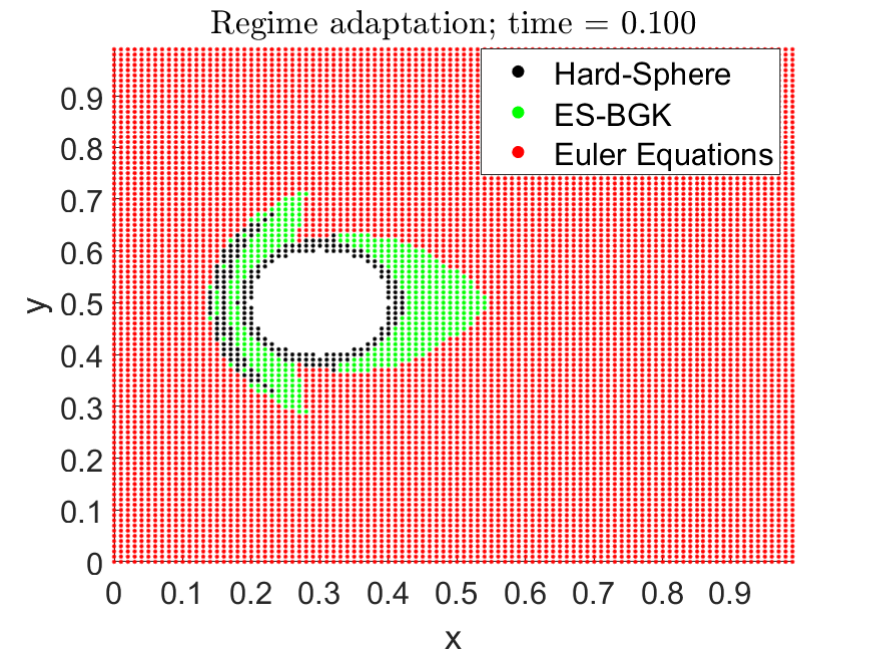}
    \end{subfigure}
    \hfill
    \begin{subfigure}[b]{0.32\textwidth}
        \centering
        \includegraphics[width=\textwidth, trim={0.9cm 0cm 0.9cm 0cm}]{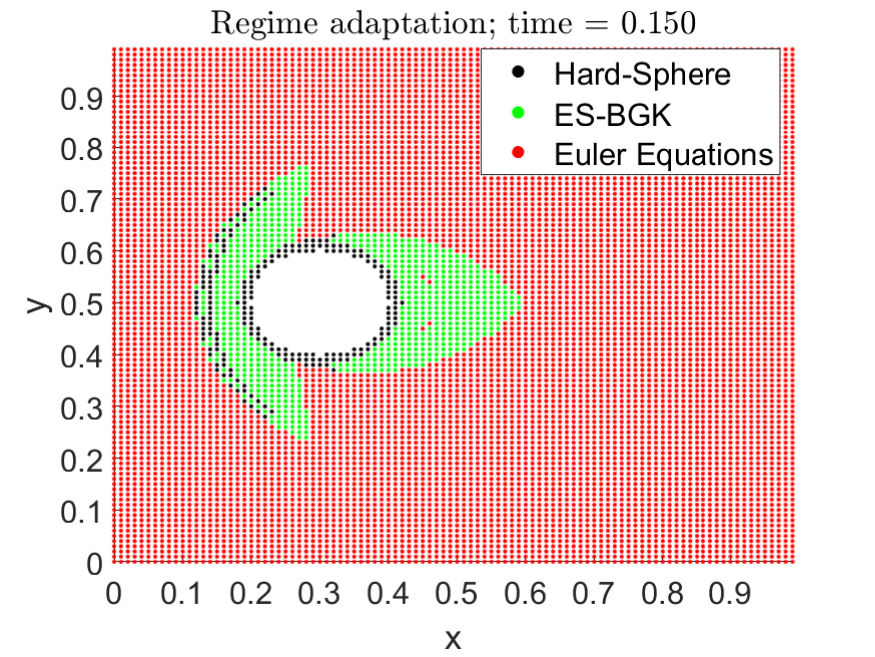}
    \end{subfigure}

    \begin{subfigure}[b]{0.32\textwidth}
        \centering
        \includegraphics[width=\textwidth, trim={0.9cm 0cm 0.9cm 0cm}]{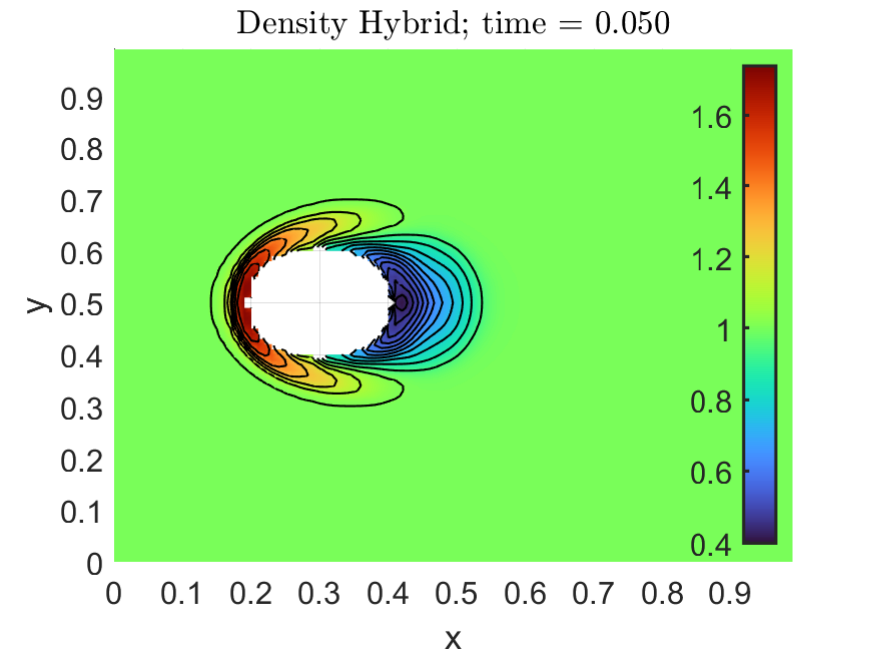}
    \end{subfigure}
    \hfill
    \begin{subfigure}[b]{0.32\textwidth}
        \centering
        \includegraphics[width=\textwidth, trim={0.9cm 0cm 0.9cm 0cm}]{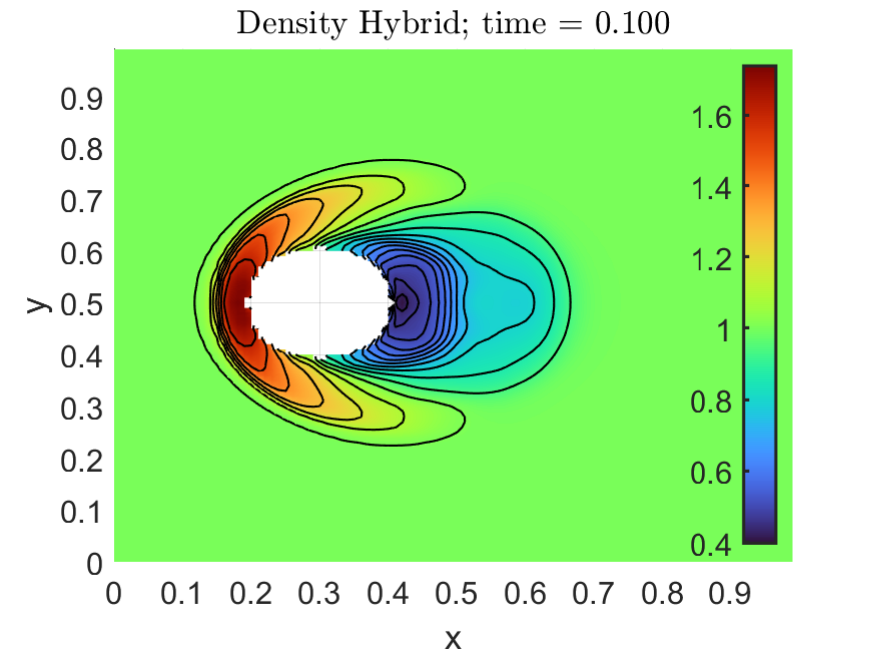}
    \end{subfigure}
    \hfill
    \begin{subfigure}[b]{0.32\textwidth}
        \centering
        \includegraphics[width=\textwidth, trim={0.9cm 0cm 0.9cm 0cm}]{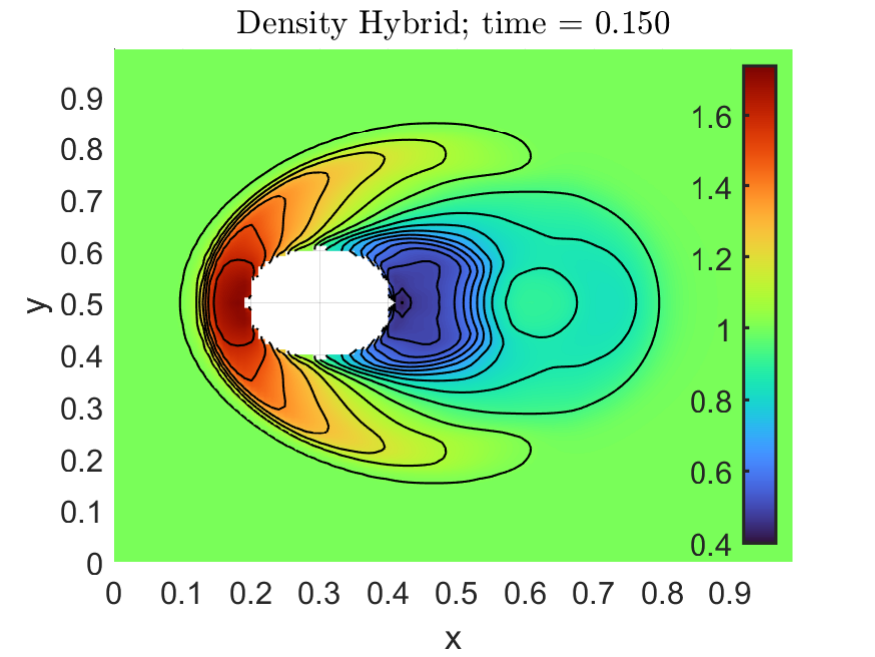}
    \end{subfigure}

    \begin{subfigure}[b]{0.32\textwidth}
        \centering
        \includegraphics[width=\textwidth, trim={0.9cm 0cm 0.9cm 0cm}]{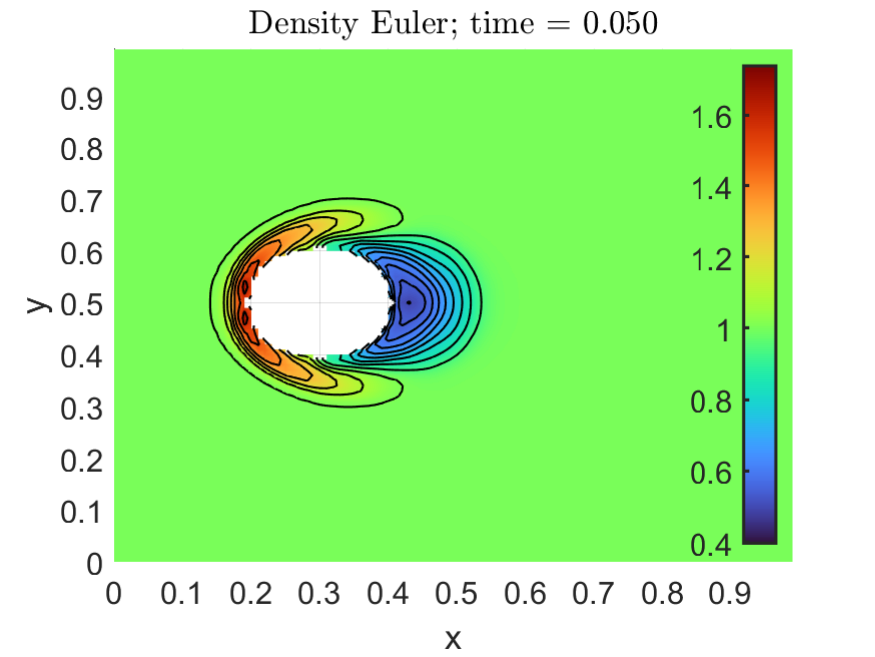}
    \end{subfigure}
    \hfill
    \begin{subfigure}[b]{0.32\textwidth}
        \centering
        \includegraphics[width=\textwidth, trim={0.9cm 0cm 0.9cm 0cm}]{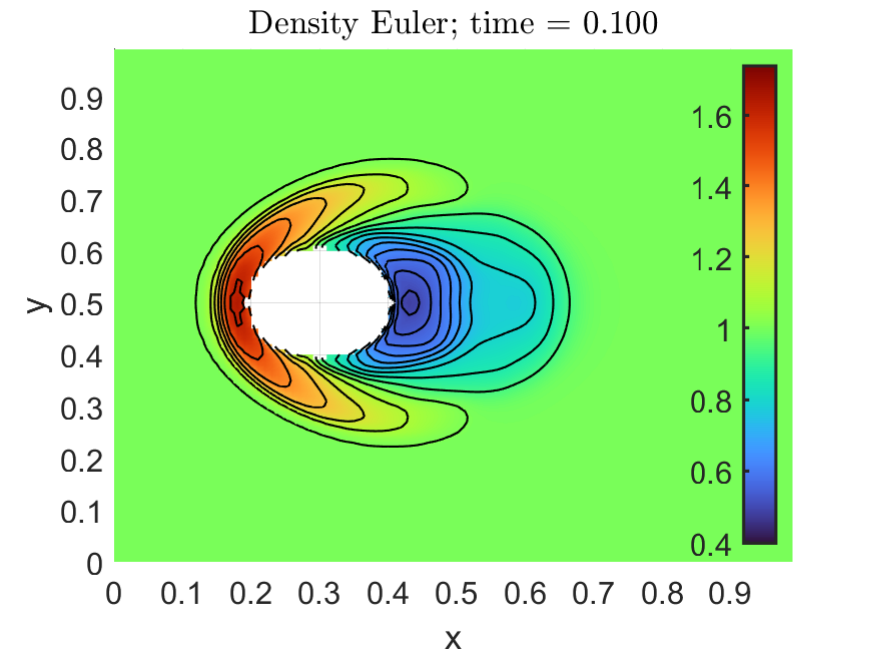}
    \end{subfigure}
    \hfill
    \begin{subfigure}[b]{0.32\textwidth}
        \centering
        \includegraphics[width=\textwidth, trim={0.9cm 0cm 0.9cm 0cm}]{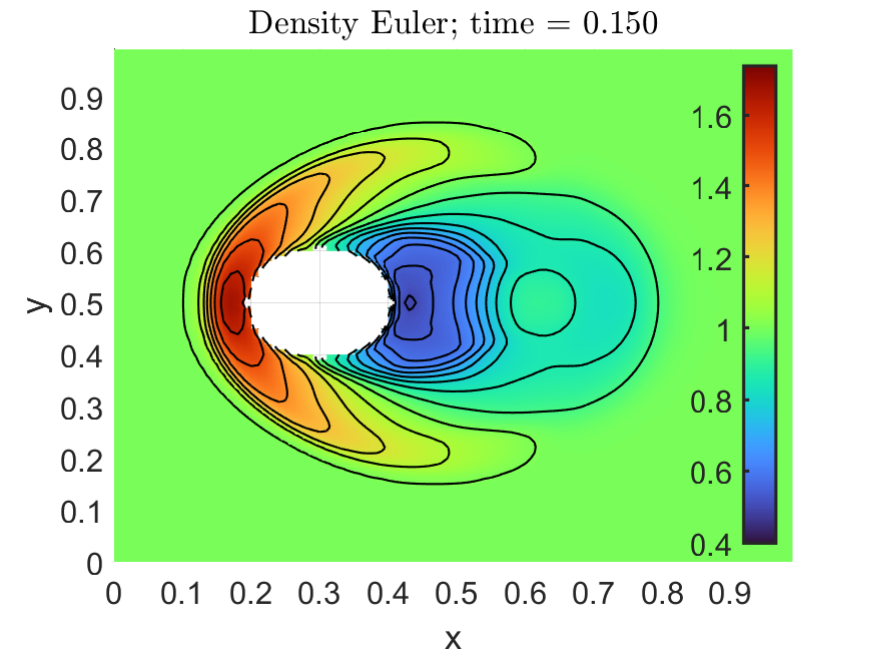}
    \end{subfigure}

    \begin{subfigure}[b]{0.32\textwidth}
        \centering
        \includegraphics[width=\textwidth, trim={0.9cm 0cm 0.9cm 0cm}]{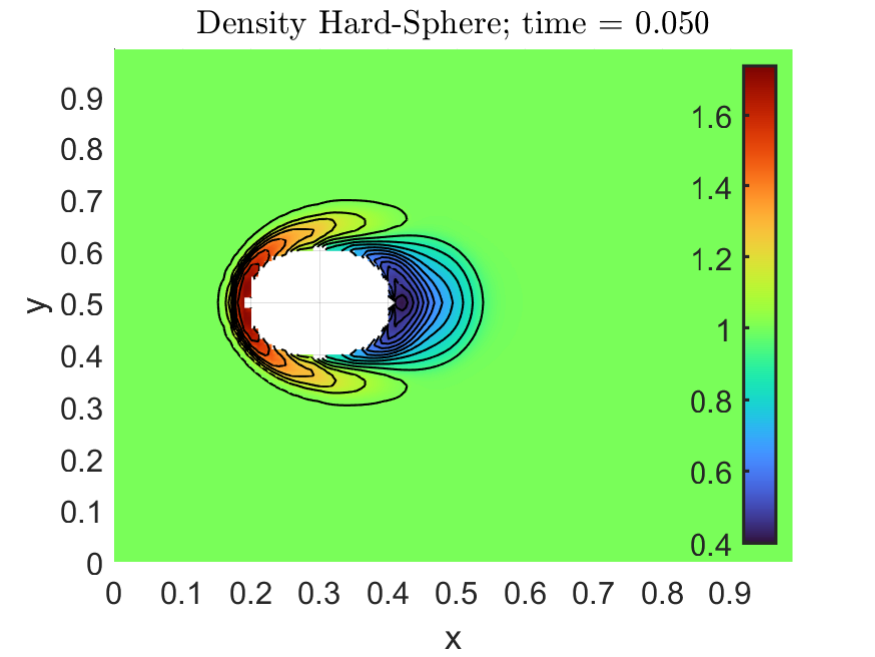}
    \end{subfigure}
    \hfill
    \begin{subfigure}[b]{0.32\textwidth}
        \centering
        \includegraphics[width=\textwidth, trim={0.9cm 0cm 0.9cm 0cm}]{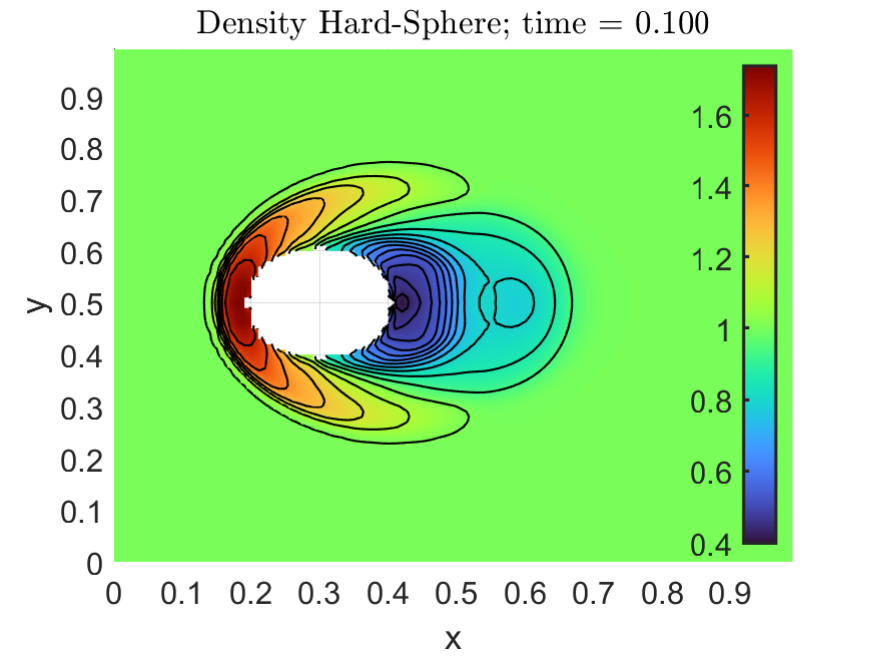}
    \end{subfigure}
    \hfill
    \begin{subfigure}[b]{0.32\textwidth}
        \centering
        \includegraphics[width=\textwidth, trim={0.9cm 0cm 0.9cm 0cm}]{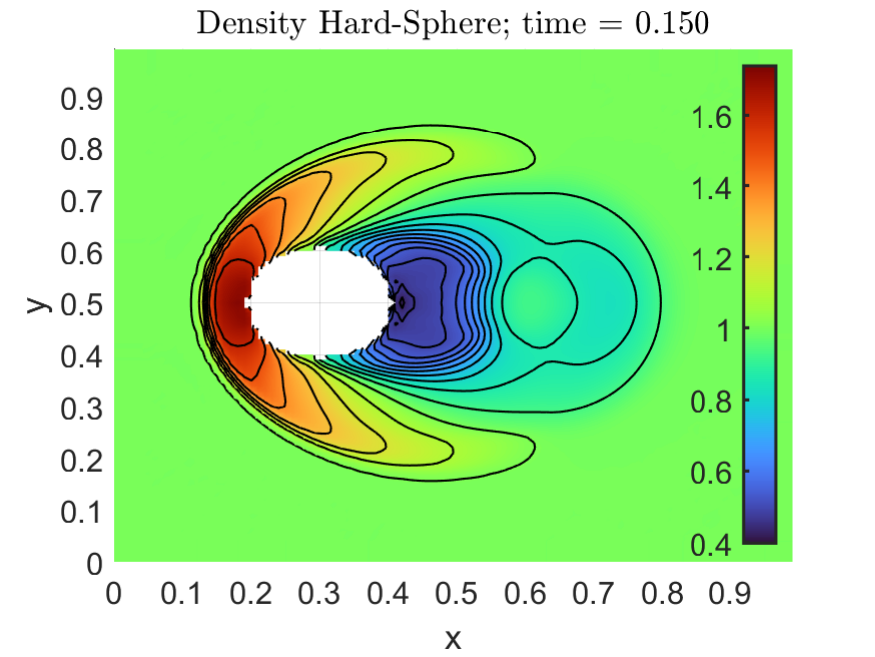}
    \end{subfigure}

    \caption{Time evolution of density for \hyperref[subTestFluidObstacle]{\textbf{Test 4}}. In the first row is displayed the regime adaptation at three different times. Red, green and black dots, identify, respectively, cells updated using Euler equations, ES-BGK equation and Boltzmann equation. In the second, third and fourth row are displayed, respectively, the time evolution of density computed using hybrid scheme, Euler equations and Boltzmann equation.
    %\textbf{\textcolor{red}{WRONG CAPTION}}
    }
    \label{fig_density_sim4}
\end{figure}

\begin{figure}
    \begin{subfigure}[b]{0.32\textwidth}
        \centering
        \includegraphics[width=\textwidth, trim={0.9cm 0cm 0.9cm 0cm}]{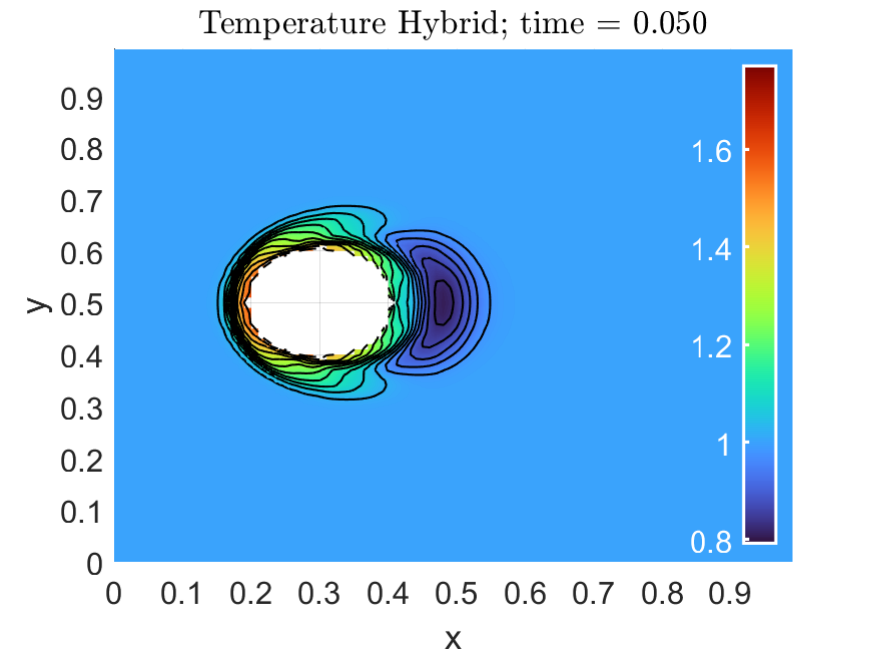}
    \end{subfigure}
    \hfill
    \begin{subfigure}[b]{0.32\textwidth}
        \centering
        \includegraphics[width=\textwidth, trim={0.9cm 0cm 0.9cm 0cm}]{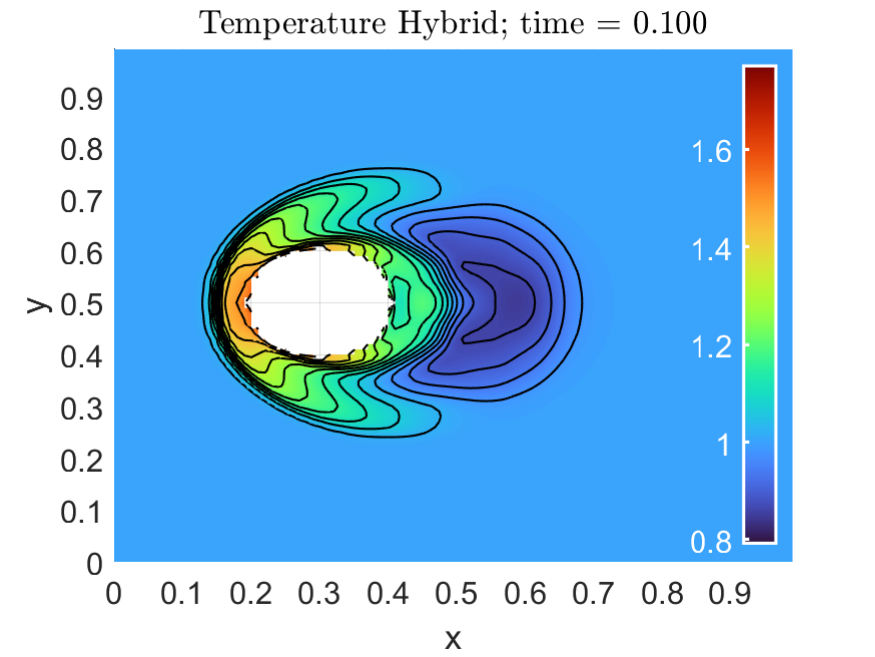}
    \end{subfigure}
    \hfill
    \begin{subfigure}[b]{0.32\textwidth}
        \centering
        \includegraphics[width=\textwidth, trim={0.9cm 0cm 0.9cm 0cm}]{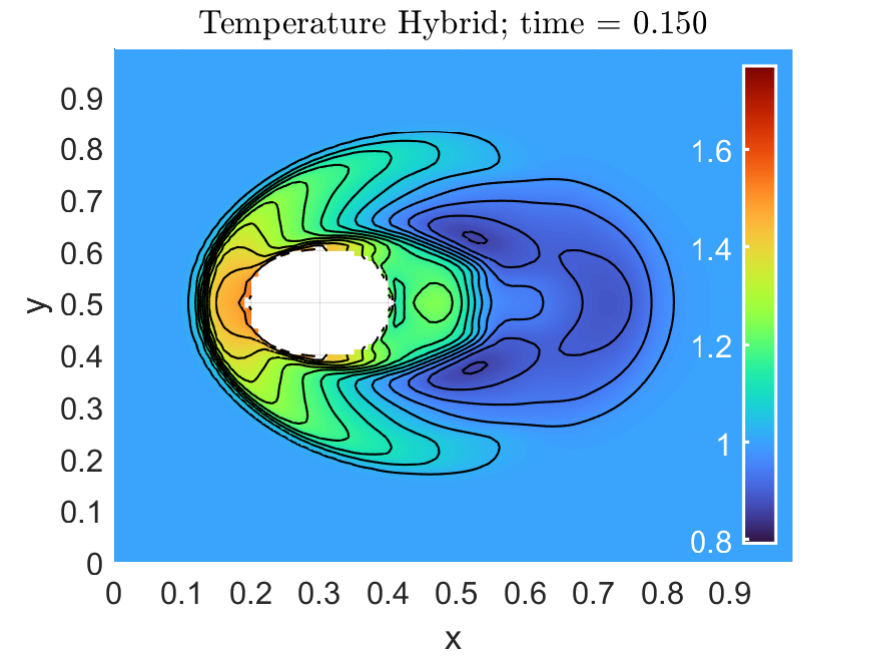}
    \end{subfigure}

    \begin{subfigure}[b]{0.32\textwidth}
        \centering
        \includegraphics[width=\textwidth, trim={0.9cm 0cm 0.9cm 0cm}]{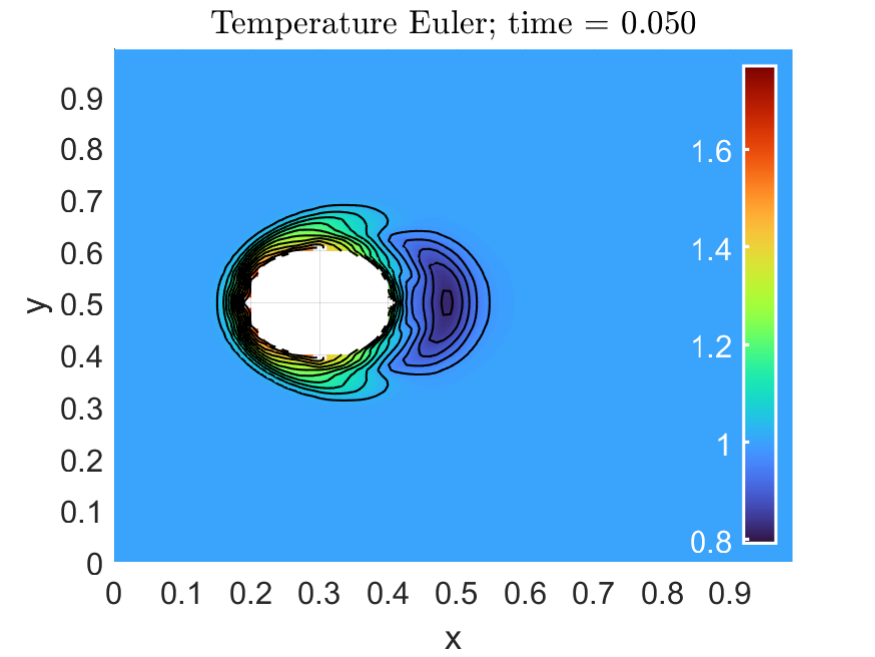}
    \end{subfigure}
    \hfill
    \begin{subfigure}[b]{0.32\textwidth}
        \centering
        \includegraphics[width=\textwidth, trim={0.9cm 0cm 0.9cm 0cm}]{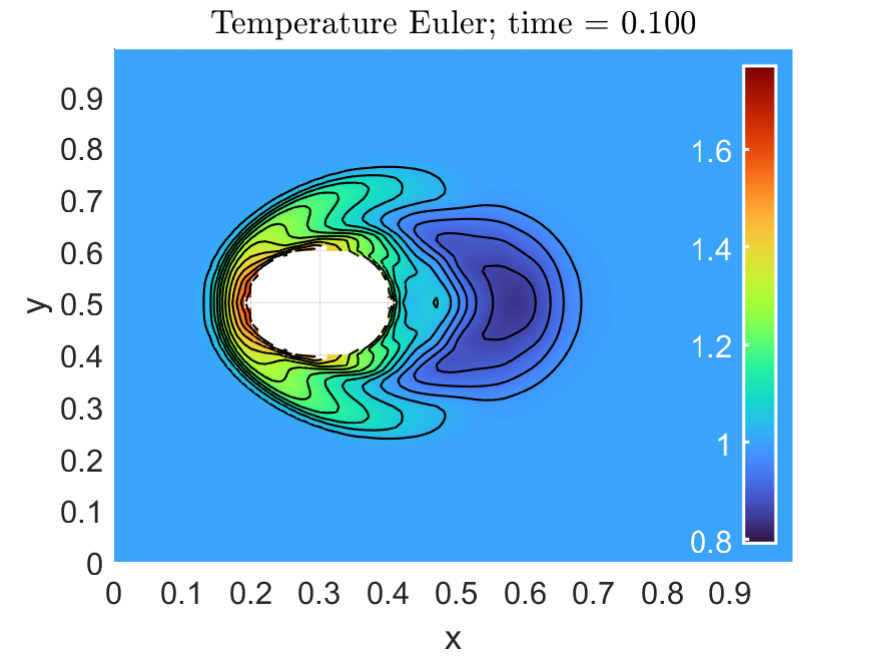}
    \end{subfigure}
    \hfill
    \begin{subfigure}[b]{0.32\textwidth}
        \centering
        \includegraphics[width=\textwidth, trim={0.9cm 0cm 0.9cm 0cm}]{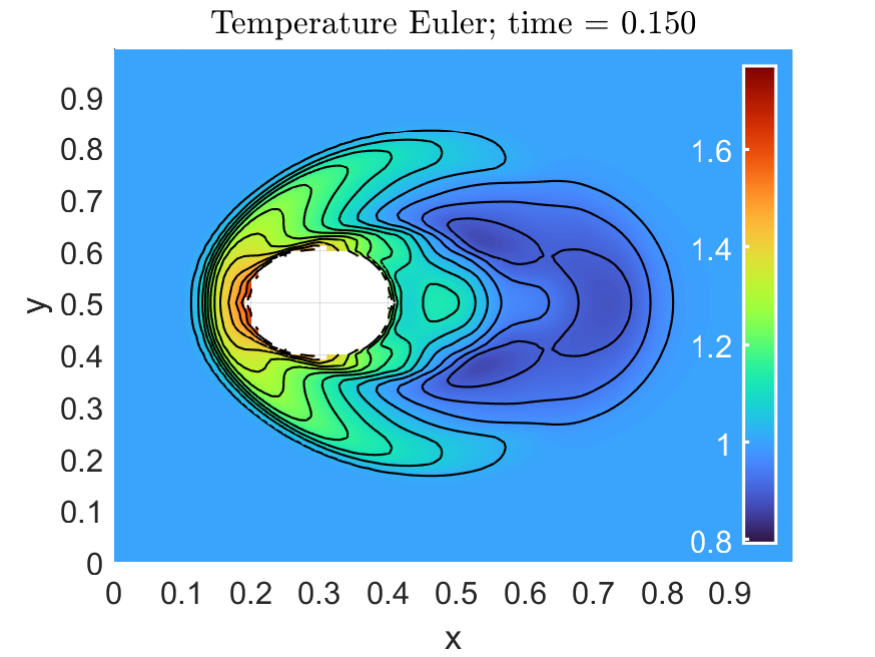}
    \end{subfigure}

    \begin{subfigure}[b]{0.32\textwidth}
        \centering
        \includegraphics[width=\textwidth, trim={0.9cm 0cm 0.9cm 0cm}]{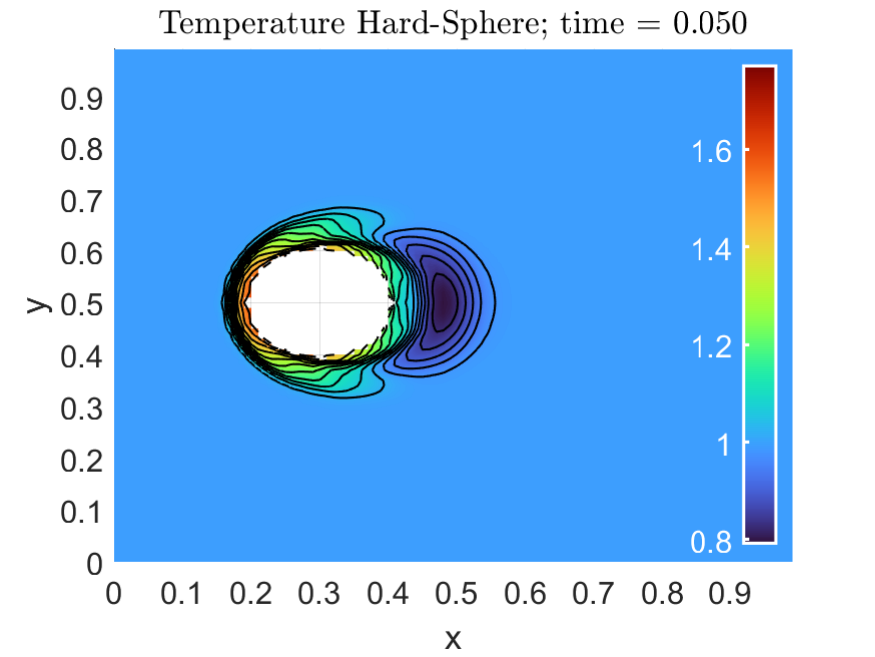}
    \end{subfigure}
    \hfill
    \begin{subfigure}[b]{0.32\textwidth}
        \centering
        \includegraphics[width=\textwidth, trim={0.9cm 0cm 0.9cm 0cm}]{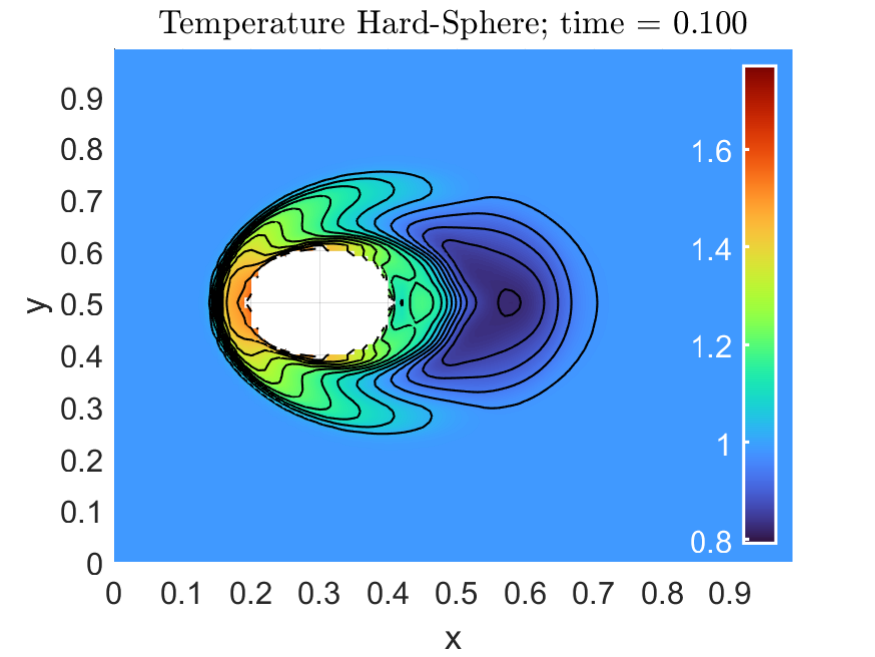}
    \end{subfigure}
    \hfill
    \begin{subfigure}[b]{0.32\textwidth}
        \centering
        \includegraphics[width=\textwidth, trim={0.9cm 0cm 0.9cm 0cm}]{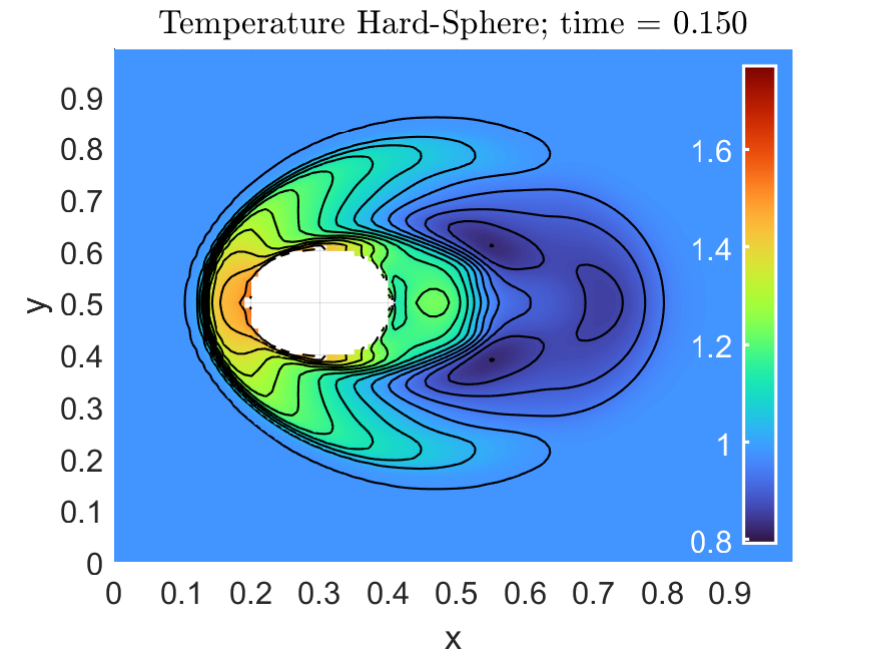}
    \end{subfigure}

    \caption{Time evolution of temperature for \hyperref[subTestFluidObstacle]{\textbf{Test 4}}. In the first, second and third row  are displayed, respectively, the time evolution of temperature computed using hybrid scheme, Euler equations and Boltzmann equation.} %In the plot of the temperature for the HardSphere simulation it seems that there is like a bump in the front region if we look at the contour plot. However I think it is due to the fact that I am using the same contour lines for all simulations. In fact if I set the countour levels adapted to the maximum and the minimum of just the hard sphere, the contour lines will not show this bump
    \label{fig_temperature_sim4}
\end{figure}

To simulate the presence of the object, we have defined a circular domain (white cells in Figure \ref{fig_density_sim4} and Figure \ref{fig_temperature_sim4}) in which the integrators are never executed. These cells, are part of the stencils of the closest cells of the integration domain, thus we have assigned to them a constant and uniform value of pressure, density, temperature and velocity, given by: $P_*=1.5, \rho_*=1.0, T_*=P_* / \rho_*$ and $\vec{u}_*=\vec{0}$. This circular shape domain, is surrounded by a small region (2 cells thick) in which the integrator is always the one associated to the full Boltzmann equation with Boltzmann operator. We identify by $(x_c, y_c)=(3 N_x/10, N_y/2)$ the coordinates of the centre of the round shape, and we identify all the cells satisfying the following condition $\left(x_i-x_c\right)^2 + \left(y_i-y_c\right)^2\leq R_c^2$ as part of the object (white cells), where $R_c=N_x/10$ represents the \emph{radius} of the round obstacle.

The fluid in the integration domain is initialized with the following values of pressure, density, temperature and velocity: $P=1.0, \rho=1.0, T=P / \rho$ and $\vec{u}_*=(1.0, 0.0, 0.0)$. The distribution functions are initialized everywhere as a Maxwellian whose moments correspond to the initial values of the Euler equations.
% The geometric schematization is reported in Figure \ref{fig7}, and we use the same symbol convention as the one used in Figure \ref{fig4}. The dimensions of the object in the Figure \ref{fig7} do not correspond to those actually implemented in the code, but the correct ones are reported below: \\
% In the dashed (blue) regions of this last Figure, Algorithm \ref{alg1} is executed at each timestep, and then the respective regime integrator is used to compute the solution. The white cells constitute the fixed object, and numerical integration is never performed in them. Since these cells constitute the stencil cells of the green rhomboid cells, it is necessary to assign to them a temperature, pressure, density and velocity value, which are given respectively by $P_*=2.5, \rho_*=1.0, T_*=P_* / \rho_*$ and $\vec{u}_*=\vec{0}$. Inside the green cells the hard-sphere operator is always integrated. In all other (blue dashed) cells, Algorithm \ref{alg1} is executed. The red arrows schematize the direction of the fluid flow. The yellow dashed circle represents the ideal round shape of the obstacle. Of course, since we are dealing with a cartesian grid, the round shape can only be sufficiently approximated. In particular, increasing the value of $R_c$ then the cartesian approximation of the round shape will become better and better.

 In Figure \ref{fig_density_sim4} we show the density profile with $\epsilon = 10^{-4}$, at time $t=0.050, 0.100, 0.150$ for the hybrid scheme, for the  Euler equations, and for the Boltzmann operator (this last one computed using CUDA GPU parallelization), as well as the domain indicators for the hybrid scheme (red, green and black dots indicate cells computed using, respectively, Euler, ES-BGK and Boltzmann equations). It is important to note that some kinetic cells (green or black), after a certain amount of time, become hydrodynamic (red) again as the fluid moves toward equilibrium (Maxwellian distribution). In Figure \ref{fig_temperature_sim4} is reported the temperature profile for the same instants in time, for the hybrid scheme, the Euler and Boltzmann equations. We estimated, as explained in Section \ref{subNumSim}, the total duration of the numerical integration performed using just the full Boltzmann operator on CPU over the entire domain, and the ratio between the execution time for the full Boltzmann and the hybrid scheme (Euler equations, ES-BGK operator and Boltzmann operator) is equal to 2.6.
%2.4.

The thresholds used are: $\eta_0=1\times 10^{-3}$, $\eta_1=2.8\times10^{-5}$, $\delta_0=10^{-4}$. The time step is $0.1 dx$.

% \begin{figure}
%     \centering
%     \includegraphics[width=0.9\linewidth, trim={0cm 0cm 0cm 0cm}]{Figure8 - Rotate.pdf}
%     \caption{Schematic view of the numerical setup for \hyperref[subTestFluidObstacle]{\textbf{Test 4}}, where the symbol conventions used are the ones of Figure \ref{fig4}. We identify by $(x_c, y_c)=(3 N_x/10, N_y/2)$ the coordinates of the centre of the round shape, and we identify all the cells satisfying the following condition $\left(x_i-x_c\right)^2 + \left(y_i-y_c\right)^2\leq R_c^2$ as part of the object (white cells), where $R_c$ represents the \emph{radius} of the round obstacle. In the dashed (blue) regions, Algorithm \ref{alg1} is executed at each timestep. The white cells constitute the fixed object, where temperature, pressure, density and velocity values are constant over time. Inside the green cells the Boltzmann equation with Boltzmann operator is always integrated. The red arrows schematize the direction of the fluid flow. The yellow dashed circle represents the ideal round shape of the obstacle.}
%     \label{fig7}
% \end{figure}

\subsection{Test 5: Non uniform Knudsen number}
\label{subNonUniformKnudsen}
\begin{figure}
    \centering
    \begin{subfigure}[b]{0.32\textwidth}
        \centering
        \includegraphics[width=\textwidth, trim={0.9cm 0cm 0.9cm 0cm}]{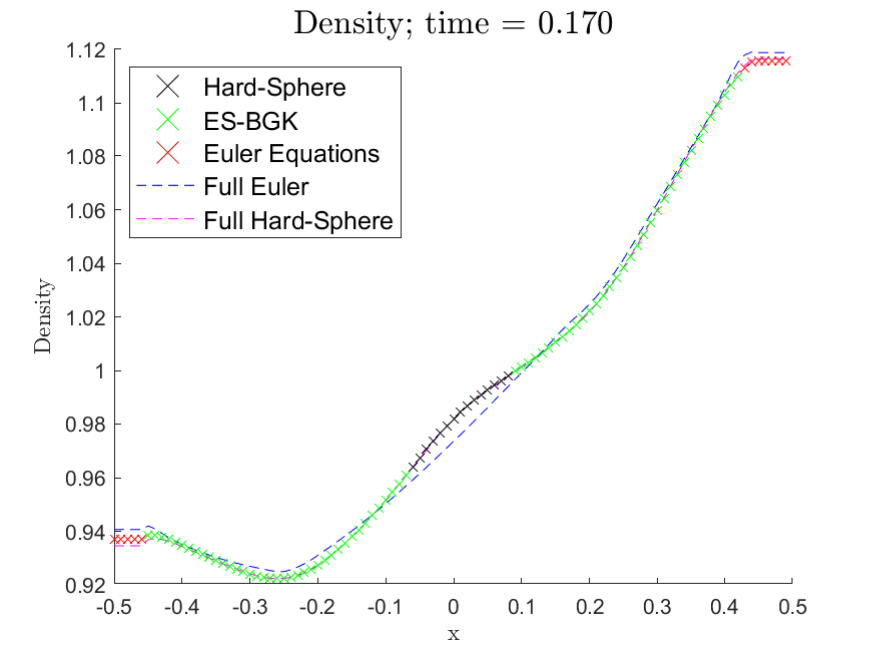}
    \end{subfigure}
    \hfill
    \begin{subfigure}[b]{0.32\textwidth}
        \centering
        \includegraphics[width=\textwidth, trim={0.9cm 0cm 0.9cm 0cm}]{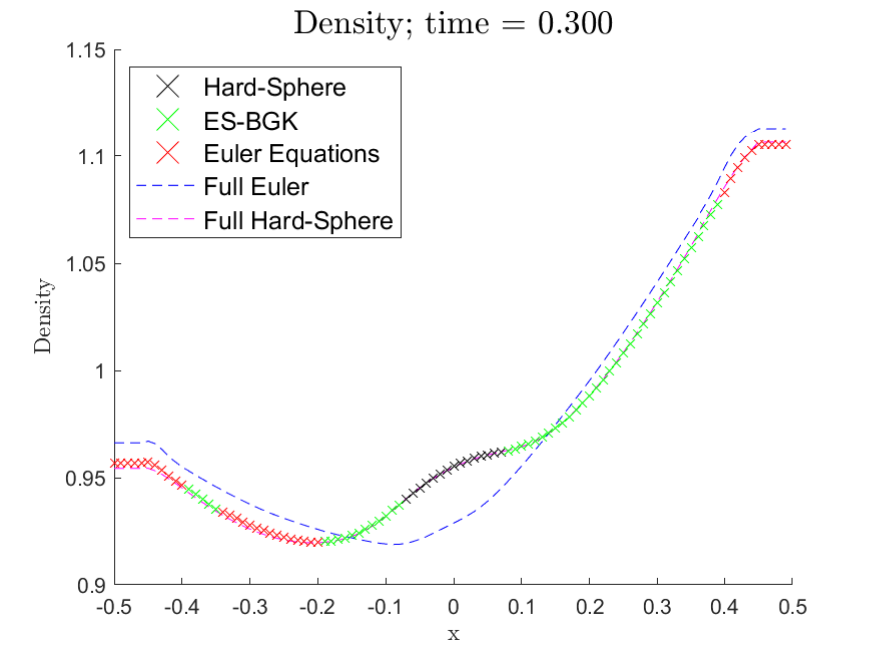}
    \end{subfigure}
    \hfill
    \begin{subfigure}[b]{0.32\textwidth}
        \centering
        \includegraphics[width=\textwidth, trim={0.9cm 0cm 0.9cm 0cm}]{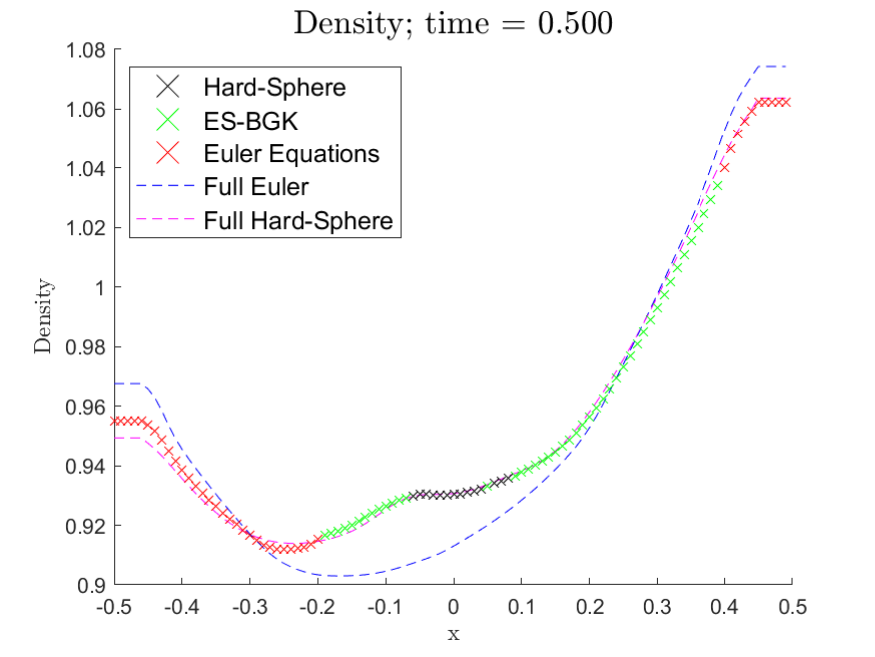}
    \end{subfigure}

    \begin{subfigure}[b]{0.32\textwidth}
        \centering
        \includegraphics[width=\textwidth, trim={0.9cm 0cm 0.9cm 0cm}]{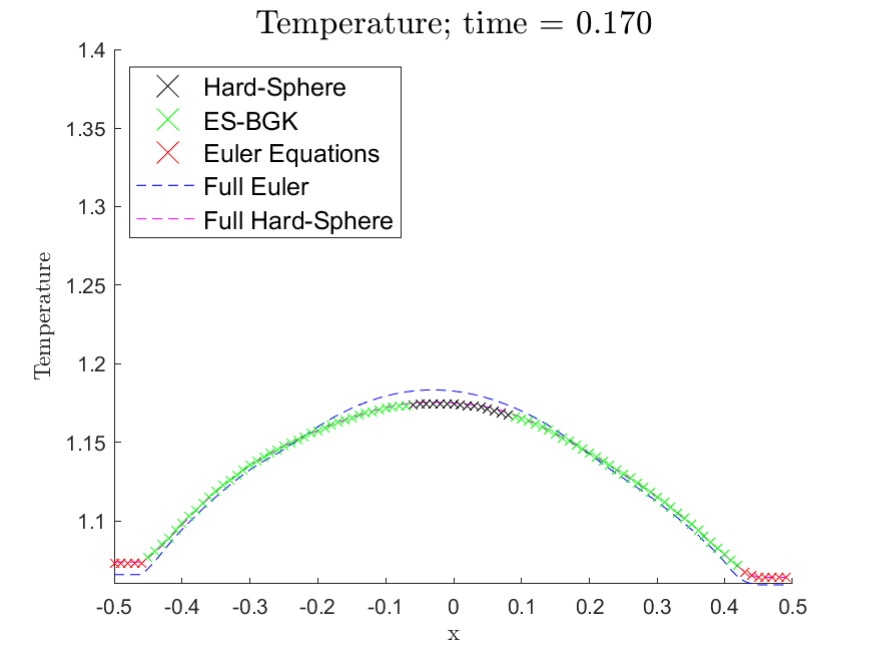}
    \end{subfigure}
    \hfill
    \begin{subfigure}[b]{0.32\textwidth}
        \centering
        \includegraphics[width=\textwidth, trim={0.9cm 0cm 0.9cm 0cm}]{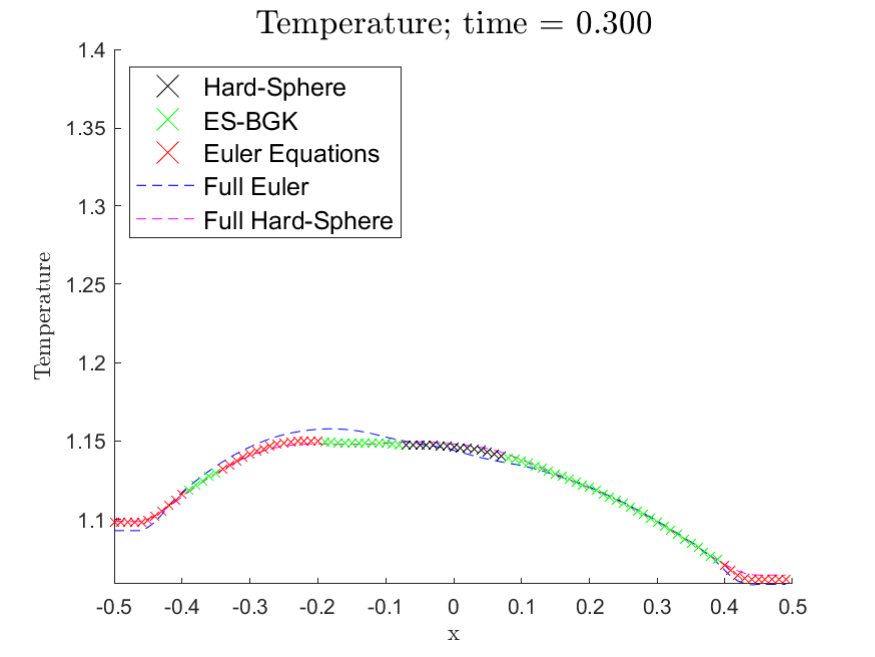}
    \end{subfigure}
    \hfill
    \begin{subfigure}[b]{0.32\textwidth}
        \centering
        \includegraphics[width=\textwidth, trim={0.9cm 0cm 0.9cm 0cm}]{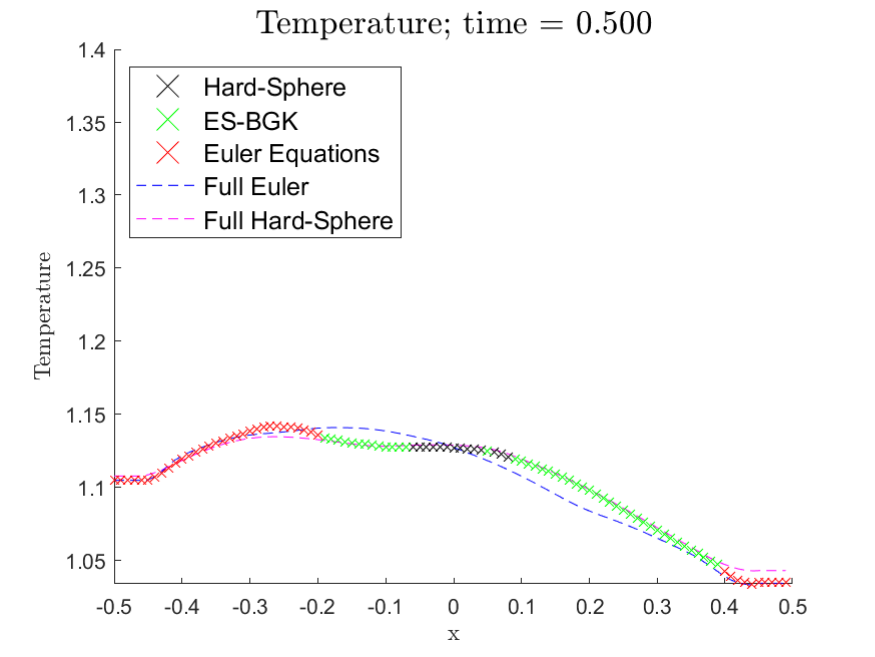}
    \end{subfigure}

    \caption{Time evolution of density for \hyperref[subNonUniformKnudsen]{\textbf{Test 5}}. In the first row the regime adaptation and the density at three different times are displayed. Red, green and black dots, identify, respectively, cells updated using Euler equations, Boltzmann equations with ES-BGK operator and Boltzmann equations with Boltzmann (Hard-Sphere) operator, using the new hybrid scheme. The magenta line is the solution obtained using only the Euler equations, and the blue line is computed using only the Boltzmann equation. The same plots, but for the temperature, are showed in the second row. It is clearly visible that the hybrid solution matches almost perfectly the one obtained with the hard-sphere solver.}
    \label{fig_density_sim5}
\end{figure}

The aim of this test is to check if the scheme is able to deal with a spatially non-uniform Knudsen number. This simulation is inspired by one of those presented in \cite{fil-rey-2014}. Here we consider a 2D spatial domain discretized with $N_x=50$ and $N_y=16$ points, with a physical space $[-L_x, L_x]$ and $[-L_y, L_y]$ of dimensions $L_x=1/2$ and $L_y$ chosen such that $dy=dx$. The macroscopic initial conditions are the following:
\begin{equation}
    \begin{cases}        
        \rho(x,y)=1+0.1 \sin{(\pi x)},\\
        u_x (x,y)= 0.125,\\
        u_y (x,y)=0,\\
        %T(x,y) = (5+2\cos{(2\pi x)})/20,
        T(x,y) = 1+0.2\cos{(\pi x}),
    \end{cases}    
\end{equation}
and the distribution function is initialized as 
$$f(0, x, \cdot) = \frac{1}{2}\left(\mathcal{M}_{\rho(x), \bm{u(x)}, T(x)}+\mathcal{M}_{\rho(x),-\bm{u(x)}, T(x)}\right),$$
so it is everywhere far from equilibrium.\\
The Knudsen number is a function of the physical space accordingly to the following smooth relation $\epsilon(x)=10^{-6}+0.5*\left(\arctan{(1+30x)}+\arctan{(1-30x)}\right)$, which means that the more external region of the domain is where we expect the Euler regime to be triggered by the hybrid scheme due to the very small Knudsen number. The regime is initially set as ES-BGK everywhere, and then we let it evolves accordingly to our scheme. The velocity domain $[-L,L]^3$, where $L=8$ is discretized with 16 points in each direction. The value of $R$ has been taken approximately equal to 3.3, and both numbers of points $A_1$ and $A_2$ discretizing the spherical integration of the Boltzmann operator (see Section \ref{subBoltzOp}) have been taken equal to 4. \\
In Figure \ref{fig_density_sim5} it is reported the density and temperature evolution computed with the hybrid scheme (red, green and black crosses are, respectively, the cells computed using Euler equations, ES-BGK model and Hard-Sphere Boltzmann model), using only Euler equations (blue dashed line) and using only Full Boltzmann Hard-Sphere model (magenta dashed line). It is evident that the Euler solution is far from the reference one (Boltzmann Hard-Sphere model) since the Knudsen number is not everywhere sufficiently small. The hybrid scheme correctly triggers the Euler regime only where the Knudsen number is small (close to the extrema of the physical space), and the region associated to the Boltzmann model (black crosses) is well confined in space in the central region. In addition the solution computed with the hybrid scheme is in perfect agreement with the reference solution (magenta line) which is a clear indicator of the correctness of the implemented methods and of the switching procedure between different regimes.

The thresholds used are: $\eta_0=4\times 10^{-2}$, $\eta_1=8\times10^{-2}$, $\delta_0=10^{-3}$. The time step is $0.1 dx$. The ratio between the execution time for the full Boltzmann and the hybrid scheme (Euler equations, ES-BGK operator and Boltzmann operator) is equal to 1.3, which means that the full Boltzmann solver is 30\% slower than the hybrid scheme.\\

\subsection{Discussion on speed-up}
As already observed, the speed-up offered by our numerical domain adaptation compared to the simulations computed using only the full hard-sphere Boltzmann solver is almost equal to 3. This is a quantity which depended on several factors, among which, the principals are, the choice of the thresholds which directly determines the sizes of the domains and the actual numerical implementation of the code.\\
In Table \ref{table:1} we report, for each simulation, the total percentage of cells used by each domain computed over the total number of cells of the numerical integration (i.e. 100\% corresponds to the total number of grid cell in a single timestep multiplied by the total number of timesteps). In the last column is indicated the ratio of time execution of a given solver by the time execution of the Euler solver. These execution times are computed for the \hyperref[subTestFluidFixedRectObs]{\textbf{Test 2}}, and using each solver over the entire spatial domain for the entire duration of the simulation. Of course, the percentage values depends on the thresholds implemented, and the ratio between the execution time strongly depend on the numerical implementation of the code and the machine performances.
\begin{center}
    \begin{table}
        \begin{tabular}{||c || c | c | c | c | c|| c ||}
        \hline
                                    & \multicolumn{5}{c||}{Total percentage} & \\
                                    & TEST 1   & TEST 2   & TEST 3   & TEST 4   & TEST 5 & Ratio Execution Time\\
        \hline
            Euler       & 82.5\%   & 94.1\%   & 92.6\%   & 93.0\%   & 25.1\% & 1\\
            ES-BGK      & 6.7\%    & 3.4\%    & 6.6\%    & 4.8\%    & 61.8\% & 4\\
            Hard-Sphere & 10.8\%   & 2.5\%    & 0.8\%    & 2.2\%    & 13.1\% & 14\\
        \hline    
        %Se calcolassi il Ratio Execution Time usando il TEST5, allora i valori numerici sarebbero 1, 7, 15
        
        \end{tabular}
        \caption{It is indicated for each simulation the total percentage of cells used by each domain computed over the total number of cells of the numerical integration (i.e. 100\% corresponds to the total number of grid cell in a single timestep multiplied by the total number of timesteps). In the last column is indicated the ratio of time execution of a given solver by the time execution of the Euler solver. These execution time are computed for the \hyperref[subTestFluidFixedRectObs]{\textbf{Test 2}}, and using each solver over the entire spatial domain for the entire duration of the simulation.}
    \end{table}
    \label{table:1}
\end{center}
	
\section{Conclusions}
\label{secCon}
In this work, we proposed a hierarchical domain decomposition strategy for solving the Boltzmann equation—and more generally, kinetic models—in multiscale regimes. The method is guided by two main criteria. The first, inspired by the works of Levermore, Morokoff, Nadiga \cite{lev-mor-nad-1998}, Filbet and Rey \cite{fil-rey-2014}, Filbet and Xiong \cite{fil-xio-2018}, Tiwari \cite{tiw-2000}, Li, Son and Wang \cite{li-son-wan-2021} and Xiong and Qiu \cite{xio-qiu-2017}, uses a Chapman–Enskog expansion to determine when to transition from a macroscopic fluid model to either the ES-BGK approximation or the full Boltzmann equation. This criterion relies solely on macroscopic quantities derived from the closure of the kinetic model and does not require evaluating the distribution function itself.
The second criterion governs the transition from the kinetic regime back to its hydrodynamic limit and is based on comparing the truncated Chapman–Enskog expansion with the corresponding equilibrium state.

A key novelty of this work is the implementation of a spatial 2D and velocity 3D hybrid fluid–kinetic solver involving three distinct physical regimes. The method combines high-order asymptotic-preserving (AP) schemes and fast spectral solvers for the Boltzmann collision operator, ensuring both stability and accuracy across a wide range of Knudsen numbers. The proposed numerical experiments demonstrate the robustness and effectiveness of our switching criteria. Notably, the kinetic regions are automatically confined both in space and time according to the evolving solution, and the hybrid scheme achieves speedups of nearly a factor of three compared to full Boltzmann solvers. Additional performance gains are expected through further optimization of the implementation, including parallelization and memory management.

Future work will focus on extending the method to multi-species gas mixtures and plasma applications, as well as enhancing computational efficiency through advanced software engineering and high-performance computing strategies.

\appendix
\section{Moment realizability matrix}\label{appendix:mrm}
    In this Appendix we give more details about the Moment Realizability Matrix and in particular on the derivation of the Equation \eqref{eqMomRealizReduced}.
    The matrix $\bm{\bar{A}}^\ve$ and the vector $\bm{\bar{B}}^\ve$ will allow us to define our hydrodynamic break down criterion.  Let us set the vector of the reduced collisional invariants for $\bm{V}=(v-\bm{u^\ve})/\sqrt{T^\ve}$,
            \[
            \bm{m} := \left (1, \bm{V}, \left (\frac{2}{5}\right )^{1/2}\left (\frac{|\bm{V}|^2}{2} - \frac{5}{2}\right ) \right ).
            \] 
            Similarly to \cite{lev-mor-nad-1998, tiw-2000}, now we define the so-called \emph{moment realizability matrix} $\bm{M}$ by setting
            \begin{equation}
                \label{eqMomRealiz_appendix}
                \bm{M} := \frac{1}{\rho^\ve}\int_{\RR^{3}} \left(\bm{m} \otimes \bm{m} \right) f^\ve(v) \, dv.
            \end{equation}
            Indeed, if the distribution function is nonnegative, then the quantity $\int \phi^2(v) f \, dv $ is also nonnegative \cite{tiw-2000}, for any function $\phi(v)$. The same is true if we consider $\bm{m}=\bm{m}(\bm{v})$ a columns vector of $c$ polynomials, and we set $\bm{\phi}=\bm{a}^\intercal \bm{m}$ for a generic $\bm{a}\in R^c$, and we have
            \begin{equation}
                \bm{a}^\intercal\left(\int_{\RR^{3}}\bm{m}\otimes \bm{m} f \, dv \right)\bm{a}=\int_{\RR^{3}} (\bm{a}^\intercal \bm{m})^2f\, dv \geq 0.
            \end{equation}
            Then if $f$ is nonnegative and not identically zero, then  for every $\bm{m}$ the $c\times c$ matrix
            \begin{equation}\label{eq1005_appendix}
                \bm{M}=\frac{1}{\rho}\int_{\RR^{3}} \bm{m}\otimes \bm{m} f\, dv,
            \end{equation}
            is positive semi-definite.\\
                        
            Then, the quantity given by Equation \refeq{eqMomRealiz_appendix} is also positive semi-definite. After some computations we obtain

            \begin{align*}
                \bm{M} =
                \begin{pmatrix}
                1                  & \bm{0}_{\RR^{3}}^\intercal                 & -(\frac{2}{5})^{\frac{1}{2}} \\
                \bm{0}_{\RR^{3}}   & \bm{I} +\bm{\bar{A}}^\ve                     & \left (\frac{2}{5}\right )^{1/2} \bm{\bar{B}}^\ve \\
                -(\frac{2}{5})^{\frac{1}{2}}                  & \left (\frac{2}{5}\right )^{1/2} (\bm{\bar{B}}^\ve)^\intercal & \bar{C}^\ve
                \end{pmatrix},   
            %\label{eqMomRealizFull}
            \end{align*}
            where $\bar{C}^\ve$ is the dimensionless fourth order moment of $f^\ve$
            \[
             \bar{C}^\ve := \frac{2}{5 \rho}\int_{\RR^{3}} \left [\frac{|\bm{V}|^2}{2} - \frac{5}{2}\right ]^2 f^\ve(v) \, dv.
            \]
            The matrix $\bm{M}$ is of the form
            \begin{equation}
                \bm{M}=
                \begin{pmatrix}
                    \bm{M_A}  & \bm{M_B}\\
                    \bm{M_B}^\intercal & M_C\\
                \end{pmatrix},
            \end{equation}
            where
            \begin{equation}
                \bm{M_A}=
                \begin{pmatrix}
                    1 &  \bm{0}_{\RR^{3}}^\intercal\\
                    \bm{0}_{\RR^{3}}^\intercal & \bm{I} +\bm{\bar{A}}^\ve\\
                \end{pmatrix},\, \, \, 
                \bm{M_B}=
                \begin{pmatrix}
                    -(\frac{2}{5})^{\frac{1}{2}}\\
                    \left (\frac{2}{5}\right )^{1/2} \bm{\bar{B}}^\ve\\
                \end{pmatrix},\, \, \, 
                M_C = \bar{C}^\ve.
            \end{equation}
            Assuming that the scalar quantity $M_C =\bar{C}^\ve>0$ (we will verify later that this assumption is verified in our circumstances) then the matrix $\bm{M}$ is positive semi-definite if and only if its Schur complement $X$ is positive semi-definite
            \begin{equation}
                \displaystyle \bm{X} = \displaystyle \bm{M_A} - \bm{M_B} M_C^{-1}\bm{M_B}^\intercal =
                \displaystyle
                \begin{pmatrix}
                    1-\frac{2}{5\bar{C}^\ve}& \frac{2}{5} \frac{\left(\bm{\bar{B}}^\ve\right)^\intercal}{\bar{C}^\ve} \\
                     \frac{2}{5} \frac{\bm{\bar{B}}^\ve}{\bar{C}^\ve} & \bm{I}+\bm{\bar{A}}^\ve - \frac{2}{5}\frac{\bm{\bar{B}}^\ve\left(\bm{\bar{B}}^\ve\right)^\intercal}{\bar{C}^\ve}
                \end{pmatrix}.
            \end{equation}
            For the sake of simplicity, let's introduce the following matrix $\bm{Q}$ for the change of basis
            \begin{equation}
                \bm{Q} = 
                \begin{pmatrix}
                    \alpha & -\frac{\chi_1^\intercal}{\chi_2} \\
                    \bm{0}_{\RR^3}      & \bm{I}
                \end{pmatrix},
            \end{equation}
            where $\bm{\chi}_1=\frac{2}{5} \frac{\bm{\bar{B}}^\ve}{\bar{C}^\ve}$, $\chi_2 = 1-\frac{2}{5\bar{C}^\ve}$ and $\alpha^2 = \frac{\bar{C}^\ve}{1 - \frac{2}{5\bar{C}^\ve}}$. Then the operation $\bm{Q}^\intercal \bm{X} \bm{Q}:=\bm{W}$ corresponds to a change of coordinate (since $X$ can be viewed as quadratic form) and this operation does not change the semipositiveness of the matrix $\bm{X}$, i.e. the signs of the eigenvalues of the matrices $\bm{X}$ and $\bm{W}$ are the same. So we obtain that
            \begin{equation}\label{eq1004_appendix}
                \bm{W} := \bm{Q}^\intercal \bm{X} \bm{Q} = 
                \begin{pmatrix}
                    \bar{C}^\ve & \bm{0}_{\RR^3}^\intercal\\
                    \bm{0}_{\RR^3} & \bm{I}+\bm{\bar{A}}^\ve - \frac{2}{5}\frac{\bm{\bar{B}}^\ve\left(\bm{\bar{B}}^\ve\right)^\intercal}{\bar{C}^\ve} - \frac{\bm{\chi}_1 \bm{\chi}_1^\intercal}{\chi_2}
                \end{pmatrix}.
            \end{equation}
            It is important to note that doing that we have assumed the validity of the following relation $\alpha^2 = \frac{\bar{C}^\ve}{1 - \frac{2}{5\bar{C}^\ve}}$, which requires that $\bar{C}^\ve>\frac{2}{5}$. Now we want to prove that, in our framework, $\bar{C}^\ve=1>\frac{2}{5}$. If the distribution function is at equilibrium, i.e. $f=\mathcal{M}_{\rho,\bm{u}, T}$ then a direct computation gives that $\bar{C}^\ve=1$. More in general, if we assume that the distribution function $f$ has a small deviation from the local Maxwellian $\mathcal{M}_{\rho,\bm{u}, T}$ then we can approximate it by \cite{tiw-2000}
            \begin{equation}
                f = \mathcal{M}_{\rho,\bm{u}, T}(1+\phi),
            \end{equation}
            where
            \begin{equation}
                \phi=\frac{<\bm{v}-\bm{u}, \bm{q}>}{\rho T^2}\left[\frac{|\bm{v}-\bm{u}|^2}{5T}-1\right]+
                \frac{1}{2\rho T^2}(\bm{v}-\bm{u})^\intercal \bm{\tau} (\bm{v}-\bm{u}), 
            \end{equation}
            where $\bm{\tau}$ and $\bm{q}$ are, respectively, the stress tensor and the heat flux
            \begin{align}
                \tau_{ij}&=\int_{\RR^3}(v_i-u_i)(v_j-u_j)f(t,x,v)dv - p \delta_{ij},\\
                q_j&=\frac{1}{2}\int_{\RR^3}(v_j-u_j)|v-u|^2 f(t,x,v) dv.
            \end{align}
            Using this last expression for $\phi$, it can be proved \cite{tiw-2000} that if $f=\mathcal{M}_{\rho,\bm{u}, T}(1+\phi)$ then $\bar{C}^\ve=1$. This confirms that in our framework $\bar{C}^\ve=1>\frac{2}{5}$.\\

            Since $\bar{C}^\ve=1$, then from Equation \eqref{eq1004_appendix} the matrix $\bm{W}$ is positive semi-definite, if and only if, the following quantity is semipositive definite
            \begin{equation}
                \mathcal{V}_{\ve}:=\bm{I}+\bm{\bar{A}}^\ve - \frac{2}{5}\frac{\bm{\bar{B}}^\ve\left(\bm{\bar{B}}^\ve\right)^\intercal}{\bar{C}^\ve} - \frac{\bm{\chi}_1 \bm{\chi}_1^\intercal}{\chi_2}=
                \bm{I} + \bm{\bar{A}}^\ve - \frac{\frac{2}{5\bar{C}^\ve}}{1-\frac{2}{5\bar{C}^\ve}}\bm{\bar{B}}^\ve\otimes \bm{\bar{B}}^\ve.
            \end{equation}
            Substituting $\bar{C}^\ve=1$ in this last equation we get
            \begin{equation}\label{eqMomRealizReduced_appendix}
                \mathcal{V}_{\ve} = \bm{I} + \bm{\bar{A}}^\ve - \frac{2}{3}\bm{\bar{B}}^\ve\otimes \bm{\bar{B}}^\ve.
            \end{equation}
            From this derivation, the initial matrix $\bm{M}$ defined in Equation \eqref{eq1005_appendix}, is positive semi-definite if and only if the matrix $\mathcal{V}_{\ve}$ is positive semi-definite. Then the idea is to monitor the evolution of the matrix $\mathcal{V}_{\ve}$ \cite{tiw-1998, lev-mor-nad-1998}: indeed the moment realizability criterion states that the fluid dynamic description has broken down when the perturbation is too large that $\mathcal{V}_{\ve}$, and thus $\bm{M}$, is no longer positive semi-definite.
 
	\section*{Acknowledgments}
    The authors received funding from the European Union's Horizon Europe research and innovation program under the Marie Skłodowska-Curie Doctoral Network DataHyking (Grant No. 101072546). 
    DC and TR are also supported by the French government, through the UniCA$_{JEDI}$ Investments in the Future project managed by the National Research Agency (ANR) with the reference number ANR-15-IDEX-01.	 
    DC would like to thanks Tommaso Tenna for the useful discussions about the manuscript. The research of LP has been supported by the Royal Society under the Wolfson Fellowship “Uncertainty quantification, data-driven simulations and learning of multiscale complex systems governed by PDEs”. This work has been written within the activities of GNCS group of INdAM (Italian National Institute of High Mathematics). LP also acknowledges the  partial support 
by European Union -
NextGenerationEU through the Italian Ministry of University and Research as
part of the PNRR – Mission 4 Component 2, Investment 1.3 (MUR Directorial
Decree no. 341 of 03/15/2022), FAIR “Future” Partnership Artificial Intelligence Research”, Proposal Code PE00000013 - CUP DJ33C22002830006) and by MIUR-PRIN Project 2022, No. 2022KKJP4X “Advanced numerical methods for time dependent parametric partial differential equations with applications”.

  \bibliographystyle{acm}
  \bibliography{biblio}

\end{document}